\documentclass{amsart}
\usepackage{amsmath}
\usepackage{amssymb}
\usepackage{amsthm}

\title{On the tensor product of two basic representations
of $U_v(\hat{sl}_e)$}

\author{Susumu Ariki}
\address{Research Institute for Mathematical Sciences, %
         Kyoto University, Kyoto 606-8502. Japan.}
\email{ariki@kurims.kyoto-u.ac.jp}

\author{Victor Kreiman}
\address{Department of Mathematics, %
University of Georgia, Athens, GA 30602, USA}
\email{vkreiman@math.uga.edu}

\author{Shunsuke Tsuchioka}
\address{Research Institute for Mathematical Sciences, %
         Kyoto University, Kyoto 606-8502. Japan.}
\email{tshun@kurims.kyoto-u.ac.jp}

\date{}

\newtheorem{thm}{Theorem}[section]
\newtheorem{defn}[thm]{Definition}
\newtheorem{prop}[thm]{Proposition}
\newtheorem{lemma}[thm]{Lemma}
\newtheorem{cor}[thm]{Corollary}
\newtheorem{Ex}[thm]{Example}
\newcommand{\res}{\operatorname{res}}
\newcommand{\up}{\operatorname{up}}
\newcommand{\roof}{\operatorname{roof}}
\newcommand{\down}{\operatorname{down}}
\newcommand{\base}{\operatorname{base}}
\newcommand{\ceil}{\operatorname{ceil}}
\newcommand{\floor}{\operatorname{floor}}
\newcommand{\Z}{\mathbb Z}
\renewcommand{\H}{\mathcal H}
\newcommand{\e}{\tilde e}
\newcommand{\f}{\tilde f}

\begin{document}
\maketitle

\begin{abstract}
Let $\{B(\Lambda_m) \mid m\in\Z/e\Z\}$ be the set of level one
$\mathfrak{g}(A^{(1)}_{e-1})$-crystals, and consider the
realization of $B(\Lambda_m)$ using $e$-restricted partitions. We
prove a purely Young diagrammatic criterion for an element of
$B(\Lambda_0)^{\otimes d_1}\otimes B(\Lambda_m)^{\otimes d_2}$ to
be in the component $B(d_1\Lambda_0+d_2\Lambda_m)$. As an
application, we give a non-recursive characterization of simple
modules of the Hecke algebra of type $B$. In the course of the proof,
we also obtain a combinatorial description of
the second type of Kashiwara's Demazure crystal in $B(\Lambda_m)$.
\end{abstract}

\section{Introduction}

Let $\H_n(Q,q)$ be the Hecke algebra of type $B$ defined over an
algebraically closed field $F$ of characteristic $\ell$. The
$F$-algebra $\H_n(Q,q)$ is generated by $T_0,\dots,T_{n-1}$
subject to the quadratic relations $(T_0-Q)(T_0+1)=0$,
$(T_i-q)(T_i+1)=0$, for $1\le i<n$, and the type $B$ braid
relations. Let $q$ be a power of a prime $p\ne\ell$. Motivated by
a desire to generalize their famous work on the classification of
simple $FGL_n(q)$-modules to other classical groups, Dipper and
James initiated the study of modular representations of Hecke
algebras of type $B$, where $q$ is an arbitrary element in $F$.
They proved a certain Morita equivalence theorem \cite[Theorem
4.14]{DJ} and as a result, they classified simple
$\H_n(Q,q)$-modules in the case when $-Q\not\in q^{\mathbb Z}$
\cite[Theorem 5.6]{DJ}. Suppose that $q\ne1$ and $-Q\in q^{\mathbb
Z}$. Then the classification of simple $\H_n(Q,q)$-modules was
achieved in \cite[Theorem 4.2]{A2}, which completed the previous
work \cite{AM}. The classification is given for cyclotomic Hecke
algebras associated with $G(r,1,n)$, which is defined by replacing
$(T_0-Q)(T_0+1)=0$ with $(T_0-v_1)\cdots(T_0-v_r)=0$ in the above
definition.\footnote{By the Morita equivalence theorem for cyclotomic Hecke
algebras proven by Dipper and Mathas \cite[Theorem 4.7]{DM}, we
may assume that $v_i\in q^{\mathbb Z}$, for $1\le i\le r$.}
We note here that Geck-Rouquier theory provides us with
another approach for classifying simple
$\H_n(Q,q)$-modules.\footnote{For the approach in \cite{Gr},
see \cite{A3}.}
The advantage of their approach is that it works for arbitrary
finite Hecke algebras. It is also worth mentioning that
Jacon generalized the theory to cyclotomic Hecke algebras
associated with $G(r,1,n)$. See \cite{Ge1} and \cite{J1}, \cite{J2}.
On the other hand, control of actual modules is rather difficult
in their approach, particularly in the
cyclotomic case.\footnote{Recently Geck has
proved that finite Hecke algebras are cellular \cite{Ge2}. Hence,
they have better control of actual modules than before for
finite Hecke algebras.} Hence, we have needed our approach
in applications such as determination
of representation type, and we are pursuing our
direction further.\footnote{We hope that a better understanding of
the two approaches will lead to the merging of both theories.}

Let $e$ be the multiplicative order of $q\ne1$,
$\mathfrak g$ the Kac-Moody Lie algebra of type $A^{(1)}_{e-1}$,
$\{\Lambda_i \mid i\in \mathbb Z/e\mathbb Z\}$ the fundamental weights.
We realize the Kashiwara crystal $B(\Lambda_i)$ on the set of
$e$-restricted partitions. Suppose that $v_i=q^{\gamma_i}$, for
$1\le i\le r$. Then, our classification theorem asserts that
simple modules are parametrized by the subset
$$
B(\Lambda_{\gamma_1}+\cdots+\Lambda_{\gamma_r})\subset
B(\Lambda_{\gamma_1})\otimes\cdots\otimes B(\Lambda_{\gamma_r}).
$$
In particular, if $-Q=q^m$ then simple $\H_n(Q,q)$-modules are
parametrized by $B(\Lambda_0+\Lambda_m)\subset
B(\Lambda_0)\otimes B(\Lambda_m)$. Further, when
$\lambda\otimes\mu\in B(\Lambda_0+\Lambda_m)$,
we can construct the corresponding simple module
$D^{(\mu,\lambda)}$ as follows.
Let $S^{(\mu,\lambda)}$ be the Specht module for
$\H_n(Q,q)$ constructed by Dipper, James and Murphy
in \cite{DJM}. $S^{(\mu,\lambda)}$ is equipped with
an invariant symmetric bilinear form.
Then $D^{(\mu,\lambda)}$ is the module obtained from
$S^{(\mu,\lambda)}$ by factoring out the radical of
the bilinear form.

A bipartition $(\mu,\lambda)$ is called {\bf Kleshchev} if
$\lambda\otimes\mu\in B(\Lambda_0+\Lambda_m)$. The set of
Kleshchev bipartitions may be computed by applying Kashiwara
operators to the empty bipartition, but this does not give us an
effective method of determining whether a given bipartition is
Kleshchev or not.

The first purpose of this article is to give a non-recursive
characterization of Kleshchev bipartitions. Our result is that
$\lambda\otimes\mu\in B(\Lambda_0+\Lambda_m)$ if and only if
$\roof(\mu)\subset\tau_m(\base(\lambda))$, where $\roof$, $\base$
and $\tau_m$ are explicit operations on abacus displays.
The definition of $\roof$ and $\base$ requires repeated application
of up and down operations respectively, but $\roof$ and $\base$ are easily
computable from a given partition.\footnote{Using this result, the first author
and Jacon have settled a conjecture in \cite{DJM} affirmatively.
See \cite{AJ}.}

The characterization of $B(\Lambda_0+\Lambda_m)$ as a subset of
$B(\Lambda_0)\otimes B(\Lambda_m)$ is a purely crystal theoretic
question.  Due to a result of Littelmann, this characterization
can be expressed in terms of his path model. Our strategy is to
interpret his result in terms of the combinatorics of partitions.
In his result, the initial direction and the final direction of a
Lakshmibai-Seshadri path play an important role, and the crucial
step in proving our theorem is to find a Young diagrammatic
interpretation of these directions. Fortunately, the
interpretation of the initial direction was already given in
\cite{KLMW}. Here, we give the interpretation of the final
direction. This suffices for proving our result for $m=0$.
Combined with arguments which interpret Littelmann's condition for
different dominant integral weights, we reach our
theorem.\footnote{In the path model, an $e$-restricted partition
is given by a sequence of $e$-cores and rational numbers. We show
that the Mullineux map in the modular representation theory of the
symmetric group and the Hecke algebra of type $A$ is given by
conjugation of the $e$-cores. See Proposition \ref{ceil and floor
under Mullineux} and the accompanying remark.}

The second purpose of this article is to describe the crystal
$B^w(\Lambda_m)$ for $w\in W$ in the same way that, in \cite{KLMW},
$B_y(\Lambda_m)$ is described for $y\in W$. The work is motivated by
standard monomial theory \cite{LS}, \cite{L4}. In the Grassmannian case,
see \cite{KL} for a self-contained presentation in the spirit of
the classical work of Hodge and Pedoe \cite{H}, \cite{HP}, and \cite{KLMW2}
for discussion of a similar approach for the affine Grassmannian.

The initial and final directions of a Lakshmibai-Seshadri
path are related to the two types of Demazure crystals $B_y(\Lambda)$
and $B^w(\Lambda)$, for an integral dominant weight $\Lambda$. We
explain the relationship in detail in section 6. The result for
the initial direction is due to Littelmann, and the result for the
final direction is due to Kashiwara and Sagaki, who proved the
result independently. We think that this self-contained
explanation of the results benefits those who have an interest in
Littelmann's path model.

The project started when the first author learned the idea of
using Littelmann's result and the existence of
\cite{KLMW} from Mark Shimozono.
We are grateful to him. We are also grateful to Kashiwara for
his permission to include his proof of the above mentioned
result in this paper. Finally, the first author thanks
Naito and Sagaki for explaining to him basic facts about Littelmann's
path model, and Mathas and Fayers for explaining to him their results for
$e=2$ and $e=3$, which give a different characterization of
Kleshchev bipartitions without using Littelmann's result.
We discuss their results in the last section.

\section{Preliminaries}
We assume that the reader is familiar with
the theory of Kashiwara crystals. The three books
\cite{HK}, \cite{Jo} and \cite{K1} are standard references.
Throughout the paper, we always consider
$\mathfrak g(A^{(1)}_{e-1})$-crystals, for fixed $e\ge2$.

Let $\{\Lambda_m \mid m\in \Z/e\Z\}$ be
the set of fundamental weights. We denote by $B(\Lambda_m)$
the Kashiwara crystal associated with $\Lambda_m$.
Recall that a {\bf partition} $\lambda$ is a sequence of non-increasing
integers
$$
\lambda_0\ge\lambda_1\ge\cdots
$$
which has only a finite number of nonzero elements. We denote
$\lambda_0$ by $a(\lambda)$. When $\lambda_{\ell-1}>0$ and
$\lambda_\ell=0$, we denote $\lambda=(\lambda_0,\dots,\lambda_{\ell-1})$
and denote $\ell$ by $\ell(\lambda)$.
A partition is called {\bf $e$-restricted} if
$0\le\lambda_i-\lambda_{i+1}<e$, for all $i$.

We shall recall the realization of $B(\Lambda_m)$ in terms of
$e$-restricted partitions.
Let $\lambda$ be a partition. We color the nodes of $\lambda$ with
the $e$ colors $\Z/e\Z$ by the following rule: let
$x(a,b)$ be the node located on the $a^{th}$ row and the $b^{th}$ column.
Then $x(a,b)$ has color $m-a+b+e\Z$. The number $m-a+b$ is called
the {\bf content} of $x(a,b)$, and the color $m-a+b+e\Z$ is called
the {\bf residue} of $x(a,b)$. Let $\lambda\subset\mu$ be a pair
of partitions such that the number of nodes differs by one.
Suppose that the residue of the node $x=\mu\setminus\lambda$ is $i$.
Then we call $x$ an {\bf addable $i$-node} of $\lambda$ and
a {\bf removable $i$-node} of $\mu$.

Let $B$ be the set of $e$-restricted partitions.
We color the nodes of $\lambda\in B$ as above, and define
$$
wt(\lambda)=\Lambda_m-\sum_{i\in\mathbb Z/e\mathbb Z}
N_i(\lambda)\alpha_i
$$
where $N_i(\lambda)$ is the number of $i$-nodes in $\lambda$.
In order to define two operators $\f_i$ and $\e_i$
on $B\sqcup\{0\}$, we read addable $i$-nodes and
removable $i$-nodes from
the first row to the last row and record the result as
a sequence of $A$'s and $R$'s.
Then we apply an algorithm which we call {\bf $RA$-deletion}.
Choose any $R \cdots A$, where the middle $\cdots$ means the letters
which have been already deleted, and change it to $\cdot$ $\cdots$ $\cdot$.
We repeat this procedure as many times as possible.
The final sequence is of the form
$$
\cdots A \cdots A \cdots A \cdots R \cdots R \cdots R \cdots R \cdots
$$
where $\cdots$ is a sequence of dots of length greater than or
equal to $0$. The final sequence is uniquely determined
(see \cite[Lemma 11.2]{A1}). The nodes which appear in the final
sequence are called {\bf addable normal $i$-nodes} and
{\bf removable normal $i$-nodes}.
We define $\f_i\lambda$ to be the partition
obtained from $\lambda$ by adding the node which corresponds to
the rightmost $A$ in the final sequence.
If there is no $A$ in the final sequence,
we set $\f_i\lambda=0$.
Similarly, we define $\e_i\lambda$ to be
the partition obtained from $\lambda$ by removing the node
which corresponds to the leftmost $R$ in the final sequence, and $0$ if no $R$
exists in the final sequence. Finally, we define
$\f_i0=0$ and $\e_i0=0$. Define
$$
\varphi_i(\lambda)=\max\{k\in\Z_{\ge0} \mid \f_i^k\lambda\ne0\}, \;\;
\epsilon_i(\lambda)=\max\{k\in\Z_{\ge0} \mid \e_i^k\lambda\ne0\}.
$$
In other words, $\varphi_i(\lambda)$ is the number of $A$'s in the
final sequence, and $\epsilon_i(\lambda)$ is the number of $R$'s in
the final sequence.

The set $B$ with the additional
data $wt$, $\epsilon_i$, $\varphi_i$, $\e_i$ and $\f_i$
is a realization of the crystal $B(\Lambda_m)$. This result
is due to Misra and Miwa. See \cite[Theorem 11.11]{A1}.
We denote the empty partition in $B(\Lambda_m)$ by $\emptyset_m$.

It is convenient to work with the abacus display of $\lambda$.
The set of {\bf beta numbers of charge} $m$ associated with $\lambda$
is, by definition, the set $J$ of decreasing integers
$$
j_0>j_1>j_2>\cdots>j_k>\cdots
$$
defined by $j_k=\lambda_k+m-k$, for $k\ge0$.
It has the property that $j_k=m-k$, for $k>>0$.
We consider an abacus with $e$ runners
\begin{equation*}
\begin{aligned}
\cdot& \hphantom{=}& \cdot& \hphantom{=}& \cdots& \hphantom{=}&\cdot\\
\cdot& \hphantom{=}& \cdot& \hphantom{=}& \cdots& \hphantom{=}&\cdot\\
0 & \hphantom{=}& 1   &\hphantom{=}& \dots & \hphantom{=}& e-1 \\
e & \hphantom{=}& e+1 &\hphantom{=}& \cdots& \hphantom{=}&2e-1\\
\cdot& \hphantom{=}& \cdot& \hphantom{=}& \cdots& \hphantom{=}&\cdot\\
\cdot& \hphantom{=}& \cdot& \hphantom{=}& \cdots& \hphantom{=}&\cdot
\end{aligned}
\end{equation*}
and put beads on the numbers $\{j_k \mid k\ge0\}$. This is the {\bf abacus display}
of charge $m$ associated with $\lambda$.

\begin{Ex}
\label{J records contents}
Let $e=3$, $m=0$, and $\lambda=(4,2,1)$.

To read $J$ from $\lambda$, we look at each row and
find the content of the node
which is adjacent to the right end of the row.
\begin{center}
\begin{tabular}{cccccc}
$\times$ & $\times$ & $\times$ & $\times$ & $4$\\
$\times$ & $\times$ & $1$      &          &    \\
$\times$ & $-1$ &&&\\
$-3$ &&&&\\
$-4$ &&&&\\
$\cdot$&&&&\\
$\cdot$&&&&
\end{tabular}
\end{center}
Thus, $J=\{4,1,-1,-3,-4,\dots\}$, and the abacus display of $\lambda$ is
as follows.
\begin{equation*}
\begin{aligned}
\cdot& \hphantom{=}& \cdot&\hphantom{=}&\cdot\\
-6   & \hphantom{=}& -5   &\hphantom{=}&-4  \\
-3   & \hphantom{=}&      &\hphantom{=}&-1  \\
     & \hphantom{=}&  1   &\hphantom{=}&    \\
     & \hphantom{=}&  4   &\hphantom{=}&
\end{aligned}
\end{equation*}
\end{Ex}

We call $j\in J$ with $j+e\Z=i+1$ a {\bf removable $i$-integer}, and
$j\in J$ with $j+e\Z=i$ an {\bf addable $i$-integer}. The Kashiwara
operators $\e_i$ and $\f_i$ in terms of $J$ are given by the same
procedure as above. We change the sequence $j_0,j_1,\dots$ to a
sequence of dots, $R$'s, and $A$'s, and apply the $RA$-deletion as
many times as possible. Note that a removable or addable integer
$j\in J$ may not correspond to a removable or addable node of
$\lambda$. However, this happens precisely when
$\lambda_k=\lambda_{k+1}$. In this case, the content of the node
which is adjacent to the right end of row $k$ is a removable
$i$-integer, and the content of the node which is adjacent to the
right end of row $k+1$ is an addable $i$-integer. In
$RA$-deletion, these two adjacent values are removed from the
final sequence.

The following definition is given in \cite{KLMW}.

\begin{defn}
Let $\lambda\in B(\Lambda_m)$ and $J$ the corresponding set of
beta numbers of charge $m$. Let $U(J)$ be the set of beads which
we may slide up by one in their runners. In other words,
$$
U(J)=\{x\in J \mid x-e\not\in J\}.
$$
If $U(J)=\emptyset$ then define $\up(\lambda)=\lambda$.
Suppose $U(J)\ne\emptyset$. Then
set $p=\max\,U(J)$ and consider
$$
V(J)=\{x>p \mid x\not\in p+e\Z, x-e\in J, x\not\in J\}.
$$
Set $q=\min\,V(J)$. Then we define $\up(J)$ to be
the set $(J\setminus\{p\})\cup\{q\}$. That is, we obtain
$\up(J)$ by moving the bead $p$ to $q$. We denote the corresponding
partition by $\up(\lambda)$.
\end{defn}

\begin{Ex}
Let $e=3$, $m=2$ and $\lambda=(3,2,1)$. Then the abacus display
of $\lambda$ is
\begin{equation*}
\begin{aligned}
\cdot& \hphantom{=}& \cdot&\hphantom{=}&\cdot\\
-3   & \hphantom{=}& -2   &\hphantom{=}&-1  \\
     & \hphantom{=}&  1   &\hphantom{=}&    \\
 3   & \hphantom{=}&      &\hphantom{=}& 5
\end{aligned}
\end{equation*}
Then $U(J)=\{3,5\}$ and $p=5$. Thus $V(J)=\{6\}$ and $q=6$.
Therefore, $\up(J)$ is given by
\begin{equation*}
\begin{aligned}
\cdot& \hphantom{=}& \cdot&\hphantom{=}&\cdot\\
-3   & \hphantom{=}& -2   &\hphantom{=}&-1  \\
     & \hphantom{=}&  1   &\hphantom{=}&    \\
 3   & \hphantom{=}&      &\hphantom{=}&    \\
 6   & \hphantom{=}&      &\hphantom{=}&
\end{aligned}
\end{equation*}
Thus, $\up(\lambda)=(4,2,1)$.
\end{Ex}

\begin{lemma}
\label{up operation}
Let $\lambda\in B(\Lambda_m)$.
\item[(1)]
$\lambda\subset \up(\lambda)$.
\item[(2)]
$\up(\lambda)$ is $e$-restricted.
\item[(3)]
$\ell(\up(\lambda))=\ell(\lambda)$.
\end{lemma}
\begin{proof}
(1) Let $j'_0>j'_1>\cdots$ be the beta numbers of charge $m$ associated
with $\up(\lambda)$. We set $j_{-1}=\infty$. Then, there exists $s\ge-1$
such that $j_s>q\ge j_{s+1}$. $q\not\in J$ implies that $q>j_{s+1}$.
Since $q>p$, there also exists $t>s$ such that $j_t=p$. Then, for $t>s+1$,
$$
\begin{cases}j'_k=j_k\quad (0\le k\le s)\\
             j'_{s+1}=q>j_{s+1} \\
             j'_k=j_{k-1}>j_k \quad(s+1<k<t)\\
             j'_t=j_{t-1}>j_t=p \\
             j'_k=j_k \quad(k\ge t+1)\end{cases}
$$
If $t=s+1$, replace the middle three lines with $j'_{s+1}=q>j_t=p$.
In any case, $j_k'\ge j_k$, for all $k$. This implies the result.

(2) We only have to check the effect of removing $p$. We want to show
$j_t'-j_{t+1}'\le e$.
Since $\lambda$ is $e$-restricted and $p-e\not\in J$, there exists
$x\in\{p-e+1,\dots,p-1\}\cap J$. Note that $j_{t+1}$ is the largest element of
$J$ which is smaller than $j_t=p$. Thus we have $x\le j_{t+1}=j_{t+1}'$.

Suppose first that $x+e\not\in J$. Then
$q\le x+e$, which implies that
$$
x\le j'_{t+1}<j'_t\le j_{s+1}'=q\le x+e.
$$
Thus, $\up(\lambda)$ is $e$-restricted.

Suppose next that $x+e\in J$. Then $j_t=p<x+e$ implies
$j_t'=j_{t-1}\le x+e$. Thus
$x\le j'_{t+1}<j'_t\le x+e$ and
$\up(\lambda)$ is $e$-restricted.

(3) Let $s\in\Z$ be such that $\Z_{\leq s}\subset J$ and $s+1\not\in J$.
Then $\ell(\lambda)=|\{x\in J \mid x>s\}|$. As $p>s$ and
$p$ moves to $q>p$, we have $\Z_{\leq s}\subset J'$ and $s+1\not\in J'$,
which implies $\ell(\up(\lambda))=|\{x\in J' \mid x>s\}|$,
and $\ell(\up(\lambda))=\ell(\lambda)$.
\end{proof}

We remark that we may deduce $\lambda\subset\up(\lambda)$ from
$|J'\cap\mathbb Z_{\ge a}|\ge |J\cap\mathbb Z_{\ge a}|$, for
all $a\in\mathbb Z$. In fact, if there existed
$k\ge0$ such that $j_0'=j_0,\dots,j'_{k-1}=j_{k-1}$ and
$j'_k<j_k$, then we would obtain
$|J'\cap\mathbb Z_{\ge j_k}|<|J\cap\mathbb Z_{\ge j_k}|$,
a contradiction.

If we apply the up operation successively, then we reach
$U(J)=\emptyset$ after finitely many steps. To see this, choose
$s$ such that $\mathbb Z_{\le s}\subset J$. Then $\mathbb Z_{\le
s}\subset\up(J)$. Thus, $\mathbb Z_{\le s}$ remains untouched
during the successive applications of up operations. Let $N$ be
the number of elements in $\{x\in J \mid x>s\}$ and $K=\mathbb
Z_{\le s}\cup\{s+ke \mid 1\le k\le N\}$. We write $J\le J'$ if
$j_k\le j'_k$, for all $k\ge0$. Note that if $J$ is the set of
beta numbers associated with an $e$-restricted partition and of
the form $J=\mathbb Z_{\le s}\cup\{j_0,\dots,j_{N-1}\}$, where
$j_0>\cdots>j_{N-1}>s$, then $J\le K$. Thus, we have $\up^i(J)\le
K$, for all $i\ge0$. As the sequence $J, \up(J), \up^2(J),\dots$
is strictly increasing as long as $U(J)\ne\emptyset$, we reach
$U(J)=\emptyset$ after finitely many steps.

This allows us to define $\roof(J)$ as follows.

\begin{defn}
Let $\lambda\in B(\Lambda_m)$ and let $J$ be as before.
Apply the up operation to $J$
until $U(J)=\emptyset$. We denote the resulting
$\up^{max}(J)$ by $\roof(J)$, and denote the corresponding partition by
$\roof(\lambda)$.
\end{defn}
\noindent Note that by definition, $\roof(\lambda)$ is an
$e$-core.

\begin{defn}
Let $\lambda\in B(\Lambda_m)$ and $J$ the corresponding set of
beta numbers of charge $m$. Let $U(J)$ be as before.
If $U(J)=\emptyset$ then define $\down(\lambda)=\lambda$.
Suppose $U(J)\ne \emptyset$. Then
set $p'=\min\,U(J)$ and consider
$$
W(J)=\{x>p'-e \mid x\in J, x+e\not\in J\}\cup\{p'\}.
$$
Set $q'=\min\,W(J)$. Then we define $\down(J)=
(J\setminus\{q'\})\cup\{p'-e\}$. That is, we obtain
$\down(J)$ by moving the bead $q'$ to $p'-e$. We denote the corresponding
partition by $\down(\lambda)$.
\end{defn}

\begin{lemma}
\label{down operation}
Let $\lambda\in B(\Lambda_m)$.
\item[(1)]
$\lambda\supset \down(\lambda)$.
\item[(2)]
$\down(\lambda)$ is $e$-restricted.
\item[(3)]
$a(\down(\lambda))=a(\lambda)$.
\end{lemma}
\begin{proof}
(1) Let $j'_0>j'_1>\cdots$ be the beta numbers of charge $m$ associated
with $\down(\lambda)$. Then, there exists $s\ge0$ such that
$j_s>p'-e>j_{s+1}$,
and there exists $0\le t\le s$ such that $j_t=q'$. Now,
$$
\begin{cases}j'_k=j_k\quad (0\le k<t)\\
             j'_t=j_{t+1}<j_t=q' \\
             j'_k=j_{k+1}<j_k \quad(t<k<s)\\
             j'_s=p'-e<j_s\\
             j'_k=j_k \quad(k\ge s+1)\end{cases}
$$
We replace the middle three lines with $j'_t=p'-e<j_t=q'$ when
$t=s$. Thus $j_k'\le j_k$, for all $k$. This implies the result.

(2) We only have to consider the effect of removing $q'$ as before. We want to
show $j_{t-1}'-j_t'\le e$. Note that
there exists $x\in\{p'-e+1,\dots,p'-1\}\cap J$ since $\lambda$ is
$e$-restricted and $p'-e\not\in J$.

Suppose first that $q'\ne p'$. Then $p'\ge j_{t-1}'$ since
$p'>q'$ and $j_{t-1}'=j_{t-1}$
is the smallest element of $J$ which is greater than $j_t=q'$.
Thus
$$
p'-e=j_s'\le j_t'<j_{t-1}'\le p'.
$$

Suppose next that $q'=p'$. There exists $x\in\{p'-e+1,\dots,p'-1\}\cap J$
as before. As $x<p'=j_t$ and $x\in J$, we have
$x\le j_{t+1}$. On the other hand,
$q'=p'$ implies that $x+e\not\in J$ is impossible. Thus,
$j_t=p'<x+e$ implies $j_{t-1}\le x+e$ and
$$
x\le j_{t+1}=j_t'<j_{t-1}'=j_{t-1}\le x+e.
$$

(3) As $a(\lambda)=j_0-m$ and $a(\down(\lambda))=j_0'-m$, we show
$j_0=j_0'$. If $p'<j_0$ then $q'\leq p'<j_0$.
If $p'=j_0$ then $j_0-e\not\in J$ and, since $\lambda$ is $e$-restricted,
there exists $x\in J$ such that $j_0-e<x<j_0$. Then, as $x+e\not\in J$,
$x\in W(J)$ and $q'\leq x<j_0$. Hence $q'<j_0$ in both cases and
$q'$ moves to $p'-e<q'$. Thus $j_0=j_0'$.
\end{proof}

As before, we may deduce $\lambda\supset\down(\lambda)$ from
$|J'\cap\mathbb Z_{\ge a}|\le |J\cap\mathbb Z_{\ge a}|$ for all
$a\in\mathbb{Z}$.

We apply the down operation successively.
It is easy to see that we reach $U(J)=\emptyset$ after finitely many steps:
the size of the corresponding partition
strictly decreases as long as $U(J)\ne\emptyset$.
In section $7$, we need a better understanding of how the value $p'$ changes
during the process. Thus, we analyze it in detail here.

Suppose that we apply the down operation to $J^{old}$ to
obtain $J^{new}$ and that $U(J^{old})\ne\emptyset$ and
$U(J^{new})\ne\emptyset$.
Since ${p'}^{old}\not\in U(J^{new})$ implies ${p'}^{new}\ne {p'}^{old}$,
we have either ${p'}^{new}>{p'}^{old}$ or ${p'}^{new}<{p'}^{old}$.

Suppose that ${p'}^{new}<{p'}^{old}$.
If ${p'}^{new}-e\not\in J^{old}$ then ${p'}^{new}\not\in J^{old}$
as ${p'}^{new}\in J^{old}$ would imply ${p'}^{new}\ge {p'}^{old}$.
Hence ${p'}^{new}\in J^{new}\setminus J^{old}$ and we have
${p'}^{new}={p'}^{old}-e$.

The set $U(J)$ changes in the following way. Let
$q'=\min W(J^{old})$.
\begin{itemize}
\item[(a)]
If $q'<{p'}^{old}$ then $q'-ke\in J^{old}$, for all $k\ge0$,
and $q'+e\not\in J^{old}$. Hence,
$$
U(J^{new})\setminus\{{p'}^{new}\}\subset U(J^{old})\setminus\{{p'}^{old}\}.
$$
\item[(b)]
If $q'={p'}^{old}$ then
$$
U(J^{new})\setminus\{{p'}^{new}\}\subset
\left(U(J^{old})\setminus\{{p'}^{old}\}\right)
\cup\{{p'}^{old}+e\}.
$$
\end{itemize}

Next suppose that $p'$ starts decreasing at $p_0=\min U(J_0)$
and stops decreasing at $p_N=\min U(J_N)$.
By the above consideration, the innovation of $p'$ is given by
the recursion ${p'}^{new}={p'}^{old}-e$, so
$p_k=p_0-ke$, for $0\le k\le N$. Denote $q_k=\min W(J_k)$.
Define $s\ge0$ by $q_k=p_k$, for $0\le k<s$, and $q_s\ne p_s$.
We shall show by induction on $k$ that
$$
U(J_k)\cap \mathbb Z_{\le p_0}=\{p_0-ke\},\;\;
\text{for $0\le k\le N$.}
$$
For $0\le k<s$, $J_{k+1}$ is obtained from $J_k$ by
sliding the bead at $p_0-ke$ up to $p_0-(k+1)e$. Thus,
if $k\ge1$ and $x\le p_0$ is such that $x\in J_{k+1}$ and
$x+e\mathbb Z=p_0+e\mathbb Z$, then
$x\le p_0-(k+1)e$.
Suppose that ${p'}^{old}+e\in U(J^{new})$ occured at $k\ge1$.
Thus, ${p'}^{new}=p_0-(k+1)e$ and ${p'}^{old}+e=p_0-(k-1)e$.
Let $x={p'}^{old}+e$. Then $x\le p_0$ satisfies
$x\in J_{k+1}$ and $x+e\mathbb Z=p_0+e\mathbb Z$ but
$x>p_0-(k+1)e$. Thus, ${p'}^{old}+e\not\in U(J^{new})$ and
$$
U(J_{k+1})\setminus\{p_{k+1}\}\subset
U(J_k)\setminus\{p_k\}\subset \mathbb Z_{\ge p_0+1}
$$
by the induction hypothesis.

For $s\le k\le N$, we have $q_k<p_k$. If $k=s$
then this is by definition. Suppose that
$q_k<p_k$. Then $q_k-e\in J_k$ and
$q_k\not\in J_{k+1}$, $p_{k+1}-e=p_k-2e<q_k-e$
imply $q_k-e\in W(J_{k+1})$. Hence,
$$
q_{k+1}\le q_k-e<p_k-e=p_{k+1}.
$$
Therefore, we have
$$
U(J_{k+1})\setminus\{p_{k+1}\}\subset
U(J_k)\setminus\{p_k\}\subset \mathbb Z_{\ge p_0+1},
$$
for $s\le k<N$. We have proved that
$U(J_k)\cap \mathbb Z_{\le p_0}=\{p_k\}$,
for $0\le k\le N$.

Now, set $J^{old}=J_N$ and ${p'}^{old}=p_N$. Then we obtain $J^{new}$
from $J^{old}$ by moving $q_N$ to $p_N-e$. Suppose that $U(J^{new})\ne\emptyset$.
Then ${p'}^{new}>{p'}^{old}$.

We claim that
${p'}^{new}>p_0$. In fact, as $p_N-e\not\in U(J^{new})$, we have either
${p'}^{new}\in U(J_N)\setminus\{p_N\}$ or ${p'}^{new}=q_N+e$. In
the former case, $U(J_N)\cap \mathbb Z_{\le p_0}=\{p_N\}$ implies
${p'}^{new}>p_0$. Suppose that ${p'}^{new}=q_N+e\le p_0$. If $q_k<p_k$ for
some $k\le N$, then $q_N<p_N$, and $p_N<q_N+e\not\in J_N$
implies $q_N+e\not\in J^{new}$, which contradicts ${p'}^{new}\in J^{new}$. If
$q_k=p_k$ for $0\le k\le N$, then $p_0-ke$ is not contained in
$J^{new}$, for $0\le k\le N$. So
$q_N+e=p_0-(N-1)e\not\in J^{new}$ either. We have proved that
${p'}^{new}>p_0$.

As we reach $U(J)=\emptyset$ after finitely many steps, we may
define $\base(J)$ as follows.

\begin{defn}
Let $\lambda\in B(\Lambda_m)$ and let $J$ be as before.
Apply the down operation to $J$ until
$U(J)=\emptyset$. We denote the resulting
$\down^{max}(J)$ by $\base(J)$, and denote the corresponding partition by
$\base(\lambda)$.
\end{defn}
\noindent Note that $\base(\lambda)$ is an $e$-core by definition.

\section{Weyl group action}

Let $B$ be a $\mathfrak g$-crystal and $W$ the corresponding
Weyl group. In our case of $B(\Lambda_m)$, $W$ is the
Coxeter group generated by $\{s_i \mid i\in \Z/e\Z\}$ subject to
$s_i^2=1$,
$s_is_{i+1}s_i=s_{i+1}s_is_{i+1}$ and $s_is_j=s_js_i$ otherwise.

\begin{thm}[{\cite[Theorem 9.4.1]{K1}}]
Let $B$ be a normal crystal. Then the following defines
a $W$-action on $B$.
$$
s_ib=\begin{cases}\f_i^{wt(b)(h_i)}b\quad (\text{if}\;\;wt(b)(h_i)\ge0)\\
                  \e_i^{-wt(b)(h_i)}b \quad (\text{if}\;\;wt(b)(h_i)\le0)\end{cases}
$$
Further, $wt(s_ib)=s_i(wt(b))=wt(b)-wt(b)(h_i)\alpha_i$.
\end{thm}

Recall that $B(\Lambda_m)$ is a normal crystal. Hence, we have a
$W$-action and
$$
\e_i^{max}\lambda=\e_i^{\epsilon_i(\lambda)}\lambda,\quad
\f_i^{max}\lambda=\f_i^{\varphi_i(\lambda)}\lambda.
$$

\begin{defn}
Let $\lambda\in B(\Lambda_m)$. We say that $\lambda$ is an
{\bf $s_i$-core} if
$x\in U(J)$ implies $x+e\mathbb Z\ne i$ and $x+e\mathbb Z\ne i+1$.
\end{defn}

Thus, $\lambda$ is an $e$-core if and only if it is an $s_i$-core,
for all $i\in\mathbb Z/e\mathbb Z$.

\begin{lemma}
\label{Weyl group action}
Suppose that $\lambda\in B(\Lambda_m)$.
\item[(1)]
Let $A_i(\lambda)$ and $R_i(\lambda)$ be the set of
addable $i$-nodes and the set of removable $i$-nodes of
$\lambda$ respectively. Then
$$
wt(\lambda)(h_i)=|A_i(\lambda)|-|R_i(\lambda)|.
$$
\item[(2)]
Assume that $\lambda$ is an $s_i$-core. Then
either
\begin{itemize}
\item[(i)]
$A_i(\lambda)=\emptyset$ and
$s_i\lambda=\e_i^{max}\lambda=
\lambda\setminus\{\text{all removable $i$-nodes}\}$, or
\item[(ii)]
$R_i(\lambda)=\emptyset$ and
$s_i\lambda=\f_i^{max}\lambda=
\lambda\cup\{\text{all addable $i$-nodes}\}$.
\end{itemize}
\end{lemma}
\begin{proof}
(1) is proved by induction on $|\lambda|$. If $\lambda=\emptyset_m$,
$\Lambda_m(h_i)=\delta_{im}$ proves the result.
Suppose that $\lambda=\mu\cup\{x\}$ and the residue of $x$ is $j$.
Thus $wt(\lambda)=wt(\mu)-\alpha_j$.
Note that
$$
wt(\lambda)(h_i)=
\begin{cases}
wt(\mu)(h_i)\quad (j\ne i, i\pm1)\\
wt(\mu)(h_i)+1 \quad(j=i\pm1)\\
wt(\mu)(h_i)-2 \quad(j=i)\end{cases}
$$
Checking how $A_i(\mu)$ and $R_i(\mu)$ change when  $x$ is added,
we obtain the result.

(2) For a hook $\Gamma=(a,1^r)$, the $a$ nodes consist the arm of
$\Gamma$ and the $r$ nodes consist the leg of $\Gamma$.
The residue of the lowest node of the leg is called
the residue of $\Gamma$. Let $J$ be the set of beta
numbers of charge $m$ associated with $\lambda$.
Recall that sliding a
bead in $J$ on the $i^{th}$ runner up by one
is the same as removing an $e$-hook $\Gamma$
whose residue is $i$. Suppose that there exist
$x\in A_i(\lambda)$ and $y\in R_i(\lambda)$ such that
$x$ is in the $j^{th}$ row of $\lambda$ and $y$ is in the
$k^{th}$ row of $\lambda$. If $j<k$ then we may remove at least one
$e$-hook of residue $i$ from $\lambda$. Similarly,
if $j>k$ then we may remove at least one $e$-hook of
residue $i+1$ from $\lambda$.
Since $\lambda$ is an $s_i$-core, both cannot occur.
In other words, one of $A_i(\lambda)$ or
$R_i(\lambda)$ must be empty.
Thus, $RA$-deletion does not occur, which
implies that
either $\epsilon_i(\lambda)=|R_i(\lambda)|$ and $\varphi_i(\lambda)=0$,
or $\epsilon_i(\lambda)=0$ and $\varphi_i(\lambda)=|A_i(\lambda)|$
respectively. Now the result follows from (1).
\end{proof}

We show that this Weyl group action coincides that of \cite{KLMW}
on $e$-cores.

\begin{lemma}
\label{lemma for cores}
Let $\lambda$ be an $s_i$-core, $J$ the corresponding set of beta numbers
of charge $m$. We denote by $s_iJ$ the set of beta numbers of
charge $m$ associated with $s_i\lambda$.
\item[(1)]
If $i\ne e-1$ then $s_iJ$ is obtained by switching the $i^{th}$ and
$(i+1)^{th}$ runners.
\item[(2)]
The $(e-1)^{th}$ runner of $s_{e-1}J$ is obtained from
the $0^{th}$ runner of $J$ by sliding up by one.
Similarly, the $0^{th}$ runner of $s_{e-1}J$ is obtained from the
$(e-1)^{th}$ runner of $J$ by sliding down by one.
\item[(3)]
$s_i\lambda$ is an $s_i$-core.
\end{lemma}
\begin{proof}
(1) If the length of the $i^{th}$ runner of $J$ exceeds that of
$(i+1)^{th}$ runner by $k$, these $k$ beads correspond to
addable $i$-nodes of $\lambda$. Thus, Lemma \ref{Weyl group action} (2)
implies that $s_i\lambda$ is obtained from $\lambda$ by adding
all the addable $i$-nodes. The resulting $s_iJ$ is the same as the
one which is obtained by switching the two runners.
If the length of the $(i+1)^{th}$ runner of $J$ exceeds that of
$i^{th}$ runner by $k$, these $k$ beads correspond to
removable $i$-nodes of $\lambda$. Thus, Lemma \ref{Weyl group action} (2)
again implies that $s_iJ$ is obtained from $J$ by
switching the two runners.

The proof of (2) is entirely similar to that of (1) and
(3) is an obvious consequence of (1) and (2).
\end{proof}

The following proposition seems to be well-known, but we could
not find a reference.

\begin{prop}
\label{core}
The set of $e$-cores in $B(\Lambda_m)$ coincides the $W$-orbit
through $\emptyset_m$.
\end{prop}
\begin{proof}
We can prove that an $e$-core belongs to $W\emptyset_m$ by
induction on $|\lambda|$. Let $x$ be the right end of the
last row of $\lambda$, and let $i$ be the residue of $x$.
Set $\mu=\e_i^{max}\lambda$. Then
$|\mu|<|\lambda|$ since $x$ is a
removable normal $i$-node, and
$\lambda=\f_i^{max}\mu$ since $\lambda$ is an $e$-core.
Since the set of $e$-cores is stable under $W$-action by
Lemma \ref{lemma for cores} (3), $\mu$ is again an $e$-core,
so $\mu\in W\emptyset_m$ by the induction hypothesis.
Thus, we have $\lambda=s_i\mu\in W\emptyset_m$.
Since a non-empty $W$-stable subset of a $W$-orbit must
coincide with the $W$-orbit itself,
we have the result.
\end{proof}

\begin{defn}
Let $W_m$ be the subgroup of $W$ generated by
$\{s_i \mid i\ne m\}$.
We denote by $W/W_m$ the set of
distinguished coset representatives.
\end{defn}

As $W_m$ is the Coxeter group of type $A_{e-1}$, $W_m$ has the
longest element. Thus the following definition makes sense.

\begin{defn}
We denote by $w_m$ the longest element of $W_m$.
\end{defn}

Recall that $W$ becomes a poset by the Bruhat-Chevalley order.
We write $u\le v$, for $u, v\in W$.
By virtue of Proposition \ref{core}, each $e$-core $\lambda\in B(\Lambda_m)$
can be written in the form $\lambda=w\emptyset_m$, for $w\in W/W_m$, in
a unique manner.

\section{Demazure crystal}

Following \cite{K1} and \cite{K3}, we introduce two types of Demazure crystals.

\begin{defn}
Let $y, w\in W$ and let $y=s_{i_1}\cdots s_{i_\ell}$ be a reduced
expression for $y$. Then we define $B_y(\Lambda_m)$ and
$B^w(\Lambda_m)$ as follows.
\begin{equation*}
\begin{split}
B_y(\Lambda_m)&=
\{\f_{i_1}^{a_1}\cdots\f_{i_\ell}^{a_\ell}\emptyset_m \mid
(a_1,\dots,a_\ell)\in (\Z_{\ge0})^\ell\}\setminus\{0\}, \\
B^w(\Lambda_m)&=\{b\in B(\Lambda_m) \mid G_v(b)\in U_v^-(\mathfrak g)u_{w\Lambda_m}\}.
\end{split}
\end{equation*}
\end{defn}

By \cite[Proposition 9.1.3, 9.1.5]{K1},
$B_y(\Lambda_m)$ does not depend on the choice of the reduced
expression. For the notations $G_v(b)$ and $u_{w\Lambda_m}$, see \S6.

The following are fundamental properties of the Demazure crystals.
The results hold for any dominant integral weight.

\begin{prop}[{\cite[Proposition 3.2.3, 3.2.4, 4.3, 4.4]{K3}}]
\label{basic properties}
\item[(1)]
$\e_iB_y(\Lambda_m)\subset B_y(\Lambda_m)\cup\{0\}$ and
$\f_iB^w(\Lambda_m)\subset B^w(\Lambda_m)\cup\{0\}$.
\item[(2)]
If $s_iy<y$ then $B_y(\Lambda_m)=\cup_{k\ge0}
\f_i^kB_{s_iy}(\Lambda_m)\setminus\{0\}$.
\item[(3)]
If $s_iw>w$ then $B^w(\Lambda_m)=\cup_{k\ge0}
\e_i^kB^{s_iw}(\Lambda_m)\setminus\{0\}$.
\item[(4)]
Let $y, w\in W/W_m$. Then the following are equivalent.
\begin{itemize}
\item[(i)]
$y\ge w$.
\item[(ii)]
$B^w(\Lambda_m)\cap B_y(\Lambda_m)\neq\emptyset$.
\item[(iii)]
$B_w(\Lambda_m)\subset B_y(\Lambda_m)$.
\item[(iv)]
$w\emptyset_m\in B_y(\Lambda_m)$.
\item[(v)]
$B^y(\Lambda_m)\subset B^w(\Lambda_m)$.
\item[(vi)]
$y\emptyset_m\in B^w(\Lambda_m)$.
\end{itemize}
\end{prop}

Next theorem is the main result of \cite{KLMW}. However,
the proof we will give is slightly different from
the original: see Theorem \ref{ceiling theorem},
Theorem \ref{roof lemma} and Corollary \ref{roof is ceil}.

\begin{thm}[{\cite[Theorem 1.1]{KLMW}}]
\label{KLMW}
In the partition realization of $B(\Lambda_m)$, we have
$$
B_y(\Lambda_m)=\{\lambda\in B(\Lambda_m) \mid \roof(\lambda)\subset y\emptyset_m\}.
$$
\end{thm}

\begin{prop}
\label{two orders}
Let $\lambda=u\emptyset_m$ and $\mu=v\emptyset_m$, for $u, v\in W$.
\item[(1)]
If $u\le v$ then $\lambda\subset\mu$.
\item[(2)]
If $\lambda\subset\mu$ and $u, v\in W/W_m$ then $u\le v$.
\end{prop}
\begin{proof}
(1) We prove this by induction on $\ell(v)$.
Let $v=s_is_{i_2}\cdots s_{i_\ell}$ be a reduced expression.
Then $u$ is a subword of the expression.

First we suppose that the leftmost $s_i$ does not
appear in this subword. Then $u\le s_iv$ and the induction hypothesis
implies that
$$
\lambda=u\emptyset_m\subset s_iv\emptyset_m=s_i\mu.
$$
Write $w=s_iv$. Then $w<s_iw$ since $s_iv<v$.
If $w^{-1}\alpha_i$ were a negative root, then
the standard argument would show that $w>s_iw$.
Hence $w^{-1}\alpha_i$ is a positive root. In other words,
$v^{-1}\alpha_i$ is a negative root and
$\langle \Lambda_m, v^{-1}h_i\rangle\le0$. We have
$$
wt(s_i\mu)=wt(s_iv\emptyset_m)=s_iv\Lambda_m=
v\Lambda_m-\langle \Lambda_m, v^{-1}h_i\rangle\alpha_i.
$$
Hence $wt(s_i\mu)-wt(\mu)\in \sum_{j\in\Z/e\Z}\Z_{\ge0}\alpha_j$.
Note that
$$
\begin{cases}
wt(\mu)&=\Lambda_m-\sum_{j\in\Z/e\Z} N_j(\mu)\alpha_j,\\
wt(s_i\mu)&=\Lambda_m-\sum_{j\in\Z/e\Z} N_j(s_i\mu)\alpha_j.
\end{cases}
$$
Thus $|s_i\mu|\le|\mu|$. In particular, $\mu$ is obtained from
$s_i\mu$ by adding all addable $i$-nodes by
Lemma \ref{Weyl group action} (2). Hence
$\lambda\subset s_i\mu\subset \mu$.

Next suppose that the leftmost $s_i$ appears in the subword for $u$.
Then $s_iu\le s_iv$ and the induction hypothesis implies
$s_i\lambda\subset s_i\mu$. Note that $s_iu<u$ and $s_iv<v$.
Thus, the same argument as above shows that
$\lambda$ and $\mu$ are obtained from $s_i\lambda$ and $s_i\mu$
by adding all addable $i$-nodes, respectively.
If an addable $i$-node of
$s_i\lambda$ is contained in $s_i\mu$, it is contained in $s_i\mu$ and
hence in $\mu$.
If an addable $i$-node of $s_i\lambda$ is not
contained in $s_i\mu$, then it is also an addable $i$-node of
$s_i\mu$. Thus, it is contained in $\mu$. We have proved
$\lambda\subset\mu$.

(2) We prove this by induction on $\ell(v)$ as above.
If $\lambda=\mu$ then there is nothing to prove.
So assume that $\lambda\ne\mu$.
Pick a removable node of the skew shape $\mu/\lambda$ and
denote its residue by $i$. As $\mu$ is an $e$-core,
$s_i\mu\subset\mu$ and $s_i\mu\ne\mu$. Thus we have
$s_iv<v$ by (1).

We show that we have either $\lambda\subset s_i\mu$ or
$s_i\lambda\subset s_i\mu$.
Suppose that $\lambda\not\subset s_i\mu$.
Then any node $x\in \lambda\setminus s_i\mu
\subset\mu\setminus s_i\mu$ is
a removable $i$-node of $\lambda$. Thus
$s_i\lambda\subset s_i\mu$ follows.
Hence, we consider these two cases.

First suppose that $\lambda\subset s_i\mu$.
Then the induction hypothesis
implies that $u\le s_iv$, as $u$ is a distinguished coset
representative. Thus $u\le v$.

Next suppose that $s_i\lambda\subset s_i\mu$ and
$\lambda\not\subset s_i\mu$. Then
$s_i\lambda\supset\lambda$ does not occur. Hence,
$s_i\lambda\subset\lambda$ and $s_i\lambda\ne\lambda$,
which implies $s_iu<u$ as before.

Write $s_iu=u't$, where $u'\in W/W_m$ and $t\in W_m$.
Let $u'=s_{i_1}\cdots s_{i_p}$ and
$t=s_{j_1}\cdots s_{j_q}$ be reduced
expressions of $u'$ and $t$ respectively. Then, as
$$
u=s_is_{i_1}\cdots s_{i_p}s_{j_1}\cdots s_{j_q}\;\;
\text{and}\;\;
\ell(u)=\ell(s_iu)+1=\ell(u')+\ell(t)+1=p+q+1,
$$
this is a reduced expression of $u$. Since $u$ is a
distinguished coset representative, we have $q=0$ and
$s_iu$ is distinguished. Now the induction hypothesis
implies $s_iu\le s_iv$. As $s_i(s_iv)>s_iv$,
we have $u\le v$ as desired.
\end{proof}

\begin{cor}
\label{cor to roof lemma}
Write $\roof(\lambda)=y_\lambda\emptyset_m$, for a unique
$y_\lambda\in W/W_m$. Then
$$
y_\lambda=\min\,\{y\in W \mid \lambda\in B_y(\Lambda_m)\}
$$
with respect to the Bruhat-Chevalley order.
\end{cor}
\begin{proof}
If $\lambda\in B_y(\Lambda_m)$ then Theorem \ref{KLMW}
shows that $\roof(\lambda)\subset y\emptyset_m$. Then
Proposition \ref{two orders} implies that $y_\lambda\le y$.
As $\roof(\lambda)\subset y_\lambda\emptyset_m$, we have
$\lambda\in B_{y_\lambda}(\Lambda_m)$ and $y_\lambda$ is
the unique minimal element of $\{y\in W \mid \lambda\in B_y(\Lambda_m)\}$.
\end{proof}

\section{Littelmann's path model}

Littelmann introduced a realization of Kashiwara crystals
in terms of $W$. \cite[\S1]{NS2} is a concise review of
the path model. The results of this section hold for a
general dominant integral weight, but we state them only for
$\Lambda_m$.

\begin{defn}
\label{defn of weight order}
Let $\mu\ne\nu\in W\Lambda_m$ be two weights.
If there exists a sequence of positive real roots
$\beta_1,\dots,\beta_r$ such that
$$
\langle s_{\beta_{j-1}}\cdots s_{\beta_1}\mu, h_{\beta_j}\rangle\in \Z_{<0},
$$
for $1\le j\le r$ and $\nu=s_{\beta_r}s_{\beta_{r-1}}\cdots s_{\beta_1}\mu$,
then we write $\mu>\nu$. Here, $h_{\beta_j}$ is the coroot of $\beta_j$.

Let $0<a<1$ be a rational number. A sequence
$$
\mu,\;s_{\beta_1}\mu,\;s_{\beta_2}s_{\beta_1}\mu,\;\cdots,\;
s_{\beta_r}s_{\beta_{r-1}}\cdots s_{\beta_1}\mu=\nu
$$
with $r$ maximal is called an {\bf $a$-chain} if
$$
\langle s_{\beta_{j-1}}\cdots s_{\beta_1}\mu, h_{\beta_j}\rangle\in a^{-1}\Z_{<0},
$$
for all $j$.
\end{defn}

If $\mu=y\Lambda_m$ and $\nu=w\Lambda_m$ for $y,w\in W/W_m$, then
$\mu>\nu$ is equivalent to $y>w$.

\begin{lemma}[{\cite[Lemma 4.1]{L2}}]
\label{weight order}
\item[(1)]
If $\mu\ge\nu$ is such that $\mu(h_i)<0$ and $\nu(h_i)\ge0$, then $s_i\mu\ge\nu$.
\item[(2)]
If $\mu\ge\nu$ is such that $\mu(h_i)\le0$ and $\nu(h_i)>0$, then $\mu\ge s_i\nu$.
\end{lemma}

Let $0=a_0<a_1<\cdots<a_s=1$ and $\nu_1,\dots,\nu_s\in W\Lambda_m$.
We consider a piecewise linear path $\pi(t)$, for $0\le t\le1$,
which takes values in the dual space of the Cartan subalgebra:
$$
\pi(t)|_{[a_{j-1},a_j]}=\sum_{k=1}^{j-1}(a_k-a_{k-1})\nu_k+(t-a_{j-1})\nu_j.
$$
In other words, we start with the origin, and change direction from
$\nu_j$ to $\nu_{j+1}$ at $t=a_j$, for $1\le j<s$.

\begin{defn}
The piecewise linear path $\pi(t)$ given by
$(\nu_1,\dots,\nu_s;a_0,\dots,a_s)$ as above,
is a {\bf Lakshmibai-Seshadri path}, if
the following hold for all $j$.
\begin{itemize}
\item[(i)]
$a_j$ is a rational number and $\nu_j>\nu_{j+1}$.
\item[(ii)]
There exists an $a_j$-chain for $\nu_j>\nu_{j+1}$.
\end{itemize}
We denote the set of Lakshmibai-Seshadri paths
by $\mathbb B(\Lambda_m)$.
\end{defn}

We call Lakshmibai-Seshadri paths {\bf LS paths} for short.

\begin{defn}
Let $\pi\in \mathbb B(\Lambda_m)$ be given by
$(\nu_1,\dots,\nu_s;a_0,\dots,a_s)$. We call $\nu_1$ the
{\bf initial direction} of $\pi$ and denote it by $i(\pi)$.
Similarly, we call $\nu_s$ the {\bf final direction} and
denote it by $f(\pi)$.
\end{defn}

\begin{defn}
We say that $\pi(t)$ satisfies the {\bf integrality condition} if
the minimum value of $\pi(t)(h_i)$ is an integer, for all $i$.
\end{defn}

\begin{lemma}[{\cite[Lemma 4.5(d)]{L2}}]
\label{integrality condition}
The LS-paths satisfy the integrality condition.
\end{lemma}

Define $Q=\min\{\pi(t)(h_i) \mid 0\le t\le1\}$.
We shall define the operators $\e_i$ and $\f_i$ on
$\mathbb B(\Lambda_m)\sqcup\{0\}$.
First of all, we set $\e_i\pi=0$ if $Q>-1$, and
$\f_i\pi=0$, if $Q>\pi(1)(h_i)-1$.

Suppose that $Q\le -1$. Then define
\begin{equation*}
\begin{split}
t_1&=\min\{t\in [0,1] \mid \pi(t)(h_i)=Q\}\\
t_0&=\max\{t\in [0,t_1] \mid \pi(t)(h_i)|_{[0,t]}\ge Q+1\}
\end{split}
\end{equation*}
and reflect the path $\pi(t)$ for the interval $[t_0,t_1]$ to define:
$$
(\e_i\pi)(t)=\begin{cases}
\pi(t) \quad(0\le t\le t_0)\\
s_i(\pi(t)-\pi(t_0))+\pi(t_0) \quad(t_0\le t\le t_1)\\
\pi(t)+\alpha_i\quad(t_1\le t\le 1)\end{cases}
$$

Suppose that $Q\le\pi(1)(h_i)-1$. Then define
\begin{equation*}
\begin{split}
t_0&=\max\{t\in [0,1] \mid \pi(t)(h_i)=Q\}\\
t_1&=\min\{t\in [t_0,1] \mid \pi(t)(h_i)|_{[t,1]}\ge Q+1\}
\end{split}
\end{equation*}
and define:
$$
(\f_i\pi)(t)=\begin{cases}
\pi(t) \quad(0\le t\le t_0)\\
s_i(\pi(t)-\pi(t_0))+\pi(t_0) \quad(t_0\le t\le t_1)\\
\pi(t)-\alpha_i\quad(t_1\le t\le 1)\end{cases}
$$
We then define $wt(\pi)=\pi(1)$ and
$$
\epsilon_i(\pi)=-Q, \quad \varphi_i(\pi)=\pi(1)(h_i)-Q.
$$

Then, by \cite[Corollary 6.4.27]{Jo} or \cite[Theorem 4.1]{K5}, the set
$\mathbb B(\Lambda_m)$ with the additional
data $wt$, $\epsilon_i$, $\varphi_i$, $\e_i$ and $\f_i$
is a realization of the crystal $B(\Lambda_m)$. The isomorphism
of the two realizations, one by $e$-restricted partitions, the other
by the LS-paths, is unique. Thus, we identify the two realizations and
sometimes write $\lambda=(\nu_1,\dots,\nu_s;a_0,\dots,a_s)$, for
an $e$-restricted partition $\lambda$. We denote $\nu_1$ and
$\nu_s$ by $i(\lambda)$ and $f(\lambda)$ respectively.

The following is one of the key results we use in this paper.

\begin{thm}[{\cite[Theorem 10.1]{L3}}]
\label{Littelmann's theorem}
Let
$$
\pi=\pi^{(1)}\otimes\cdots\otimes \pi^{(r)}\in
B(\Lambda_{m_1})\otimes\cdots B(\Lambda_{m_r}).
$$
Then $\pi$
belongs to $B(\Lambda_{m_1}+\cdots\Lambda_{m_r})$
if and only if there exists a sequence
$$
w_1^{(1)}\ge\cdots\ge w_{N_1}^{(1)}\ge w_1^{(2)}\ge\cdots\ge
w_{N_2}^{(2)}\ge\cdots\cdots\ge w_{N_r}^{(r)}
$$
in $W$ such that
$$
\pi^{(k)}=(w_1^{(k)}\Lambda_{m_k},\dots,w_{N_k}^{(k)}\Lambda_{m_k};
a_0^{(k)},\dots,a_{N_k}^{(k)}),
$$
for $1\le k\le r$.
\end{thm}

Recall that $w_0$ is the longest element of $W_0$.

\begin{cor}
\label{the result to use}
Let $\pi=\pi^{(1)}\otimes\cdots\otimes \pi^{(r)}\in
B(\Lambda_0)^{\otimes d}\otimes B(\Lambda_m)^{\otimes r-d}$, and
write $w\Lambda_0=f(\pi^{(d)})$ and
$w'\Lambda_m=i(\pi^{(d+1)})$, for $w\in W/W_0$ and
$w'\in W/W_m$ respectively.
Then $\pi$
belongs to $B(d\Lambda_0+(r-d)\Lambda_m)$ if and only if
\begin{itemize}
\item[(a)]
$f(\pi^{(k)})\ge i(\pi^{(k+1)})$, for $1\le k<d$,
\item[(b)]
$ww_0\ge w'$,
\item[(c)]
$f(\pi^{(k)})\ge i(\pi^{(k+1)})$, for $d<k<r$.
\end{itemize}
\end{cor}
\begin{proof}
If $\pi$ belongs to $B(d\Lambda_0+(r-d)\Lambda_m)$, then
Theorem \ref{Littelmann's theorem} gives a non-increasing sequence in $W$,
which implies that conditions (a) to (c) hold.

Suppose that conditions (a) to (c) hold.
Consider the elements $w\in W/W_0$
such that $w\Lambda_0$ appears as one of the direction vectors
of $\pi^{(1)},\dots,\pi^{(d)}$. Multiplying them with $w_0$ simultaneously,
we can find the desired sequence in $W$.
Thus, Theorem \ref{Littelmann's theorem}
implies that $\pi$ belongs to $B(d\Lambda_0+(r-d)\Lambda_m)$.
\end{proof}

Our purpose is to interpret this result in terms of Young diagrams.
To achieve this goal, we first have to find which partitions correspond
to $f(\pi)$ and $i(\pi)$ when $\pi$ corresponds to a partition $\lambda$.

For this, we need to use the approach to path models in \cite{K5} and
\cite[chapter 8]{K1}.

\begin{defn}
Let $B$ and $B'$ be crystals. A map $\psi:B\rightarrow B'$ is
called a {\bf crystal morphism of amplitude $h$} if
\begin{itemize}
\item[(i)]
$wt(\psi(b))=h wt(b)$, $\epsilon_i(\psi(b))=h\epsilon_i(b)$ and
$\varphi_i(\psi(b))=h\varphi_i(b)$,
\item[(ii)]
$\psi(\e_ib)=\e_i^h\psi(b)$ and
$\psi(\f_ib)=\f_i^h\psi(b)$, for all $b\in B$.
\end{itemize}
\end{defn}

\begin{defn}
\begin{itemize}
\item[(1)]
$U_v^-(\mathfrak g)$ is a module
over the Kashiwara algebra, which defines a crystal.
This is the crystal $B(\infty)$ and
$$
\epsilon_i(b)=\max\{k\in \mathbb Z_{\ge0} \mid \e_i^kb\ne0\},\;\;
\varphi_i(b)=\epsilon_i(b)+wt(b)(h_i).
$$
\item[(2)]
Define, for $a\in\mathbb Z$,
$$
wt(a)=a\alpha_i,\;\;
\epsilon_j(a)=\begin{cases} -a\quad(j=i)\\
                            -\infty\quad(j\ne i)\end{cases},\;\;
\varphi_j(a)=\begin{cases}  a\quad(j=i)\\
                            -\infty\quad(j\ne i)\end{cases}
$$
and
$$
\e_j(a)=\begin{cases} a+1\quad(j=i)\\
                            0\quad(j\ne i)\end{cases},\;\;
\f_j(a)=\begin{cases}  a-1\quad(j=i)\\
                             0\quad(j\ne i)\end{cases}.
$$
Then $\mathbb Z$ becomes a crystal. This is the crystal $B_i$.
\item[(3)]
Let $\Lambda$ be a weight, and define
$$
wt(t_\Lambda)=\Lambda,\;\;\epsilon_i(t_\Lambda)=\varphi_i(t_\Lambda)=-\infty,
\;\;
\e_it_\Lambda=\f_it_\Lambda=0.
$$
Then $\{t_\Lambda\}$ is the crystal $T_\Lambda$.
\item[(4)]
Define another crystal structure on the underlying set of $B(\infty)$ by
redefining $(wt,\epsilon_i,\varphi_i,\e_i,\f_i)$ by
$$
wt^{new}=-wt^{old},\;\epsilon_i^{new}=\varphi_i^{old},\;
\varphi_i^{new}=\epsilon_i^{old},\;
\e_i^{new}=\f_i^{old},\;\f_i^{new}=\e_i^{old}.
$$
This crystal is denoted by
$B(-\infty)$. It may be considered as the crystal arising from
the positive part $U_v^+(\mathfrak g)$. We have
$$
\varphi_i(b)=\max\{k\in \mathbb Z_{\ge0} \mid \f_i^kb\ne0\},\;\;
\epsilon_i(b)=\varphi_i(b)-wt(b)(h_i).
$$
\end{itemize}
\end{defn}

We fix an infinite sequence
$\underline i=(\cdots,i_k,\cdots,i_2,i_1)$ such
that $i_k\ne i_{k+1}$, for all $k$, and that $i$ appears
infinitely many times in the sequence, for all $i$.
Then we can realize $B(\infty)$ as a subcrystal of
$\mathbb Z^\infty_{\underline i}=\cdots\otimes B_{i_k}\otimes\cdots
\otimes B_{i_2}\otimes B_{i_1}$ \cite[Theorem 2.2.1]{K3}.
This is the {\bf Kashiwara embedding} and the
{\bf polyhedral realization} associated with $\underline i$.

\begin{prop}[{\cite[Proposition 8.1.3]{K1}}]
\label{Sh for Binfty}
For all $h\in \mathbb N$,
there exists a unique crystal morphism
$S_h: B(\infty)\rightarrow B(\infty)$ of amplitude $h$.
$S_h$ is an injective map.
In any polyhedral realization, we have
$$
S_h(\cdots,a_k,\cdots,a_2,a_1)=(\cdots,ha_k,\cdots,ha_2,ha_1).
$$
\end{prop}

In fact, this is proved by defining $S_h$ by the above formula
in the polyhedral realization of $B(\infty)$
and showing that this is a crystal morphism of amplitude $h$.
Define $S_h: B(\infty)\otimes T_\Lambda \rightarrow
B(\infty)\otimes T_{h\Lambda}$ by
$b\otimes t_\Lambda\mapsto S_h(b)\otimes t_{h\Lambda}$.
This is again a crystal morphism of amplitude $h$.

\begin{prop}[{\cite[Corollary 8.1.5]{K1}}]
\label{Sh for hw crystal}
Let $\Lambda$ be dominant integral. Then
there exists a unique crystal morphism
$S_h: B(\Lambda)\rightarrow B(h\Lambda)$ of amplitude $h$,
for all $h\in \mathbb N$. Further, we have the
following commutative diagram.
\[
\begin{array}{ccc}
B(\Lambda) & \stackrel{S_h}{\longrightarrow} & B(h\Lambda) \\
\cap & & \cap \\
B(\infty)\otimes T_{\Lambda} & \stackrel{S_h}{\longrightarrow} &
B(\infty)\otimes T_{h\Lambda}
\end{array}
\]
\end{prop}

Let $\lambda\in B(\Lambda_m)$.
Using the canonical embedding $B(h\Lambda_m)\subset B(\Lambda_m)^{\otimes h}$,
we can write
$$
S_h(\lambda)=\lambda^{(1)}\otimes\cdots\otimes \lambda^{(h)}.
$$
We denote
$$
S_h(\lambda)^{1/h}=
{\lambda^{(1)}}^{\otimes 1/h}\otimes\cdots\otimes{\lambda^{(h)}}^{\otimes 1/h},
$$
and replace $(\mu^{\otimes 1/h})^{\otimes k}$ with
$\mu^{\otimes k/h}$, for any $\mu$ that appears in
$\lambda^{(1)},\dots, \lambda^{(h)}$. In this way, we may write
$$
S_h(\lambda)^{1/h}={\nu_1}^{\otimes a_1}\otimes{\nu_2}^{\otimes(a_2-a_1)}
\otimes\cdots\otimes{\nu_s}^{\otimes(1-a_{s-1})},
$$
where $a_0=0<a_1<\cdots<a_s=1$ are rational integers and
$\nu_1,\dots,\nu_s$ are pairwise distinct $e$-restricted partitions.
Then the following theorem holds.

\begin{thm}[{\cite[Proposition 8.3.2]{K1}}]
If $h$ is sufficiently divisible then
\item[(1)]
$\nu_j=w_j\emptyset_m$, for a unique $w_j\in W/W_m$.
\item[(2)]
$a_j$ and $\nu_j$ all stabilize.
\end{thm}

\begin{thm}[{\cite[Proof of Theorem 8.2.3]{K1}}]
\label{stability}
Given sufficiently divisible $h$, we write
$$
S_h(\lambda)^{1/h}={\nu_1}^{\otimes a_1}\otimes{\nu_2}^{\otimes(a_2-a_1)}
\otimes\cdots\otimes{\nu_s}^{\otimes(1-a_{s-1})}
$$
as above, and define $\pi_\lambda$ to be the path given by
$(wt(\nu_1),\dots,wt(\nu_s);a_0,\dots,a_s)$. Then
$wt(\nu_j)=w_j\Lambda_m$, for $1\le j\le s$, and the following hold.
\item[(1)]
$\pi_\lambda$ is a LS-path.
\item[(2)]
The map $B(\Lambda_m)\rightarrow \mathbb B(\Lambda_m)$ defined by
$\lambda\mapsto \pi_\lambda$ is an isomorphism of crystals.
\end{thm}

The proof of \cite[Proposition 8.3.2]{K1} also gives a very explicit
inductive algorithm to compute the $e$-cores $\nu_j$ as follows.

Recall that the tensor product rule for $B(\Lambda_m)^{\otimes r}$ is
given by the following rule:

Let $\lambda^{(1)}\otimes\cdots\otimes\lambda^{(r)}\in
B(\Lambda_m)^{\otimes r}$. Then, starting with $\lambda^{(r)}$,
we read addable and removable
$i$-nodes of each $\lambda^{(k)}$ from the first row to the last row,
for $k=r,r-1,\dots,1$ succesively. We then apply the $RA$-deletion
to the resulting sequence of dots, $A$'s and $R$'s.

\begin{lemma}
\label{rule to compute paths}
Suppose that
$$
S_h(\lambda)^{1/h}={\nu_1}^{\otimes a_1}\otimes{\nu_2}^{\otimes(a_2-a_1)}
\otimes\cdots\otimes{\nu_s}^{\otimes(1-a_{s-1})}.
$$
and that $\f_i\lambda\ne0$. Then $(a_{k+1}-a_k)h$ are positive integers
and we write
$$
{\nu_1}^{\otimes a_1h}\otimes{\nu_2}^{\otimes(a_2-a_1)h}
\otimes\cdots\otimes{\nu_s}^{\otimes(1-a_{s-1})h}=
\mu_1\otimes\cdots\otimes \mu_h.
$$
Then we may write
$$
\f_i^h(\mu_1\otimes\cdots\otimes \mu_h)=
\f_i^{c_1}\mu_1\otimes\f_i^{c_2}{\mu_2}\otimes\cdots\otimes
\f_i^{c_h}\mu_h,
$$
for some non-negative integers $c_j$ such that $\sum_{j=1}^h c_j=h$.
Then, for some multiple $h'$ of $h$, we have
\begin{multline*}
S_{h'}(\f_i\lambda)^{1/{h'}}=
\bigl(({s_i\mu_1}^{\otimes(c_1/\varphi_i(\mu_1))h'/h}
\otimes{\mu_1}^{\otimes(1-c_1/\varphi_i(\mu_1))h'/h})
\otimes\cdots\\
\cdots\cdots\otimes
({s_i\mu_h}^{\otimes(c_h/\varphi_i(\mu_h))h'/h}
\otimes{\mu_h}^{\otimes(1-c_h/\varphi_i(\mu_h))h'/h})\bigr)^{1/h'}.
\end{multline*}
\end{lemma}

\begin{Ex}
Let $m=0$ and $e=3$. Then $\lambda=(3,1^2)$ is an $e$-core.
Thus $S_h(\lambda)^{1/h}=(3,1^2)$, for all $h$.
Consider $\lambda'=(3,1^3)=\f_0\lambda$.
Then, $\varphi_0(\lambda)=3$ and we have, for $h$ which is divisible by $3$,
$$
S_{h}(\lambda')^{1/h}=(4,2,1^2)^{\otimes 1/3}\otimes (3,1^2)^{\otimes 2/3}.
$$
\end{Ex}

\begin{defn}
Suppose that
$$
S_h(\lambda)^{1/h}={\nu_1}^{\otimes a_1}\otimes{\nu_2}^{\otimes(a_2-a_1)}
\otimes\cdots\otimes{\nu_s}^{\otimes(1-a_{s-1})},
$$
for sufficiently divisible $h$. Then we call
$\nu_1$ the {\bf ceiling} of $\lambda$
and denote it by $\ceil(\lambda)$. Similarly, we call $\nu_s$ the
{\bf floor} of $\lambda$ and denote it by $\floor(\lambda)$.
\end{defn}

We have $wt(\ceil(\lambda))=i(\lambda)$ and $wt(\floor(\lambda))=f(\lambda)$
by the definitions and Theorem \ref{stability}(2).

\begin{Ex}
Let $m=0$, $e=3$, and $\lambda=(2^2,1)$. Then, for $h$ which is divisible by 6,
$$
S_h(\lambda)^{1/h}=(5,3,1)^{\otimes 1/3}\otimes (4,2)^{\otimes 1/6}\otimes
(2)^{\otimes 1/2}.
$$
Thus, $\ceil(\lambda)=(5,3,1)$ and $\floor(\lambda)=(2)$.
\end{Ex}

Note that in this paper we define $\ceil(\lambda)$ in a different
manner than \cite{KLMW}, because
we follow a slightly different line of proof.
That the two definitions give the same $e$-core
follows from Theorem \ref{ceiling theorem} and
Corollary \ref{roof is ceil} below, which prove Theorem \ref{KLMW}, and
Corollary \ref{cor to roof lemma}.

Fayers pointed out that $\ceil(\lambda)$ and $\floor(\lambda)$
behave well under the Mullineux map. Let us review
the Mullineux map quickly.
Let $\H_n(q)$ be the Hecke algebra of type $A$.
This is the $F$-algebra generated by $T_1,\dots,T_{n-1}$
subject to the quadratic relations $(T_i-q)(T_i+1)=0$ and
the type $A$ braid relations. Let
$\tau$ be the involution of $\H_n(q)$
defined by $T_i\mapsto qT_i^{-1}$.
The simple $\H_n(q)$-modules are
$\{D^\lambda \mid \text{$\lambda$ is $e$-restricted}\}$.
Then the Mullineux map is defined by $D^{m(\lambda)}=(D^\lambda)^\tau$.
In \cite[Theorem 7.1]{LLT}, it is observed that the description
of the Mullineux map obtained by Brundan and Kleshchev may be
expressed in terms of the crystal $B(\Lambda_0)$.
Shifting the residues, the Mullineux map may be described by
$B(\Lambda_m)$ also.

\begin{prop}
Suppose that $\lambda\in B(\Lambda_m)$ is such that
$\lambda=\f_{m+i_1}\cdots \f_{m+i_n}\emptyset$.
Then we have
$m(\lambda)=\f_{m-i_1}\cdots \f_{m-i_n}\emptyset$.
\end{prop}

\begin{cor}
$\epsilon_{m+i}(\lambda)=\epsilon_{m-i}(m(\lambda))$ and
$\varphi_{m+i}(\lambda)=\varphi_{m-i}(m(\lambda))$.
\end{cor}
\begin{proof}
If $\e_{m+i}^k\lambda\ne0$ then
$\e_{m-i}^km(\lambda)=m(\e_{m+i}^k\lambda)\ne0$.
Thus we have
$\epsilon_{m+i}(\lambda)\le\epsilon_{m-i}(m(\lambda))$.
Similarly, we have
$\varphi_{m+i}(\lambda)\le\varphi_{m-i}(m(\lambda))$.
Then we also have $\epsilon_{m+i}(m(\lambda))\le\epsilon_{m-i}(\lambda)$
and $\varphi_{m+i}(m(\lambda))\le\varphi_{m-i}(\lambda)$.
Hence the equalities hold.
\end{proof}

\begin{prop}
\label{ceil and floor under Mullineux}
Let $\lambda\in B(\Lambda_m)$. Then
$\ceil(m(\lambda))$ and $\floor(m(\lambda))$ are the
conjugate partitions of $\ceil(\lambda)$ and
$\floor(\lambda)$ respectively.
\end{prop}
\begin{proof}
We may assume that $m=0$ without loss of generality.
We prove by induction on $|\lambda|$ that if
$S_h(\lambda)=\nu_1\otimes\cdots\otimes \nu_h$
for sufficiently divisible $h$, then
$$
S_h(m(\lambda))=\nu_1'\otimes\cdots\otimes \nu_h',
$$
where $\nu_k'$ is the conjugate partition of $\nu_k$,
for all $k$.

If $|\lambda|=0$ there is nothing to prove.
Assume that the assertion holds for $\lambda$ and that
$\f_i\lambda\ne0$. Note that $\nu_k$ are $e$-cores
and thus $\nu_k$ has removable $i$-nodes only, or
addable $i$-nodes only.
If $\nu_k$ has $n_{i,k}$ removable $i$-nodes
then $\nu_k'$ has $n_{i,k}$ removable $(-i)$-nodes,
and similarly, if $\nu_k$ has $n_{i,k}$ addable $i$-nodes
then $\nu_k'$ has $n_{i,k}$ addable $(-i)$-nodes.
This implies that if
$$
S_h(\f_i\lambda)
=\f_i^{c_1}\nu_1\otimes\cdots\otimes\f_i^{c_h}\nu_h,
$$
then
$$
S_h(m(\f_i\lambda))
=\f_{-i}^{c_1}\nu_1'\otimes\cdots\otimes\f_{-i}^{c_h}\nu_h'.
$$
Now, to obtain $S_{h'}(\f_i\lambda)$, for sufficiently divisible
$h'$, we replace $\f_i^{c_k}\nu_k$ with
$$
(s_i\nu_k)^{\otimes(c_k/\varphi_i(\nu_k))h'/h}\otimes
\nu_k^{\otimes(1-c_k/\varphi_i(\nu_k))h'/h},
$$
for all $k$.
As $\varphi_i(\nu_k)=\varphi_{-i}(\nu_k')$, the assertion
holds for $\f_i\lambda$.
\end{proof}

Recall that the Mullineax map is given by conjugation of a partition
when $\H_n(q)$ is semisimple. The proof of
Proposition \ref{ceil and floor under Mullineux} shows that
the Mullineax map is always given by conjugation, if we work
in the right model -- the path model.

The descriptions of $\ceil(\lambda)$ and $\floor(\lambda)$
are a crucial part of our main results. In the case when $e=2$,
we have closed formulas for them.

\begin{prop}
\label{floor and ceil for e=2}
Assume that $e=2$ and that $\lambda\in B(\Lambda_m)$.
Let $a(\lambda)$ be the length of the first row, and let
$\ell(\lambda)$ be the length of the first column.
Then
$$
\ceil(\lambda)=(\ell(\lambda),\ell(\lambda)-1,\dots,1),\quad
\floor(\lambda)=(a(\lambda),a(\lambda)-1,\dots,1).
$$
\end{prop}
\begin{proof}
We prove both formulas by induction on the size of $\lambda$.
As $\lambda$ is $2$-restricted, the last node of the first column
is removable. Let $i$ be its residue. Let
$\mu=\e_i^{max}\lambda=\e_i^t\lambda$.
Then by the induction hypothesis, we have
$$
\ceil(\mu)=(\ell(\lambda)-1,\ell(\lambda)-2,\dots,1).
$$
Observe that there exists an addable normal
$i$-node on the first column of
$\ceil(\mu)$. Thus
all normal $i$-nodes are addable and the addable $i$-node on
the first column of
$\ceil(\mu)$ is the first addable $i$-node to be changed into
a removable $i$-node when $\f_i^t$
is applied to $\mu$. Thus, Lemma \ref{rule to compute paths} implies that
$$
\ceil(\lambda)=s_i(\ell(\lambda)-1,\ell(\lambda)-2,\dots,1)=
(\ell(\lambda),\ell(\lambda)-1,\dots,1).
$$
Hence, the formula for $\ceil(\lambda)$ is proved.

Next assume that the formula for $\floor(\lambda)$ is already proved.
Consider the addable node on the first row. Let $i$ be its residue.
Then, this addable node is a normal $i$-node. The induction
hypothesis implies that $\floor(\lambda)$ has addable normal $i$-nodes.
First suppose that $\varphi_i(\lambda)>1$. Then $\f_i\lambda$ differs
from $\lambda$ at some node which lies in the second row or below.
Thus $a(\f_i\lambda)=a(\lambda)$.
Let $h$ be sufficiently divisible. Then $\varphi_i(\lambda)>1$
implies that we do not apply
$\f_i^{max}=\f_i^{h\varphi_i(\lambda)}$ to $S_h(\lambda)$
when computing $S_h(\f_i\lambda)$. Since the addable $i$-nodes
of $\floor(\lambda)$ are the last addable normal $i$-nodes to be
changed into removable $i$-nodes, that we do not apply
$\f_i^{max}$ to $S_h(\lambda)$ implies that
$\floor(\f_i\lambda)=\floor(\lambda)$
by Lemma \ref{rule to compute paths}. Hence we have proved the formula
in this case. Second suppose that $\varphi_i(\lambda)=1$. Then
$\f_i\lambda$ differs from $\lambda$ at
the addable $i$-node on the first row. Thus
$a(\f_i\lambda)=a(\lambda)+1$. Let $h$ be sufficiently divisible.
Then $\varphi_i(\lambda)=1$ implies that we apply
$\f_i^{max}$ to $S_h(\lambda)$ when computing $S_h(\f_i\lambda)$.
As the addable $i$-nodes of $\floor(\lambda)$ are all normal, this
implies that
$$
\floor(\f_i\lambda)=s_i\floor(\lambda)=
s_i(a(\lambda),a(\lambda)-1,\dots,1)=
(a(\lambda)+1,a(\lambda),\dots,1).
$$
Hence, we have proved the formula in this case also.
\end{proof}

\section{Description of Demazure crystals}

\begin{lemma}
\label{lemma for i and f}
\item[(1)]
Suppose that $i(\pi)(h_i)<0$. Then
\begin{itemize}
\item[(a)]
$\e_i\pi\ne0$.
\item[(b)]
$i(\e_i^{max}\pi)=s_ii(\pi)<i(\pi)$.
\end{itemize}
\item[(2)]
Suppose that $f(\pi)(h_i)>0$. Then
\begin{itemize}
\item[(a)]
$\f_i\pi\ne0$.
\item[(b)]
$f(\f_i^{max}\pi)=s_if(\pi)>f(\pi)$.
\end{itemize}
\end{lemma}
\begin{proof}
(1) (a) $i(\pi)(h_i)<0$ implies that $Q=\min\{\pi(t)(h_i) \mid 0\le t\le 1\}<0$.
Thus, Lemma \ref{integrality condition} implies that $Q\le -1$ and
$\e_i\pi\ne0$.

(b) By (a), $i(\e_i^{max}\pi)(h_i)\ge0$. Then
$$
s_ii(\e_i^{max}\pi)\ge i(\e_i^{max}\pi).
$$
On the other hand, we have either $i(\e_i^{max}\pi)=i(\pi)$ or
$s_ii(\pi)$. As $i(\pi)(h_i)<0$ we have $s_ii(\pi)<i(\pi)$ and
we have the result.

(2) (a) $f(\pi)(h_i)>0$ implies that $Q\le \pi(1)(h_i)-1$ by the integrality
condition again. Thus $\f_i\pi\ne0$.

(b) The proof is similar to that of (1).
\end{proof}

We thank Dr. Sagaki for showing us the proof of the following theorem.
The proof for the first equality
works for dominant integral weights in general.

\begin{thm}[{\cite[Theorem 2]{L4}}]
\label{ceiling theorem}
Suppose $y\in W/W_m$. Then
$$
B_y(\Lambda_m)=\{\lambda\in B(\Lambda_m) \mid i(\lambda)\le y\Lambda_m\}
=\{\lambda\in B(\Lambda_m) \mid \ceil(\lambda)\subset y\emptyset_m\}.
$$
\end{thm}
\begin{proof}
We only have to prove the first equality. The second equality follows from the
remark at the end of Definition \ref{defn of weight order} and Proposition
\ref{two orders}. We prove
$$
B_y(\Lambda_m)\supset\{\lambda\in B(\Lambda_m) \mid i(\lambda)\le y\Lambda_m\}
$$
by induction on $\ell(y)$. If $y=1$ then $B_y(\Lambda_m)=\{\emptyset_m\}$
and $i(\lambda)\le\Lambda_m$ implies that $\lambda=\emptyset_m$.
Thus the statement is true.

Let $y=s_is_{i_2}\cdots s_{i_\ell}$ be a reduced expression. First
we remark that $s_iy\Lambda_m=y\Lambda_m$ is impossible: otherwise
$s_iy=yu$, for some $u\in W_m$, which implies
$\ell(s_iy)<\ell(y)\le\ell(yu)=\ell(s_iy)$, a contradiction.
Thus $y\Lambda_m(h_i)<0$ and $s_iy\in W/W_m$.

Assume that $i(\lambda)\le y\Lambda_m$.
If $i(\lambda)(h_i)\ge0$ then Lemma \ref{weight order} (1)
implies that $i(\lambda)\le s_iy\Lambda_m$. Hence, by the induction
hypothesis and the fact that $s_iy<y$, Proposition \ref{basic properties} (4)
implies that
$\lambda\in B_{s_iy}(\Lambda_m)\subset B_y(\Lambda_m)$.
If $i(\lambda)(h_i)<0$ then
Lemma \ref{lemma for i and f} (1) implies that
$i(\e_i^{max}\lambda)=s_ii(\lambda)<i(\lambda)$.
Since $s_iy\Lambda_m<y\Lambda_m$ and $i(\lambda)\le y\Lambda_m$,
we have $s_ii(\lambda)\le s_iy\Lambda_m$. The induction hypothesis
then implies that $\e_i^{max}\lambda\in B_{s_iy}(\Lambda_m)$.
Now, $\lambda\in B_y(\Lambda_m)$ by
Proposition \ref{basic properties} (2).

The opposite inclusion is easy to prove. In fact, if
$\lambda\in B_y(\Lambda_m)$ then we may write
$\lambda=\f_i^{a_1}\cdots \f_{i_\ell}^{a_\ell}\emptyset_m$.
We apply $\f_i^{a_1}\cdots \f_{i_\ell}^{a_\ell}$ to the path
associated with $\emptyset_m$ and we obtain
$i(\lambda)=y'\Lambda_m$, for some $y'\le y$. Hence
$i(\lambda)\le y\Lambda_m$.
\end{proof}

Theorem \ref{KLMW} is proved by the theorem below,
which is called the \lq\lq roof lemma\rq\rq in \cite{KLMW}.

\begin{thm}[{\cite[Lemma 3.3]{KLMW}}]
\label{roof lemma}
Let $\lambda\in B(\Lambda_m)$. Denote the residue of the
removable node on the last row by $i$. Then
$$
\roof(\lambda)\supset\roof(\e_i^{max}\lambda)=s_i\roof(\lambda).
$$
\end{thm}

\begin{cor}
\label{roof is ceil}
$\roof(\lambda)=\ceil(\lambda)$.
\end{cor}
\begin{proof}
Note that Lemma \ref{rule to compute paths} implies that
$$
\ceil(\lambda)\supset\ceil(\e_i^{max}\lambda)=s_i\ceil(\lambda).
$$
Thus induction on the size of
$\lambda$ proves the result.
\end{proof}

Hence, Theorem \ref{KLMW} follows from Theorem \ref{ceiling theorem} and
Corollary \ref{roof is ceil}.

Our aim is to prove a similar result for $B^w(\Lambda_m)$. For
this, we need a desciption of $B^w(\Lambda_m)$ which is similar to
the description of $B_y(\Lambda_m)$ in Theorem \ref{ceiling
theorem}. Fortunately, such a result exists. We thank Kashiwara
and Sagaki, who kindly showed us the result. Here we follow
Kashiwara's argument. As there exists no written proof, he
permitted us to include his argument here.

Before explaining the result, which is Theorem \ref{floor theorem} below,
we recall more results from the crystal theory.

Denote the canonical basis of $U_v^-(\mathfrak g)$
by $\{G_v(b) \mid b\in B(\infty)\}$.
Let $\Lambda$ be a dominant integral weight.
Then the irreducible highest weight module with highest weight $\Lambda$
has the basis $\{G_v(b)u_\Lambda \mid b\in B(\infty)\}\setminus\{0\}$, where
$u_\Lambda$ is a highest weight vector. When $wt(G_v(b)u_\Lambda)=w\Lambda$,
for $w\in W$, we denote $G_v(b)u_\Lambda$ by $u_{w\Lambda}$.

\begin{prop}[{\cite[Proposition 4.1]{K3}}]
\label{down to h/w crystal}
\begin{itemize}
\item[(1)]
Let $G_v(b)u_{w\Lambda}\ne0$, for $b\in B(\infty)$. Then
$G_v(b)u_{w\Lambda}=G_v(b')u_\Lambda$, for some $b'\in B(\infty)$.
\item[(2)]
If $G_v(b)u_{w\Lambda}=G_v(b')u_{w\Lambda}\ne0$, for
$b, b'\in B(\infty)$, then $b=b'$.
\end{itemize}
\end{prop}

Let $(L_v(\Lambda),B(\Lambda))$ be the crystal basis of
the integrable highest weight module $U_v(\mathfrak g)u_\Lambda$.
We have $\{G_v(b) \mid b\in B(\Lambda)\}=
\{G_v(b)u_\Lambda \mid b\in B(\infty)\}\setminus\{0\}$.
Then the following holds by \cite[(4.1)]{K3}.

\begin{lemma}
\label{third main theorem for extremal weight}
$\{G_v(b) \mid b\in B^w(\Lambda)\}$
is a basis of $U_v^-(\mathfrak g)u_{w\Lambda}$.
\end{lemma}

\begin{lemma}
\label{equality}
$\{G_v(b) \mid b\in B^w(\Lambda)\}= \{G_v(b)u_{w\Lambda} \mid b\in
B(\infty)\}\setminus\{0\}$.
\end{lemma}
\begin{proof}
Suppose that $b\in B^w(\Lambda)$.
Then, Lemma \ref{third main theorem for extremal weight} implies that
we may write
$$
G_v(b)=\sum_{b'\in B(\infty)}f_{b'}G_v(b')u_{w\Lambda},
$$
for some $f_{b'}\in\mathbb Q(v)$. Then
Proposition \ref{down to h/w crystal}(1) asserts that each nonzero
$G_v(b')u_{w\Lambda}$ is of the form $G_v(b")u_\Lambda$, for some
$b"\in B(\infty)$. Therefore,
$G_v(b)=G_v(b')u_{w\Lambda}$, for
some $b'\in B(\infty)$.

Suppose that $G_v(b)u_{w\Lambda}\ne0$, for some $b\in B(\infty)$.
Then $G_v(b)u_{w\Lambda}=G_v(b')$, for some $b'\in B(\Lambda)$, by
Proposition \ref{down to h/w crystal}(1) again. Since
$$
G_v(b')=G_v(b)u_{w\Lambda}\in U_v^-(\mathfrak g)u_{w\Lambda},
$$
Lemma \ref{third main theorem for extremal weight} implies that
$b'\in B^w(\Lambda)$.
\end{proof}

\begin{prop}
Assume that there exists a sequence
$$
w_1\ge w_2\ge\cdots\ge w_h=w.
$$
Then
$$
u_{w_1\Lambda}\otimes\cdots\otimes u_{w_h\Lambda}+vL_v(\Lambda)^{\otimes h}
\in B^w(h\Lambda)\subset B(h\Lambda)\subset B(\Lambda)^{\otimes h}.
$$
\end{prop}
\begin{proof}
The proof is by induction on $h$. When $h=1$,
$u_{w\Lambda}+vL_v(\Lambda)\in B^w(\Lambda)$ by Lemma
\ref{equality}, so there is nothing to prove. Suppose that $h>1$.
By the induction hypothesis, we may assume that
$$
u_{w_1\Lambda}\otimes\cdots\otimes u_{w_{h-1}\Lambda}+
vL_v(\Lambda)^{\otimes(h-1)}\in B^{w_{h-1}}((h-1)\Lambda)
\subset B^w((h-1)\Lambda).
$$
This and Lemma \ref{equality} imply that
there exists $b\in B(\infty)$ such that
$$
u_{w_1\Lambda}\otimes\cdots\otimes u_{w_{h-1}\Lambda}+
vL_v(\Lambda)^{\otimes(h-1)}=
G_v(b)u_{(h-1)w\Lambda}+
vL_v(\Lambda)^{\otimes(h-1)}.
$$
Consider $G_v(b)(u_{hw\Lambda})$. As $G_v(b)(u_{hw\Lambda})=0$ or
$G_v(b)(u_{hw\Lambda})=G_v(b')u_{h\Lambda}$, for some $b'\in B(\infty)$,
by Proposition \ref{down to h/w crystal}(1), we have
$$
G_v(b)(u_{(h-1)w\Lambda}\otimes u_{w\Lambda})=G_v(b)(u_{hw\Lambda})\in
L_v(\Lambda)^{\otimes h}.
$$
If we view $L_v(\Lambda)^{\otimes h}$ as a $\mathfrak g^{\otimes h}$-crystal
lattice and consider its weight decomposition,
$(G_v(b)u_{(h-1)w\Lambda})\otimes u_{w\Lambda}$ is one of the
weight components of $G_v(b)(u_{hw\Lambda})$. Thus
$$
(G_v(b)u_{(h-1)w\Lambda})\otimes u_{w\Lambda}\in
L_v(\Lambda)^{\otimes h}.
$$
As $(G_v(b)u_{(h-1)w\Lambda})\otimes u_{w\Lambda}
\not\in vL_v(\Lambda)^{\otimes h}$ because
$$
u_{w_1\Lambda}\otimes\cdots\otimes u_{w_{h-1}\Lambda}-
G_v(b)u_{(h-1)w\Lambda}\in vL_v(\Lambda)^{\otimes(h-1)},
$$
we may conclude that $G_v(b)(u_{hw\Lambda})\ne0$ and
$$
u_{w_1\Lambda}\otimes\cdots\otimes u_{w_h\Lambda}+
vL_v(\Lambda)^{\otimes h}=G_v(b)(u_{(h-1)w\Lambda}\otimes u_{w\Lambda})
+vL_v(\Lambda)^{\otimes h}.
$$
On the other hand, Lemma \ref{equality} implies
$$
G_v(b)(u_{(h-1)w\Lambda}\otimes u_{w\Lambda})
+vL_v(\Lambda)^{\otimes h}=
G_v(b)u_{hw\Lambda}+vL_v(\Lambda)^{\otimes h}\in B^w(h\Lambda).
$$
Thus we have proved
$$
u_{w_1\Lambda}\otimes\cdots\otimes u_{w_h\Lambda}+
vL_v(\Lambda)^{\otimes h}\in B^w(h\Lambda).
$$
\end{proof}

\begin{cor}
\label{cor to a LS condition}
Let $w\in W$. If there exists a sequence
$w_1\ge w_2\ge\cdots\ge w_h\ge w$ in $W$ such that
$\nu_i=w_i\emptyset_m$, for $1\le i\le h$, then
$$
\nu_1\otimes\cdots\otimes\nu_h\in B^w(h\Lambda_m)
\subset B(h\Lambda_m)\subset B(\Lambda_m)^{\otimes h}.
$$
\end{cor}

Define the $\mathbb Q(v)$-linear
anti-involution $*$ on $U_v^-(\mathfrak g)$
by $f_i^*=f_i$. It preserves the crystal lattice
of $U_v^-(\mathfrak g)$ \cite[Proposition 5.2.4]{K2}.
Then, as in \cite[Corollary 6.1.2]{K2},
$(G(b)^*,G(b)^*)\equiv (G(b),G(b))\equiv 1$ modulo
$v\mathbb Z[v]$ implies
$G(b)^*=\pm G(b^*)$, for some $b^*\in B(\infty)$.
Now, it is proved in \cite[Theorem 2.1.1]{K3}
that the minus sign does not occur. To summarize,
we have the following.

\begin{prop}
\label{star for negative part}
\item[(1)]
$B(\infty)^*=B(\infty)$.
\item[(2)]
$G_v(b^*)=G_v(b)^*$, for $b\in B(\infty)$.
\end{prop}

Next let $\tilde U_v(\mathfrak g)$ be
the modified quantized enveloping algebra. Namely,
$$
\tilde U_v(\mathfrak g)=\bigoplus_{\Lambda\in P}U_v(\mathfrak g)a_\Lambda
$$
such that $v^ha_\Lambda=a_\Lambda v^h=v^{\Lambda(h)}a_\Lambda$,
$a_\Lambda e_i=e_ia_{\Lambda-\alpha_i}$,
$a_\Lambda f_i=f_ia_{\Lambda+\alpha_i}$ and
$a_\Lambda a_{\Lambda'}=\delta_{\Lambda\Lambda'}a_\Lambda$.
Define the $\mathbb Q(v)$-linear anti-involution $*$ by
$$
(v^h)^*=v^{-h},\;\;e_i^*=e_i,\;\; f_i^*=f_i,\;\;a_\Lambda^*=a_{-\Lambda}.
$$
Lusztig constructed global bases for tensor products of integrable
highest weight and lowest weight $U_v(\mathfrak g)$-modules
\cite[24.3]{L}, and
showed that their inverse limits exist in $\tilde U_v(\mathfrak g)$.
Thus we have the crystal basis of $\tilde U_v(\mathfrak g)$
\cite[25.2]{L}. We denote the crystal by
$$
B(\tilde U_v(\mathfrak g))=\bigsqcup_{\Lambda\in P}
B(U_v(\mathfrak g)a_\Lambda).
$$
The global basis of $\tilde U_v(\mathfrak g)$ is also denoted by
$\{G_v(b) \mid b\in B(\tilde U_v(\mathfrak g))\}$.

\begin{thm}[{\cite[Theorem 3.1.1]{K4}}]
\label{crystal isom}
Let $\Lambda$ be an integral weight.
We choose dominant integral weights $\Lambda^+$ and $\Lambda^-$ such that
$\Lambda=\Lambda^+-\Lambda^-$. Then combining two embeddings
$$
B(\Lambda^+)\subset B(\infty)\otimes T_{\Lambda^+},\quad
B(-\Lambda^-)\subset T_{-\Lambda^-}\otimes B(-\infty)
$$
and $T_{\Lambda^+}\otimes T_{-\Lambda^-}=T_{\Lambda}$,
we have a strict embedding of crystals
$$
B(\Lambda^+)\otimes B(-\Lambda^-)\subset
B(\infty)\otimes T_{\Lambda} \otimes B(-\infty).
$$
By taking the direct limit, we have
$$
B(U_v(\mathfrak g)a_{\Lambda})\simeq
B(\infty)\otimes T_{\Lambda} \otimes B(-\infty).
$$
\end{thm}

In the remainder of this discussion, we identify $B(U_v(\mathfrak g)a_{\Lambda})$
with $B(\infty)\otimes T_{\Lambda} \otimes B(-\infty)$.
The following theorem generalizes Proposition \ref{star for negative part}.

\begin{thm}[{\cite[Theorem 4.3.2, Corollary 4.3.3]{K4}}]
\label{star}
\item[(1)]
$B(\tilde U_v(\mathfrak g))^*=B(\tilde U_v(\mathfrak g))$, and if
$b=b_1\otimes t_{\Lambda}\otimes b_2\in B(\tilde U_v(\mathfrak g))$
then
$$
b^*=b_1^*\otimes t_{-\Lambda-wt(b_1)-wt(b_2)}\otimes b_2^*.
$$
\item[(2)]
$G_v(b^*)=G_v(b)^*$, for $b\in B(\tilde U_v(\mathfrak g))$.
\end{thm}

Now we define, for $b\in B(\tilde U_v(\mathfrak g))$,
\begin{gather*}
\epsilon^*_i(b)=\epsilon_i(b^*),\;\;
\varphi^*_i(b)=\varphi_i(b^*),\;\;
wt^*(b)=wt(b^*),\\
\e_i^*b=(\e_ib^*)^*,\;\;
\f_i^*b=(\f_ib^*)^*.
\end{gather*}
Then this defines another crystal structure on $B(\tilde U_v(\mathfrak g))$,
which is called the {\bf star crystal structure}. The star crystal structure
is compatible with the original
crystal structure on $B(\tilde U_v(\mathfrak g))$
in the following sense.

\begin{thm}[{\cite[Theorem 5.1.1]{K4}}]
$\e_i^*$ and $\f_i^*$ are strict morphisms of crystals.
\end{thm}

Using the star crystal structure, we can define another Weyl group
action on $B(\tilde U_v(\mathfrak g))$. We denote the action by
$w^*b$, for $w\in W$ and $b\in B(\tilde U_v(\mathfrak g))$.

\begin{defn}
Let $B$ be a normal crystal.
An element $b\in B$ of weight $\Lambda$
is called {\bf extremal} if there exists
a subset $\{b_w\}_{w\in W}$ of $B$ such that
\begin{itemize}
\item[(i)]
$b_w=b$ if $w=1$.
\item[(ii)]
If $w\Lambda(h_i)\ge0$ then $\e_ib_w=0$ and
$\f_i^{max}b_w=b_{s_iw}$.
\item[(iii)]
If $w\Lambda(h_i)\le0$ then $\f_ib_w=0$ and
$\e_i^{max}b_w=b_{s_iw}$.
\end{itemize}
\end{defn}

\underline{}When $B=B(\Lambda)$ for a dominant integral weight $\Lambda$,
this is a natural crystal analogue of extremal weight
vectors in the highest weight module
$U_v(\mathfrak g)u_\Lambda$.

\begin{lemma}
\label{extremal elements}
Let $\Lambda$ be dominant integral. Define
$b_w=u_{w\Lambda}+vL_v(\Lambda)\in B(\Lambda)$, for $w\in W$.
\item[(1)]
The set of extremal elements
of $B(\Lambda)$ coincides with $\{b_w\}_{w\in W}$.
\item[(2)]
$\{u_{w\Lambda}\}_{w\in W}$ are extremal vectors. That is, we have the following.
\begin{itemize}
\item[(i)]
If $w\Lambda(h_i)\ge0$ then $e_iu_{w\Lambda}=0$ and
$f_i^{(w\Lambda(h_i))}u_{w\Lambda}=u_{s_iw\Lambda}$.
\item[(ii)]
If $w\Lambda(h_i)\le0$ then $f_iu_{w\Lambda}=0$ and
$e_i^{(-w\Lambda(h_i))}u_{w\Lambda}=u_{s_iw\Lambda}$.
\end{itemize}
\item[(3)]
If $s_iw<w$ and $b\in B(\infty)$ satisfies
$G_v(b)u_{w\Lambda}\ne0$ then $\epsilon_i^*(b)=0$.
\item[(4)]
Suppose that $s_iw<w$ and $b\in B(\infty)$ satisfies
$\epsilon_i^*(b)=0$. Then
$$
G_v(b)u_{w\Lambda}+vL_v(\Lambda)=
G_v(\tilde {f_i^*}^{-w\Lambda(h_i)}b)u_{s_iw\Lambda}+vL_v(\Lambda).
$$
\end{lemma}
\begin{proof}
(1) and (2) are well-known, and we only prove (3) and (4).
Note that $s_iw<w$ implies $w\Lambda(h_i)\le0$. Thus
$f_iu_{w\Lambda}=0$ by (2). If $\epsilon_i(b^*)>0$ then
$G_v(b^*)\in f_iU_v^-(\mathfrak g)$. Thus $G_v(b)\in U_v^-(\mathfrak g)f_i$
by Proposition \ref{star for negative part}. Then
$G_v(b)u_{w\Lambda}=0$, which contradicts our assumption.
We have proved $\epsilon_i^*(b)=0$.

To prove (4), note that $s_iw\Lambda(h_i)=-w\Lambda(h_i)\ge0$ and
$f_i^{(-w\Lambda(h_i))}u_{s_iw\Lambda}=u_{w\Lambda}$ by (2).
Now, $\epsilon_i(b^*)=0$ implies that
$$
\f_i^{-w\Lambda(h_i)}G_v(b^*)=f_i^{(-w\Lambda(h_i))}G_v(b^*).
$$
Hence, we have
$$
\bigl(\tilde {f_i^*}^{-w\Lambda(h_i)}G_v(b)\bigr)u_{s_iw\Lambda}=
G_v(b)f_i^{(-w\Lambda(h_i))}u_{s_iw\Lambda}=G_v(b)u_{w\Lambda}.
$$
Thus $G_v(\tilde {f_i^*}^{-w\Lambda(h_i)}b)u_{s_iw\Lambda}+vL_v(\Lambda)=
G_v(b)u_{w\Lambda}+vL_v(\Lambda)$ follows.
\end{proof}

\begin{defn}
Suppose that $\Lambda$ is dominant integral. For $w\in W$, we define
$$
B(w\Lambda)=\{b\in B(U_v(\mathfrak g)a_{w\Lambda}) \mid
\text{$b^*$ is extremal}\}.
$$
\end{defn}

We identify $B(w\Lambda)$ with a subcrystal of
$B(\infty)\otimes T_{w\Lambda} \otimes B(-\infty)$
through the crystal isomorphism given in Theorem
\ref{crystal isom}.
As the property that $b^*$ is extremal is stable under
$\e_i$ and $\f_i$, if we define
$I_{w\Lambda}$ to be the subspace of $U_v(\mathfrak g)a_{w\Lambda}$
spanned by $\{G_v(b) \mid b\not\in B(w\Lambda)\}$ then
it is a $U_v(\mathfrak g)$-submodule of $U_v(\mathfrak g)a_{w\Lambda}$.
The $U_v(\mathfrak g)$-module
$V_v(w\Lambda)=U_v(\mathfrak g)a_{w\Lambda}/I_{w\Lambda}$
is Kashiwara's extremal weight module.

\begin{thm}[{\cite[Proposition 8.2.2]{K4}}]
\label{isom for various w}
Suppose that $\Lambda$ is dominant integral.
\item[(1)]
$V_v(w\Lambda)$ is an integrable $U_v(\mathfrak g)$-module.
\item[(2)]
$B(w\Lambda)$ is the crystal graph of $V_v(w\Lambda)$.
\item[(3)]
The map $b\mapsto w^*b$, for $b\in B(\Lambda)$, defines an
isomorphism of crystals
$$
B(\Lambda)\simeq B(w\Lambda).
$$
\end{thm}

As $V_v(w\Lambda)$ is generated by the extremal vector of weight $w\Lambda$,
and integrable, $V_v(w\Lambda)$ with $w=1$ is the integrable
highest weight module
$U_v(\mathfrak g)u_\Lambda$. Hence
$B(w\Lambda)$ with $w=1$ is nothing but $B(\Lambda)$,
and there is no conflict in the notation.

Fix $\underline i$ and
let $\mathbb Z_{\underline i}$ be the polyhedral realization of
$B(\infty)$ as before. If $b\in B(\infty)$ corresponds to
$(\cdots,0,0,a_r,\cdots,a_2,a_1)\in \mathbb Z_{\underline i}$,
then the integers $a_k$ are determined by
$$
b^*=\f_{i_1}^{a_1}\f_{i_2}^{a_2}\cdots
\f_{i_r}^{a_r}u_\infty
$$
such that $\epsilon_{i_k}(\f_{i_{k+1}}^{a_{k+1}}
\f_{i_{k+2}}^{a_{k+2}}\cdots u_\infty)=0$, for all $k$.
See \cite[(2.35), (2.36)]{NZ}.

Define $S_h:B(\infty)\otimes T_{\Lambda}\otimes B(-\infty)
\rightarrow B(\infty)\otimes T_{h\Lambda}\otimes B(-\infty)$ by
$$
S_h(b_1\otimes t_\Lambda\otimes b_2)=
S_h(b_1)\otimes t_{h\Lambda}\otimes S_h(b_2).
$$
This is also a crystal morphism of amplitude $h$.

The next results are proved in \cite[Proposition 3.2, 3.5]{NS1} in
a slightly different manner.

\begin{lemma}
\label{conjugate linear}
\item[(1)]
Let $b\in B(\infty)$. Then $S_h(b)^*=S_h(b^*)$, for all $h$.
\item[(2)]
Let $b\in B(\tilde U_v(\mathfrak g))$. Then
$S_h(b)^*=S_h(b^*)$, for all $h$.
\end{lemma}
\begin{proof}
(1) We fix a polyhedral realization $\mathbb Z_{\underline i}$ of
$B(\infty)$ and denote by
$$
(\dots,0,0,a_r,\dots,a_2,a_1)
$$
the element which corresponds to $b$. Then $S_h(b)$ corresponds to
$$
(\dots,0,0,ha_r,\dots,ha_2,ha_1)
$$
by Proposition \ref{Sh for Binfty}. Thus, we have
\begin{multline*}
S_h(b)^*=\f_{i_1}^{ha_1}\f_{i_2}^{ha_2}\cdots
\f_{i_r}^{ha_r}u_\infty
=\f_{i_1}^{ha_1}\f_{i_2}^{ha_2}\cdots
\f_{i_r}^{ha_r}S_h(u_\infty)\\
=\f_{i_1}^{ha_1}\f_{i_2}^{ha_2}\cdots
\f_{i_{r-1}}^{ha_{r-1}}S_h(\f_{i_r}^{a_r}u_\infty)
=\cdots\cdots=S_h(\f_{i_1}^{a_1}\f_{i_2}^{a_2}\cdots
\f_{i_r}^{a_r}u_\infty).
\end{multline*}
Thus, $S_h(b)^*=S_h(b^*)$ as desired.

(2) Let $b=b'\otimes t_\Lambda\otimes b''$. Then $S_h(b)^*$
is equal to
$$
\bigl(S_h(b')\otimes t_{h\Lambda}\otimes S_h(b'')\bigr)^*
=S_h(b')^*\otimes t_{h(-\Lambda-wt(b')-wt(b''))}\otimes S_h(b'')^*.
$$
Since
$S_h(b^*)=S_h((b')^*)\otimes t_{h(-\Lambda-wt(b')-wt(b''))}
\otimes S_h((b'')^*)$,
$S_h(b)^*=S_h(b^*)$ follows by (1).
\end{proof}

\begin{lemma}
Let $b\in B(\tilde U_v(\mathfrak g))$. If $b^*$ is extremal,
so is $S_h(b)^*$.
\end{lemma}
\begin{proof}
By the definition of tensor product, we have
\begin{equation*}
\begin{split}
\epsilon_i(b_1\otimes t_\Lambda\otimes b_2)&=
\max(\epsilon_i(b_1),\epsilon_i(b_2)-(\Lambda+wt(b_1))(h_i)),\\
\varphi_i(b_1\otimes t_\Lambda\otimes b_2)&=
\max(\varphi_i(b_1)+(\Lambda+wt(b_2))(h_i),\varphi_i(b_2)).
\end{split}
\end{equation*}
Suppose that there exists
$\{b_w=b_w'\otimes t_{-w\Lambda}\otimes b_w''\}_{w\in W}$ such that
\begin{itemize}
\item[(i)]
$b_w^*=(b_w')^*\otimes t_{w\Lambda-wt(b_w')-wt(b_w'')}
\otimes (b_w'')^*=b^*$ if $w=1$.
\item[(ii)]
If $w\Lambda(h_i)\ge0$ then $\e_ib_w^*=0$ and
$\f_i^{max}b_w^*=b_{s_iw}^*$.
\item[(iii)]
If $w\Lambda(h_i)\le0$ then $\f_ib_w^*=0$ and
$\e_i^{max}b_w^*=b_{s_iw}^*$.
\end{itemize}
We want to show that $\{S_h(b_w)^*\}_{w\in W}$ satisfies
conditions (i) to (iii) above. As (i) is obvious,
we prove (ii) and (iii). Suppose that
$w\Lambda(h_i)\ge0$. Then $\tilde e_ib_w^*=0$ implies
$$
\epsilon_i(b_w^*)=
\max(\epsilon_i((b_w')^*),\epsilon_i((b_w'')^*)+(-w\Lambda+wt(b_w''))(h_i))
=0.
$$
By Lemma \ref{conjugate linear}, we have
\begin{equation*}
\begin{split}
\epsilon_i(S_h(b_w)^*)&=
\epsilon_i(S_h((b_w')^*)\otimes t_{h(w\Lambda-wt(b_w')-wt(b_w''))}
\otimes S_h((b_w'')^*))\\
&=\max(h\epsilon_i((b_w')^*),h\epsilon_i((b_w'')^*)+h(-w\Lambda+wt(b_w''))(h_i)).
\end{split}
\end{equation*}
Thus $\epsilon_i(S_h(b_w)^*)=h\epsilon_i(b_w^*)=0$ and
$\e_iS_h(b_w)^*=0$ follows. By a similar computation,
we have $\varphi_i(S_h(b_w)^*)=h\varphi_i(b_w^*)$, which implies
that
\begin{equation*}
\begin{split}
\f_i^{max}S_h(b_w)^*&=\f_i^{h\varphi_i(b_w^*)}S_h(b_w)^*
=\f_i^{h\varphi_i(b_w^*)}S_h(b_w^*)\\
&=S_h(\f_i^{\varphi_i(b_w^*)}b_w^*)
=S_h(b_{s_iw}^*)=S_h(b_{s_iw})^*.
\end{split}
\end{equation*}
Suppose that $w\Lambda(h_i)\le0$. Then, by similar arguments,
we have $\f_iS_h(b_w)^*=0$ and
$\e_i^{max}S_h(b_w)^*=S_h(b_{s_iw})^*$.
\end{proof}

Let $\Lambda$ be dominant integral, $w\in W$.
Since $B(w\Lambda)\simeq B(\Lambda)$ by
Theorem \ref{isom for various w}(3), we have
a unique crystal morphism
$B(w\Lambda)\rightarrow B(hw\Lambda)$ of amplitude $h$,
which we also denote by $S_h$. The following
corollary generalizes Proposition \ref{Sh for hw crystal}.

\begin{cor}
Let $\Lambda$ be dominant integral, $w\in W$. Then
we have the following commutative diagram.
\[
\begin{array}{ccc}
\label{commutativity of Sh}
B(w\Lambda) & \stackrel{S_h}{\longrightarrow} & B(hw\Lambda) \\
\cap & & \cap \\
B(\infty)\otimes T_{w\Lambda}\otimes B(-\infty)&
\stackrel{S_h}{\longrightarrow} &
B(\infty)\otimes T_{hw\Lambda}\otimes B(-\infty)
\end{array}
\]
\end{cor}

We need two formulas. In the lemma below,
(1) is taken from \cite[(3.1.1)]{K4} and
(2) is taken from \cite[Appendix]{K6}.

\begin{lemma}
\label{two formulas}
\item[(1)]
Let $b=b_1\otimes t_\Lambda\otimes b_2\in B(\tilde U_v(\mathfrak g))$. Then
$G_v(b)\in \tilde U_v(\mathfrak g)$ equals
$G_v(b_1)G_v(b_2)a_\Lambda$ plus the linear combination
$\sum X_iY_ia_\Lambda$, where $X_i\in U_v^-(\mathfrak g)_{-\alpha}$
and $Y_i\in U_v^+(\mathfrak g)_{\beta}$ such that
$ht(\alpha)<ht(wt(b_1))$ and $ht(\beta)<ht(wt(b_2))$ respectively.
In particular,
$$
G_v(b_1\otimes t_\Lambda\otimes u_{-\infty})=G_v(b_1)a_\Lambda.
$$
\item[(2)]
Let $b=b_1\otimes t_\Lambda\otimes u_{-\infty}$
and suppose that $b^*$ is extremal. Then
$$
s_i^*b=\begin{cases}
\tilde {f_i^*}^{-\Lambda(h_i)}b_1\otimes t_{s_i\Lambda}
\otimes u_{-\infty}\quad (\text{if $\epsilon_i^*(b)=0$.})\\
\tilde {e_i^*}^{max}b_1\otimes t_{s_i\Lambda}
\otimes \tilde {e_i^*}^{\Lambda(h_i)-\epsilon_i^*(b_1)}u_{-\infty}
\quad (\text{if $\varphi_i^*(b)=0$.})
\end{cases}
$$
\end{lemma}

\begin{prop}
\label{correspondence of Demazure crystal}
Suppose that $\Lambda$ is dominant integral.
\item[(1)]
If $b\in B^w(\Lambda)$ then
$w^*b\in B(\infty)\otimes t_{w\Lambda}\otimes u_{-\infty}$.
\item[(2)]
Under the isomorphism $B(\Lambda)\simeq B(w\Lambda)$ given by
$b\mapsto w^*b$,
$B^w(\Lambda)$ may be identified with
$$
\{b\in B(\infty)\otimes t_{w\Lambda}\otimes u_{-\infty} \mid
\text{$b^*$ is extremal}\}.
$$
\end{prop}
\begin{proof}
(1) We identify the extremal weight module $V_v(\Lambda)$ with
the highest weight module $U_v(\mathfrak g)u_\Lambda$ as before.
Write $G_v(b)=G_v(b')u_\Lambda$ in $U_v(\mathfrak g)u_\Lambda$. As
$$
G_v(b'\otimes t_\Lambda\otimes u_{-\infty})=G_v(b')a_\Lambda
$$
by Lemma \ref{two formulas}(1), we have
$b=b'\otimes t_\Lambda\otimes u_{-\infty}$
under the identification of the crystal of the highest weight module
$U_v(\mathfrak g)u_\Lambda$ with $B(\Lambda)$ which is defined by
the extremal weight module $V_v(\Lambda)$.

Suppose now that $b\in B^w(\Lambda)$. Then there exists
$b_1\in B(\infty)$ such that
$G_v(b)=G_v(b_1)u_{w\Lambda}$ by Lemma \ref{equality}.
Let $w=s_{i_1}\cdots s_{i_\ell}$ be a reduced expression.
Then Lemma \ref{extremal elements}(3),(4) imply that
$$
G_v(b)+vL_v(\Lambda)=
G_v(\tilde {f_{i_\ell}^*}^{a_\ell}\cdots\tilde {f_{i_1}^*}^{a_1}b_1)u_\Lambda
+vL_v(\Lambda),
$$
where $a_k=-s_{i_k}\cdots s_{i_\ell}\Lambda(h_{i_k})
=s_{i_{k+1}}\cdots s_{i_\ell}\Lambda(h_{i_k})$, such that
$$
\epsilon_{i_{k+1}}^*(\tilde {f_{i_k}^*}^{a_k}
\cdots\tilde {f_{i_1}^*}^{a_1}b_1)=
\epsilon_{i_{k+1}}(\f_{i_k}^{a_k}
\cdots\f_{i_1}^{a_1}b_1^*)=0,
$$
for $0\le k<\ell$. This implies
$G_v(b)=
G_v(\tilde {f_{i_\ell}^*}^{a_\ell}\cdots\tilde {f_{i_1}^*}^{a_1}b_1)u_\Lambda$.
Thus, by the first paragraph, we have
$$
b=\tilde {f_{i_\ell}^*}^{a_\ell}
\cdots\tilde {f_{i_1}^*}^{a_1}b_1\otimes t_\Lambda\otimes u_{-\infty}.
$$
We show by downward induction on $k$ that
$$
s_{i_{k+1}}^*\cdots s_{i_\ell}^*b=\tilde {f_{i_k}^*}^{a_k}
\cdots\tilde {f_{i_1}^*}^{a_1}b_1\otimes
t_{s_{i_{k+1}}\cdots s_{i_\ell}\Lambda}\otimes u_{-\infty}.
$$
If $k=\ell$ there is nothing to prove. Suppose that the equation holds
for $k$. As $s_{i_{k+1}}^*\cdots s_{i_\ell}^*b\in
B(s_{i_{k+1}}\cdots s_{i_\ell}\Lambda)$ by
Theorem \ref{isom for various w}(3),
$s_{i_{k+1}}\cdots s_{i_\ell}b^*$ is extremal. As
$$
wt(s_{i_{k+1}}\cdots s_{i_\ell}b^*)(h_{i_k})=
-s_{i_{k+1}}\cdots s_{i_\ell}\Lambda(h_{i_k})=-a_k\le0,
$$
we have $\varphi_{i_k}^*(s_{i_{k+1}}^*\cdots s_{i_\ell}^*b)=0$.
Thus Lemma \ref{two formulas}(2)
implies
$$
s_{i_k}^*\cdots s_{i_\ell}^*b
=\tilde {e_{i_k}^*}^{max}\tilde {f_{i_k}^*}^{a_k}
\cdots\tilde {f_{i_1}^*}^{a_1}b_1\otimes
t_{s_{i_k}\cdots s_{i_\ell}\Lambda}\otimes
\tilde {e_{i_k}^*}^{a_k-\epsilon_{i_k}
(\f_{i_k}^{a_k}\cdots\f_{i_1}^{a_1}b_1^*)}
u_{-\infty}.
$$
Since the formula $\epsilon_{i_{k+1}}^*(\tilde {f_{i_k}^*}^{a_k}
\cdots\tilde {f_{i_1}^*}^{a_1}b_1)=0$ implies
$\epsilon_{i_k}(\f_{i_k}^{a_k}\cdots\f_{i_1}^{a_1}b_1^*)=a_k$ if we
replace $k$ with $k-1$ in the formula, we have the equation for $k-1$. As a result, we have
$w^*b=b_1\otimes t_{w\Lambda}\otimes u_{-\infty}\in
B(\infty)\otimes t_{w\Lambda}\otimes u_{-\infty}$.

(2) We only have to show that if
$b=b_1\otimes t_{w\Lambda}\otimes u_{-\infty}\in B(w\Lambda)$ then
we have $(w^{-1})^*b\in B^w(\Lambda)$. Define
$a_k=s_{i_{k+1}}\cdots s_{i_\ell}\Lambda(h_{i_k})$. We show by
induction on $k$ that
$$
s_{i_k}^*\cdots s_{i_1}^*b=\tilde {f_{i_k}^*}^{a_k}
\cdots\tilde {f_{i_1}^*}^{a_1}b_1\otimes
t_{s_{i_{k+1}}\cdots s_{i_\ell}\Lambda}\otimes u_{-\infty}.
$$
If $k=0$ there is nothing to prove. Suppose that the equation holds for
$k$. As $s_{i_k}\cdots s_{i_1}b^*$ is extremal and
$$
wt(s_{i_k}\cdots s_{i_1}b^*)(h_{i_{k+1}})=
-s_{i_{k+1}}\cdots s_{i_\ell}\Lambda(h_{i_{k+1}})
=a_{k+1}\ge0,
$$
we have $\epsilon_{i_{k+1}}(s_{i_k}\cdots s_{i_1}b^*)=0$.
Thus Lemma \ref{two formulas}(2) implies the equation
for $k+1$. As a result, we have
$$
(w^{-1})^*b=\tilde {f_{i_\ell}^*}^{a_\ell}
\cdots\tilde {f_{i_1}^*}^{a_1}b_1\otimes t_\Lambda\otimes u_{-\infty}.
$$
Now, $\epsilon^*_{i_{k+1}}(\tilde {f_{i_k}^*}^{a_k}
\cdots\tilde {f_{i_1}^*}^{a_1}b_1)=0$, for $0\le k<\ell$, because
\begin{equation*}
\begin{split}
0=\epsilon^*_{i_{k+1}}(s^*_{i_k}\cdots s^*_{i_1}b)
&=\epsilon^*_{i_{k+1}}(\tilde {f_{i_k}^*}^{a_k}
\cdots\tilde {f_{i_1}^*}^{a_1}b_1\otimes
t_{s_{i_{k+1}}\cdots s_{i_\ell}\Lambda}\otimes u_{-\infty})\\
&\ge \epsilon^*_{i_{k+1}}(\tilde {f_{i_k}^*}^{a_k}
\cdots\tilde {f_{i_1}^*}^{a_1}b_1)\ge 0.
\end{split}
\end{equation*}
Thus Lemma \ref{extremal elements}(4) shows that
$$
G_v(b_1)u_{w\Lambda}=G_v(\tilde {f_{i_\ell}^*}^{a_\ell}
\cdots\tilde {f_{i_1}^*}^{a_1}b_1)u_\Lambda=G_v((w^{-1})^*b).
$$
Therefore, we
have $(w^{-1})^*b\in B^w(\Lambda)$ by Lemma \ref{equality}.
\end{proof}

\texttt{}The following is a theorem proved by Kashiwara and Sagaki independently.
The proof for the first equality works for general dominant integral weights.

\begin{thm}
\label{floor theorem}
Suppose $w\in W/W_m$. Then
$$
B^w(\Lambda_m)=\{\lambda\in B(\Lambda_m) \mid f(\lambda)\ge w\Lambda_m\}
=\{\lambda\in B(\Lambda_m) \mid \floor(\lambda)\supset w\emptyset_m\}.
$$
\end{thm}
\begin{proof}
If we write $\floor(\lambda)=u\emptyset_m$, for $u\in W/W_m$,
then $f(\lambda)=wt(\floor(\lambda))=u\Lambda_m$, and
$f(\lambda)\ge w\Lambda_m$ if and only if $u\ge w$.
Thus, the second equality follows from
Proposition \ref{two orders}. We prove the
first equality.

Suppose that $h$ is sufficiently divisible and write
$S_h(\lambda)=\nu_1\otimes\cdots\otimes\nu_h$, for $\lambda$ with
$f(\lambda)\ge w\Lambda_m$. Then
there exists a sequence $w_1\ge\cdots\ge w_h\ge w$ in $W$ such that
$\nu_i=w_i\emptyset_m$, for $1\le i\le h$.
By Corollay \ref{cor to a LS condition},
we have $S_h(\lambda)\in B^w(h\Lambda_m)$. We want to show
$\lambda\in B^w(\Lambda_m)$. Let us consider the
crystal morphism of amplitude $h$:
$$
B(\infty)\otimes T_{w\Lambda_m}\otimes B(-\infty)\longrightarrow
B(\infty)\otimes T_{hw\Lambda_m}\otimes B(-\infty).
$$
Then it induces $S_h: B(w\Lambda_m) \rightarrow B(hw\Lambda_m)$ by
Corollary \ref{commutativity of Sh}.

Write $w^*\lambda=b_1\otimes t_{w\Lambda_m}\otimes b_2\in B(w\Lambda_m)$.
Note that we have
$S_h(w^*\lambda)=w^*S_h(\lambda)$
by the uniqueness of the crystal morphism of amplitude
$h$ given in Proposition \ref{Sh for hw crystal}.
Since $S_h(\lambda)\in B^w(h\Lambda_m)$, we have
$$
S_h(b_1)\otimes t_{hw\Lambda}\otimes S_h(b_2)
=S_h(w^*\lambda)=w^*S_h(\lambda)
\in B(\infty)\otimes t_{hw\Lambda_m}\otimes u_{-\infty}
$$
by Proposition \ref{correspondence of Demazure crystal}(2),
which implies $S_h(b_2)=u_{-\infty}$.
Since $S_h:B(\infty)\rightarrow B(\infty)$ is injective
by Proposition \ref{Sh for Binfty},
we have $w^*\lambda=b_1\otimes t_{w\Lambda_m}\otimes u_{-\infty}$.
Therefore, Proposition \ref{correspondence of Demazure crystal}(2)
implies that $\lambda\in B^w(\Lambda_m)$.

Next suppose that $\lambda\in B^w(\Lambda_m)$. Then,
we have $S_h(\lambda)\in B^w(h\Lambda_m)$ by the similar argument. Take
sufficiently divisible $h$ and write
$S_h(\lambda)=\mu_1\otimes\cdots\otimes\mu_h$, for $\mu_1\geq\cdots\geq\mu_h$. Then
$S_h(\lambda)\in B^w(h\Lambda_m)$ implies that
$$
G_v(\mu_1\otimes\cdots\otimes\mu_h)\in
U_v^-(\mathfrak g)(u_{w\Lambda_m}\otimes\cdots\otimes u_{w\Lambda_m})
\subset V_v(\Lambda_m)\otimes\cdots\otimes V_v(\Lambda_m).
$$
Expand $G_v(\mu_1\otimes\cdots\otimes\mu_h)$ in the basis
$\{G_v(\nu_1)\otimes\cdots\otimes G_v(\nu_h) \mid \nu_1,\dots,\nu_h\in B(\Lambda_m)\}$.
If $G_v(\nu_1)\otimes\cdots\otimes G_v(\nu_h)$ appears in the expansion then
$\nu_1,\dots,\nu_h\in B^w(\Lambda_m)$, since
$$
U_v^-(\mathfrak g)(u_{w\Lambda_m}\otimes\cdots\otimes u_{w\Lambda_m})\subset
U_v^-(\mathfrak g)u_{w\Lambda_m}\otimes\cdots\otimes U_v^-(\mathfrak g)u_{w\Lambda_m}.
$$
In particular, we have $\mu_1,\dots,\mu_h\in B^w(\Lambda_m)$.
Write $\mu_h=y\emptyset_m$, for $y\in W/W_m$, and apply
Proposition \ref{basic properties}(4). Then $y\geq w$ and
$f(\lambda)=wt(\mu_h)\geq w\Lambda_m$ follows.
\end{proof}

\section{A property of Base}

We write $\lambda\le\mu$ for $\lambda\subset\mu$ in this and the next
sections.

Let $\lambda\in B(\Lambda_m)$ be $\lambda=(\lambda_0,\lambda_1,\dots)$.
We denote $\mu=(\lambda_1,\lambda_2,\dots)$ and write
$\lambda=\{\lambda_0\}\cup\mu$. In this section we shall show
$\base(\lambda)=\base(\{\lambda_0\}\cup\base(\mu))$.

\begin{defn}
Let $J\subset\mathbb Z$ and $x\in\mathbb Z$. Then we denote
$J\cap\mathbb Z_{\le x}$ by $J_{\le x}$.
\end{defn}

\begin{lemma}
\label{easy case}
Let $\lambda\in B(\Lambda_m)$, $J$ the corresponding set of beta
numbers of charge $m$, $j_0=\max J$. Write $K=J_{\le j_0-1}$.
Define $t=\min\{i\ge0 \mid \down^i(K)=\base(K)\}$.
\item[(1)]
Suppose that $j_0-e\not\in J$. Then the partition associated with
$\base(K)\cup\{j_0\}$ is $e$-restricted and
$\base(J)=\base(\base(K)\cup\{j_0\})$.
\item[(2)]
Suppose that $j_0-e\in J$ and fix $0\le s\le t$.
If there exists no $0\le i<s$ such that
$$
j_0-e=\min W(\down^i(K))<\min U(\down^i(K))\le j_0-1,
$$
then $\down^s(J)=\down^s(K)\cup\{j_0\}$. Furthermore,
\begin{itemize}
\item[(i)]
if $s<t$ then $U(\down^s(J))\neq\emptyset$ and
$\min U(\down^s(J))=\min U(\down^s(K))$,
\item[(ii)]
if $s=t$ then the partition associated with
$\base(K)\cup\{j_0\}$ is $e$-restricted and
$\base(J)=\base(\base(K)\cup\{j_0\})$.
\end{itemize}
\end{lemma}
\begin{proof}
Define $J_i=\down^i(K\cup\{j_0\})$ and
$K_i=\down^i(K)$, for $0\le i\le t$.

(1) We prove $J_i=K_i\cup\{j_0\}$,
$j_0-e\not\in J_i$ and $\max K_i\le j_0-1$
by induction on $i$.

When $i=0$, there is nothing to prove. Suppose that
$0\le i<t$ and that the claim holds for $i$. We want to show that
$J_{i+1}=K_{i+1}\cup\{j_0\}$,
$j_0-e\not\in J_{i+1}$ and $\max K_{i+1}\le j_0-1$.
As $i<t$, we have $U(K_i)\ne\emptyset$ and
$$
U(K_i)\subset U(J_i)=U(K_i\cup\{j_0\})\subset U(K_i)\cup\{j_0\}.
$$
As $\min U(K_i)\le\max K_i\le j_0-1$ we have
$\min U(J_i)=\min U(K_i)$, which we denote by $p'$.
Hence $p'\le j_0-1$ and $p'-e\ne j_0-e$, which implies
$j_0-e\not\in J_{i+1}$. We show that
$\min W(J_i)=\min W(K_i)$. Let $q'=\min W(J_i)$.
As $q'\le p'\le j_0-1$ and $q'\in J_i=K_i\cup\{j_0\}$,
we have $q'\in K_i$. If $q'=p'$ then $q'\in W(K_i)$.
If $q'<p'$ then $q'+e\not\in K_i$ because
$q'+e\not\in J_i$. Thus we also have $q'\in W(K_i)$.
Suppose that
there exists $p'-e<x<q'$ such that $x\in K_i$ and
$x+e\not\in K_i$. If $x+e\not\in J_i$ then
the minimality of $q'$ is contradicted. If
$x+e\in J_i$ then $x\in J_i$ and $x+e=j_0$, which contradicts
the induction hypothesis $j_0-e\not\in J_i$.
We have proved $\min W(K_i)=\min W(J_i)$.
Therefore, $\max K_{i+1}\le\max K_i\le j_0-1$ and
$$
J_{i+1}=\down(J_i)=\down(K_i)\cup\{j_0\}=K_{i+1}\cup\{j_0\}.
$$
Now, $J_t=\base(K)\cup\{j_0\}$ is associated with an $e$-restricted
partition by Lemma \ref{down operation}(2), and
$\base(J)=\base(\base(K)\cup\{j_0\})$ follows.

(2) We prove that $J_i=K_i\cup\{j_0\}$ and
$\max K_i\le j_0-1$, for $0\le i\le s$. Suppose that
$0\le i<s$ and that the claim holds for $i$. As $U(K_i)\ne\emptyset$,
we have $U(J_i)\ne\emptyset$ and
$$
p'=\min U(J_i)=\min U(K_i)\le j_0-1
$$
as before. Let $q'=\min W(J_i)$.
By the same argument as in (1),
we also have $q'\in W(K_i)$. Suppose that
there is $p'-e<x<q'$ such that $x\in K_i$ and
$x+e\not\in K_i$. If $x+e\not\in J_i$ then the minimality of
$q'$ is contradicted. If
$x+e\in J_i$ then $x+e=j_0$.
Thus
$$
j_0-e=\min W(K_i)<q'\le p'=\min U(K_i)\le j_0-1,
$$
which contradicts our assumption. Hence we have
$\min W(K_i)=\min W(J_i)$ and
$J_{i+1}=K_{i+1}\cup\{j_0\}$ follows.
We also have
$\max K_{i+1}\le\max K_i\le j_0-1$. By setting $i=s$,
we obtain $\down^s(J)=\down^s(K)\cup\{j_0\}$.

If $s<t$ then $U(K_s)\neq\emptyset$ and we have $U(J_s)\ne\emptyset$ and
$\min U(J_s)=\min U(K_s)$ by the same argument as above.
If $s=t$ then $J_t=\base(K)\cup\{j_0\}$
is associated with an $e$-restricted partition and we have
$\base(J)=\base(\base(K)\cup\{j_0\})$.
\end{proof}

Let $\lambda\in B(\Lambda_m)$, $J$, $j_0$,
$K$ and $t$ as above.

In the rest of this section we assume that $j_0-e\in J$ and
that there exists $0\le a<t$ such that $U(\down^a(J))\ne\emptyset$ and
\begin{itemize}
\item[(i)]
$\down^i(J)=\down^i(K)\cup\{j_0\}$,
for $0\le i\le a$.
\item[(ii)]
$p''=\min U(\down^a(K))$ and $q''=\min W(\down^a(K))$ satisfy
$$
p''=\min U(\down^a(J)),\quad q''=j_0-e<p''\le j_0-1.
$$
\end{itemize}
We also define $p'=\min U(\down^a(J))$ and $q'=\min W(\down^a(J))$.
Note that $q''=j_0-e\not\in W(\down^a(J))$ by
$p'=p''\ne q''$ and $\down^a(J)=\down^a(K)\cup\{j_0\}$.
Hence, $q'\neq q''$ and $\down^{a+1}(J)\ne\down^{a+1}(K)\cup\{j_0\}$.
More precisely, we have
$$
\down^{a+1}(K)=(\down^{a+1}(J)\setminus\{j_0,j_0-e\})\cup\{q'\}.
$$
Further, $q'>q''$ since $q'\le q''$ would imply $q'<p'$ and
$q'\in W(\down^a(K))$, which contradicts $q''=\min W(\down^a(K))$.
Thus we must have
$$
j_0-e<q'\le p'\le j_0-1.
$$

We also have $j_0-e\mathbb Z_{\ge0}\subset\down^a(J)$ and
$U(\down^a(J))=U(\down^a(K))$.
In fact, by $j_0-e\in\down^a(K)\subset\down^a(J)$ and $j_0-e<p'$,
$j_0-ke\in\down^a(J)$, for $k\ge1$. As $j_0\in\down^a(J)$, we
conclude that $j_0-e\mathbb Z_{\ge0}\subset\down^a(J)$.
Then
$$
U(\down^a(K))\subset U(\down^a(J))\subset U(\down^a(K))\cup\{j_0\}
$$
implies $U(\down^a(J))=U(\down^a(K))$.

\begin{defn}
Let $x\in J$.
\item[(1)]
We define the {\bf runner index of $x$}, which we denote by
$r(x)$, by
$$
1\le r(x)\le e\;\;\text{and}\;\;
x+e\mathbb Z=j_0+r(x)+e\mathbb Z.
$$
\item[(2)]
The {\bf layer level of $x$}, which we denote by $\ell(x)$, is defined by
$$
\ell(x)=-\frac{\min\{z\in j_0+e\mathbb Z | z\ge x\}-j_0}{e}.
$$
\end{defn}

The definitions are naturally understood on the abacus display
which is adjusted by $j_0$. Namely, we display $J$ on the abacus
in such a way that $j_0$ is on the rightmost runner. Then the
runner index is $1$ to $e$ from left to right, and $x$ is
$\ell(x)$ rows higher than $j_0$ in this $j_0$-adjusted abacus
display.

Define $b\ge1$ by $b=\min\{i\ge0\mid \base(J)=\base(\down^a(J))=\down^{a+i}(J)\}$,
and, this time, we define
$$
J_i=\down^i(\down^a(J))\;\;\text{and}\;\;K_i=\down^i(\down^a(K)),
$$
for $0\le i\le b$.
We set
$p_i'=\min U(J_i)$, $q_i'=\min W(J_i)$, for $0\leq i<b$.
Note that we have either $\ell(p'_i)=\ell(q'_i)$ or
$\ell(p'_i)=\ell(q'_i)-1$.
We also define
$p_i''=\min U(K_i)$ and $q_i''=\min W(K_i)$ if $U(K_i)\ne\emptyset$.

\begin{defn}
We say that $0\le i<b$ is a {\bf reset point} if
$\ell(p'_i)=\ell(q'_i)=0$.
\end{defn}

As $j_0-e<q'_0\leq p'_0\leq j_0-1$, $i=0$ is a reset point.

\begin{defn}
$U$ is the set of indices $0\le i<b$ such that $\ell(q'_i)=\ell(p'_i)$.
\end{defn}

$U$ is also the set of indices $0\le i<b$ such that $r(q'_i)\le
r(p'_i)$. Now, we analyze the relationship between $J_i$ and $K_i$
in detail. We start with an example.

\begin{Ex}
If $q'=p'$ and
\begin{equation*}
J_0:\;\;
\begin{aligned}
\times&\times&\times&\times&\times&\times\\
\times&      &      &\times&      &\times\\
\times&      &      &      &      &\times\\
\times&      &      &      &      &\times\\
      &      &      &      &      &\times\\
      &\times&      &      &\times&\times\\
\end{aligned}\qquad\qquad
K_0:\;\;
\begin{aligned}
\times&\times&\times&\times&\times&\times\\
\times&      &      &\times&      &\times\\
\times&      &      &      &      &\times\\
\times&      &      &      &      &\times\\
      &      &      &      &      &\times\\
      &\times&      &      &\times&      \\
\end{aligned}
\end{equation*}
then $K_0=J_0\setminus\{j_0\}$ and $0\in U$. We compute
$J_i$ and $K_i$, for $i>0$.

\begin{equation*}
J_1:\;\;
\begin{aligned}
\times&\times&\times&\times&\times&\times\\
\times&      &      &\times&      &\times\\
\times&      &      &      &      &\times\\
\times&      &      &      &      &\times\\
      &\times&      &      &      &\times\\
      &      &      &      &\times&\times\\
\end{aligned}\qquad\qquad
K_1:\;\;
\begin{aligned}
\times&\times&\times&\times&\times&\times\\
\times&      &      &\times&      &\times\\
\times&      &      &      &      &\times\\
\times&      &      &      &      &\times\\
      &\times&      &      &      &      \\
      &\times&      &      &\times&      \\
\end{aligned}
\end{equation*}
Thus, $K_1=(J_1\setminus\{j_0,j_0-e\})\sqcup\{q'_0\}$ and $1\in U$.

\begin{equation*}
J_2:\;\;
\begin{aligned}
\times&\times&\times&\times&\times&\times\\
\times&      &      &\times&      &\times\\
\times&      &      &      &      &\times\\
\times&\times&      &      &      &\times\\
      &      &      &      &      &\times\\
      &      &      &      &\times&\times\\
\end{aligned}\qquad\qquad
K_2:\;\;
\begin{aligned}
\times&\times&\times&\times&\times&\times\\
\times&      &      &\times&      &\times\\
\times&      &      &      &      &\times\\
\times&\times&      &      &      &      \\
      &\times&      &      &      &      \\
      &\times&      &      &\times&      \\
\end{aligned}
\end{equation*}
Thus, $K_2=(J_2\setminus\{j_0,j_0-e,j_0-2e\})\sqcup\{q'_0,q'_1\}$
and $2\in U$.

\begin{equation*}
J_3:\;\;
\begin{aligned}
\times&\times&\times&\times&\times&\times\\
\times&      &      &\times&      &\times\\
\times&\times&      &      &      &\times\\
      &\times&      &      &      &\times\\
      &      &      &      &      &\times\\
      &      &      &      &\times&\times\\
\end{aligned}\qquad\qquad
K_3:\;\;
\begin{aligned}
\times&\times&\times&\times&\times&\times\\
\times&      &      &\times&      &\times\\
\times&\times&      &      &      &      \\
\times&\times&      &      &      &      \\
      &\times&      &      &      &      \\
      &\times&      &      &\times&      \\
\end{aligned}
\end{equation*}
Thus, $K_3=(J_3\setminus\{j_0,j_0-e,j_0-2e,j_0-3e\})\sqcup\{q'_0,q'_1,q'_2\}$
and $3\not\in U$.

\begin{equation*}
J_4:\;\;
\begin{aligned}
\times&\times&\times&\times&\times&\times\\
\times&\times&      &      &      &\times\\
\times&\times&      &      &      &\times\\
      &\times&      &      &      &\times\\
      &      &      &      &      &\times\\
      &      &      &      &\times&\times\\
\end{aligned}\qquad\qquad
K_4:\;\;
\begin{aligned}
\times&\times&\times&\times&\times&\times\\
\times&\times&      &      &      &\times\\
\times&\times&      &      &      &      \\
\times&\times&      &      &      &      \\
      &\times&      &      &      &      \\
      &\times&      &      &\times&      \\
\end{aligned}
\end{equation*}
Thus, $K_4=(J_4\setminus\{j_0,j_0-e,j_0-2e,j_0-3e\})\sqcup\{q'_0,q'_1,q'_2\}$.
Note that $i=4$ is a reset point. We also have $4\in U$.

\begin{equation*}
J_5:\;\;
\begin{aligned}
\times&\times&\times&\times&\times&\times\\
\times&\times&      &      &      &\times\\
\times&\times&      &      &      &\times\\
      &\times&      &      &      &\times\\
      &      &      &      &\times&\times\\
      &      &      &      &      &\times\\
\end{aligned}\qquad\qquad
K_5:\;\;
\begin{aligned}
\times&\times&\times&\times&\times&\times\\
\times&\times&      &      &      &\times\\
\times&\times&      &      &      &      \\
\times&\times&      &      &      &      \\
      &\times&      &      &\times&      \\
      &      &      &      &\times&      \\
\end{aligned}
\end{equation*}
Thus, $K_5=(J_5\setminus\{j_0,j_0-e,j_0-2e,j_0-3e\})\sqcup\{q'_4,q'_1,q'_2\}$
and $5\in U$.

\begin{equation*}
J_6:\;\;
\begin{aligned}
\times&\times&\times&\times&\times&\times\\
\times&\times&      &      &      &\times\\
\times&\times&      &      &      &\times\\
      &\times&      &      &\times&\times\\
      &      &      &      &      &\times\\
      &      &      &      &      &\times\\
\end{aligned}\qquad\qquad
K_6:\;\;
\begin{aligned}
\times&\times&\times&\times&\times&\times\\
\times&\times&      &      &      &\times\\
\times&\times&      &      &      &      \\
\times&\times&      &      &\times&      \\
      &      &      &      &\times&      \\
      &      &      &      &\times&      \\
\end{aligned}
\end{equation*}
Thus, $K_6=(J_6\setminus\{j_0,j_0-e,j_0-2e,j_0-3e\})\sqcup\{q'_4,q'_5,q'_2\}$
and $6\in U$.

\begin{equation*}
J_7:\;\;
\begin{aligned}
\times&\times&\times&\times&\times&\times\\
\times&\times&      &      &      &\times\\
\times&\times&      &      &\times&\times\\
      &      &      &      &\times&\times\\
      &      &      &      &      &\times\\
      &      &      &      &      &\times\\
\end{aligned}\qquad\qquad
K_7:\;\;
\begin{aligned}
\times&\times&\times&\times&\times&\times\\
\times&\times&      &      &      &\times\\
\times&\times&      &      &\times&      \\
      &\times&      &      &\times&      \\
      &      &      &      &\times&      \\
      &      &      &      &\times&      \\
\end{aligned}
\end{equation*}
Thus, $K_7=(J_7\setminus\{j_0,j_0-e,j_0-2e,j_0-3e\})\sqcup\{q'_4,q'_5,q'_6\}$
and $7\in U$.

\begin{equation*}
J_8:\;\;
\begin{aligned}
\times&\times&\times&\times&\times&\times\\
\times&\times&      &      &\times&\times\\
      &\times&      &      &\times&\times\\
      &      &      &      &\times&\times\\
      &      &      &      &      &\times\\
      &      &      &      &      &\times\\
\end{aligned}\qquad\qquad
K_8:\;\;
\begin{aligned}
\times&\times&\times&\times&\times&\times\\
\times&\times&      &      &\times&      \\
\times&\times&      &      &\times&      \\
      &\times&      &      &\times&      \\
      &      &      &      &\times&      \\
      &      &      &      &\times&      \\
\end{aligned}
\end{equation*}
We finish with
$K_8=(J_8\setminus\{j_0,j_0-e,j_0-2e,j_0-3e,j_0-4e\})
\sqcup\{q'_4,q'_5,q'_6,q'_7\}$.
\end{Ex}

\begin{lemma}
\label{comparison}
Define $p'_b=p'_{b-1}-e$. Then, for each $0\le i\le b$, there exist $m_i\ge0$ and
$x_0,\dots,x_{m_i-1}\in K_i\setminus\{p'_i\}$ such that $U(J_i)=U(K_i)$ and
\begin{itemize}
\item[(a)]
$j_0-e\mathbb Z_{\ge0}\subset J_i$ and $\max J_i=j_0$.
\item[(b)]
$K_i=\left(J_i\setminus\{j_0,j_0-e,\dots,j_0-m_ie\}\right)
\sqcup\{x_0,\dots,x_{m_i-1}\}$.
\item[(c)]
If $x\in J_i$ is such that $r(x)\le r(p'_i)$ then $x\not\in U(J_i)$
unless $x=p'_i$.
\item[(d)]
If $x\in K_i$ is such that $r(x)\le r(p'_i)$ then $x\not\in U(K_i)$
unless $x=p'_i$.
\item[(e)]
$\ell(x_k)=k$, for $0\le k\le m_i-1$.
\item[(f)]
$r(p'_i)\ge r(x_0)\ge\cdots\ge r(x_{m_i-1})$.
\item[(g)]
If there exists $x\in J_i$ such that $1\le\ell(x)\le m_i$ and $r(p'_i)<r(x)<e$
then $(x+e\Z)\cap \Z_{\le j_0}\subset J_i$.
\item[(h)]
If $j_0-(k+1)e<x<x_k$, for some $0\le k\le m_i-1$, then $x\not\in J_i$
and $x\not\in K_i$.
\end{itemize}
Further, $1\le m_1\le\cdots\le m_b$.
\end{lemma}
\begin{proof}
$m_1\le\cdots\le m_b$ follows from (a), (b) and (f) because
$p''_i=p'_i\not\in j_0+e\Z$ implies that elements cannot
be added to $K_i\cap(j_0+e\Z)$, only removed.

$i=0$ is a reset point and we already know that the claims hold
when $i=0$; $m_0=0$ and (e), (f), (g) and (h) are vacant conditions.
Let $i_1$ be a reset point and
assume that the claims hold when $i\le i_1$.
Let $i_2\leq b-1$ be maximal
such that $p_i'$ decreases in the interval $i_1\leq i\leq i_2$.
We showed in section 2 that $p'_{i+1}=p'_i-e$ for $i_1\le i<i_2$ and that
$p'_{i_2+1}>p'_{i_1}$ if $i_2+1<b$.
We show that the claims hold for $i_1\le i\le i_2+1$ and
$m_1\le\cdots\le m_{i_2+1}$.
If $i_2+1<b$ then $i_2+1$ is a reset point because
$\ell(p'_{i_2+1})=0$ by $p'_{i_2+1}>p'_{i_1}$ and
$\ell(q'_{i_2+1})=0$ by (a) and (g) for $i=i_2+1$.
The condition (g) for $i=i_2+1$ is not vacant since
we already know $m_1=1$.
We repeat this process until $b$ is reached.

As we will see in the proof below, three patterns appear in the interval
$i_1\le i\le i_2+1$.
The first pattern occurs in the interval $i_1\le i<i_1+m_{i_1}$,
thus it does not occur when $i_1=0$,
and we reach $i=i_2+1$ when we are performing the second or the third pattern.
We will show that $p'_{i_2+1}-ke\not\in J_{i_2+1}$, for
$1\le k\le m_{i_2+1}$, when $i_2+1$ is a reset point. Hence,
we may assume that
$p'_{i_1}-ke\not\in J_{i_1}$, for $1\le k\le m_{i_1}$, when the first
pattern occurs at $i=i_1$.

Let $i=i_1+k$. When $k=0$, $U(K_{i_1})=U(J_{i_1})\ne\emptyset$ and
we have $x_0,\dots, x_{m-1}\in K_{i_1}\setminus\{p'_{i_1}\}$ which
satisfy (a) to (h), for $m=m_{i_1}$.
We want to show that $i_1+m\le b$ and that
the claims hold for $i_1\le i\le i_1+m$. If $m=0$ then
there is nothing to prove. Suppose that $m>0$ and
$p'_{i_1}-je\not\in J_{i_1}$, for $1\le j\le m$.
Then $x_0\ne p'_{i_1}$ and (f) for $i=i_1$ imply that
$r(x_k)<r(p'_{i_1})$, for $0\le k\le m-1$.
We shall show the following ($\dot a$) to ($\dot h$), for $0\le k\le m$,
by induction on $k$.
\begin{itemize}
\item[($\dot a$)]
$j_0-e\mathbb Z_{\ge0}\subset J_{i_1+k}$ and $\max J_{i_1+k}=j_0$.
\item[($\dot b$)]
$K_{i_1+k}=\left(J_{i_1+k}\setminus\{j_0,j_0-e,\dots,j_0-me\}\right)
\sqcup\{q'_{i_1},\dots,q'_{i_1+k-1},x_k,\dots, x_{m-1}\}$.
\item[($\dot c$)]
If $x\in J_{i_1+k}$ is such that $r(x)\le r(p'_{i_1})$ then
$x\not\in U(J_{i_1+k})$ unless $x=p'_{i_1+k}$.
\item[($\dot d$)]
If $x\in K_{i_1+k}$ is such that $r(x)\le r(p'_{i_1})$ then
$x\not\in U(K_{i_1+k})$ unless $x=p'_{i_1+k}$.
\item[($\dot e$)]
$\ell(q'_{i_1+j})=j$, for $0\le j\le k-1$.
\item[($\dot f$)]
$r(p'_{i_1})\ge r(q'_{i_1})\ge\cdots\ge
r(q'_{i_1+k-1})\ge r(x_k)\ge\cdots\ge r(x_{m-1})$.
\item[($\dot g$)]
If there exists $x\in J_{i_1+k}$ such that
$1\le\ell(x)\le m$ and $r(p'_{i_1})<r(x)<e$ then
$(x+e\mathbb Z)\cap J_{\le j_0}\subset J_{i_1+k}$.
\item[($\dot h$)]
If $\ell(x)=j$ and $r(x)<r(q'_{i_1+j})$, for some $0\le j\le k-1$,
or $\ell(x)=j$ and $r(x)<r(x_j)$, for some $k\le j\le m-1$,
then $x\not\in J_{i_1+k}$ and $x\not\in K_{i_1+k}$.
\end{itemize}
If $k\le m-1$ we also show $p'_{i_1}-je\not\in J_{i_1+k}$, for
$k+1\le j\le m$, $i_1+k<b$ and $p'_{i_1+k}=p'_{i_1}-ke$.

Before proving these claims, we explain that these imply the desired claims for
$i_1\le i\le i_1+m$. First, $r(x_k)<r(p'_{i_1})$, for $0\le k\le m-1$, implies
$x_j\ne p'_{i_1+k}, p'_{i_1+k}-e$, for $k\le j\le m-1$. We also have
$q'_{i_1+j}\ne p'_{i_1+k}, p'_{i_1+k}-e$, for $0\le j\le k-1$.
This follows from ($\dot e$)
when $0\le k\le m-1$ or $i_1+m=b$, since $p'_{i_1+k}=p'_{i_1}-ke$ in
these cases, and from $r(p'_{i_1+m})>r(p'_{i_1})\ge r(q'_{i_1+j})$ when $i_1+m$ is a reset point.
Second, if $i_1+k<b$ then $U(J_{i_1+k})=U(K_{i_1+k})$.
In fact, if $p'_{i_1+k}=p'_{i_1+k-1}-e$ then
$U(J_{i_1+k})=U(K_{i_1+k})=\{p'_{i_1+k}\}$ on runners $1,\dots,r(p'_{i_1})$
by ($\dot b$), ($\dot c$), ($\dot d$), $q'_{i_1+j}\ne p'_{i_1+k}, p'_{i_1+k}-e$,
for $0\le j\le k-1$, and $x_j\ne p'_{i_1+k}, p'_{i_1+k}-e$, for $k\le j\le m-1$.
If $p'_{i_1+k}>p'_{i_1}$ then $U(J_{i_1+k})=U(K_{i_1+k})=\emptyset$
on runners $1,\dots,r(p'_{i_1})$ by ($\dot c$), ($\dot d$) and $r(p'_{i_1+k})>r(p'_{i_1})$.
$J_{i_1+k}=K_{i_1+k}$ on runners $r(p'_{i_1})+1,\dots,e-1$ by ($\dot b$) and ($\dot f$),
and $U(J_{i_1+k})=U(K_{i_1+k})=\emptyset$ on runner $e$ by ($\dot a$), ($\dot b$)
and ($\dot f$).
Thus, $U(J_{i_1+k})=U(K_{i_1+k})$ if $i_1+k<b$.
If $i_1+m=b$ then the same proof shows that $U(K_b)=\emptyset$.
(a) to (h) for $i_1\le i\le i_1+m-1$ or $i=i_1+m$ when $i_1+m=b$
clearly follows from ($\dot a$) to ($\dot h$). When $i_1+m$ is a reset point,
$U(J_{i_1+m})=U(K_{i_1+m})$ implies (c) and (d) for $i=i_1+m$.
The other parts of (a) to (h) are obvious.

Now we prove the claims. The claims hold when $k=0$. Suppose that the claims hold for $k$
such that $0\le k\le m-1$. Thus $p'_{i_1}-je\not\in J_{i_1+k}$, for
$k+1\le j\le m$, $i_1+k<b$ and $p'_{i_1+k}=p'_{i_1}-ke$.
If $k+1\le m-1$ then $p'_{i_1}-je\not\in J_{i_1+k+1}$, for
$k+2\le j\le m$, and $p'_{i_1}-(k+1)e\in U(J_{i_1+k+1})$.
Hence, $i_1+k+1<b$ and $p'_{i_1+k+1}<p'_{i_1+k}$ implies
$p'_{i_1+k+1}=p'_{i_1}-(k+1)e$.

As $p'_{i_1+k}-e=p'_{i_1}-(k+1)e<x<x_k$, for $x\in J_{i_1+k}$,
implies $x\not\in W(J_{i_1+k})$ by ($\dot a$), ($\dot g$) and ($\dot h$), we have
$x_k\le q'_{i_1+k}\le p'_{i_1+k}$. As $\ell(p'_{i_1+k})=k$ and
$\ell(x_k)=k$, this implies
$$
\ell(q'_{i_1+k})=k\;\;\text{and}\;\;
r(q'_{i_1+k})\le r(p'_{i_1})<e.
$$
Hence, ($\dot a$) and ($\dot e$) for $k+1$ follow.

$\ell(x_k)=\ell(q'_{i_1+k})$
and $x_k\le q'_{i_1+k}$ imply $r(x_k)\le r(q'_{i_1+k})$. As
$r(x_{k+1})\le r(x_k)$, we have $r(x_{k+1})\le r(q'_{i_1+k})$.
If $k=0$ then we have proved ($\dot f$) for $k+1$.
If $k\ge1$ then we have to show $r(q'_{i_1+k})\le r(q'_{i_1+k-1})$.
Note that we have either $r(q'_{i_1+k-1})=r(p'_{i_1+k-1})$ or $r(q'_{i_1+k-1})<r(p'_{i_1+k-1})$
by ($\dot f$). If $r(q'_{i_1+k-1})=r(p'_{i_1+k-1})$ then we have
$r(q'_{i_1+j})=r(p'_{i_1+j})$ and $\ell(q'_{i_1+j})=\ell(p'_{i_1+j})$, for
$0\le j\le k-1$, by ($\dot f$).
This implies that $q'_{i_1+j}=p'_{i_1+j}$, for $0\le j\le k-1$.
Thus, $J_{i_1+k-1}$ is obtained from $J_{i_1}$ by
moving the bead $p'_{i_1}$ up to $p'_{i_1+k-1}=q'_{i_1+k-1}$.
Hence,
$$
q'_{i_1+k-1}-e\in J_{i_1+k}\;\;\text{and}\;\;
q'_{i_1+k-1}\not\in J_{i_1+k}.
$$
If $r(q'_{i_1+k-1})<r(p'_{i_1+k-1})$ then $q'_{i_1+k-1}\in J_{i_1+k-1}$
implies $q'_{i_1+k-1}-e\in J_{i_1+k-1}$ by ($\dot c$) for $k-1$.
Thus,
$q'_{i_1+k-1}-e\in J_{i_1+k}$ and $q'_{i_1+k-1}\not\in J_{i_1+k}$
follow again. Therefore, $q'_{i_1+k-1}-e\in W(J_{i_1+k})$ and we conclude
$$
q'_{i_1+k}\le q'_{i_1+k-1}-e.
$$
Then, $\ell(q'_{i_1+k})=\ell(q'_{i_1+k-1}-e)$ implies
$r(q'_{i_1+k})\le r(q'_{i_1+k-1}-e)=r(q'_{i_1+k-1})$.
We have proved ($\dot f$) for $k+1$. As $r(q'_{i_1+k})\le r(p'_{i_1})$,
($\dot g$) for $k+1$ also follow.

Now $U(J_{i_1+k})=U(K_{i_1+k})$ implies $p'_{i_1+k}\in U(K_{i_1+k})$ and
$p''_{i_1+k}=p'_{i_1+k}$. Hence it is clear that
($\dot c$) and ($\dot d$) for $k+1$ hold.

To show that $q''_{i_1+k}=x_k$, first suppose that
$$
p'_{i_1+k}-e=p'_{i_1}-(k+1)e<x<j_0-(k+1)e,
$$
for $x\in K_{i_1+k}$. Then $x\in J_{i_1+k}$ by ($\dot b$) and ($\dot f$), and
$x+e\in J_{i_1+k}$ by ($\dot g$). Using ($\dot b$) and ($\dot f$) again, we have
$x+e\in K_{i_1+k}$ and $x\not\in W(K_{i_1+k})$. If
$$
j_0-(k+1)e<x<x_k,
$$
then $x\not\in K_{i_1+k}$ by ($\dot h$), and $x\not\in W(K_{i_1+k})$
again. We have proved that $q''_{i_1+k}\ge x_k$. To see that
$q''_{i_1+k}=x_k$, it remains to show $x_k\in W(K_{i_1+k})$.

Note that
$\ell(x_k)=\ell(p'_{i_1+k})$ and $r(x_k)\le r(p_{i_1})=r(p'_{i_1+k})$
imply that
$$
p'_{i_1+k}-e<x_k\le p'_{i_1+k}.
$$
As $x_k\in K_{i_1+k}$,
$x_k\in W(K_{i_1+k})$ follows when $x_k=p'_{i_1+k}$. If $x_k< p'_{i_1+k}$,
we have to show $x_k+e\not\in K_{i_1+k}$. It is clear when $k=0$.
Suppose $k\ge1$ and $x_k+e\in K_{i_1+k}$. Thus ($\dot h$) implies
$r(x_k+e)\ge r(q'_{i_1+k-1})$. On the other hand, ($\dot f$) implies
$\ell(x_k+e)=\ell(q'_{i_1+k-1})$ and $r(x_k+e)\le r(q'_{i_1+k-1})$.
Hence $x_k+e=q'_{i_1+k-1}\in J_{i_1+k-1}$ follows.
As $r(x_k+e)\le r(p'_{i_1+k-1})$,
we have either $x_k\in J_{i_1+k-1}$ or
$x_k+e=p'_{i_1+k-1}$ by ($\dot c$) for $k-1$.
By ($\dot b$) for $k-1$, $x_k\in J_{i_1+k-1}$ does
not occur. $x_k+e=p'_{i_1+k-1}$ implies
$x_k=p'_{i_1+k}\in J_{i_1+k}$, which contradicts ($\dot b$).
Therefore, $x_k+e\not\in K_{i_1+k}$.
We have proved that
$x_k\in W(K_{i_1+k})$, and $q''_{i_1+k}=x_k$ follows.
In other words, we have proved
$$
K_{i_1+k+1}=\left(K_{i_1+k}\setminus\{x_k\}\right)\sqcup\{p'_{i_1+k+1}\}.
$$
By $J_{i_1+k+1}\sqcup\{q'_{i_1+k}\}=J_{i_1+k}\sqcup\{p'_{i_1+k+1}\}$
and (b), $K_{i_1+k+1}$ is equal to
$$
\left(J_{i_1+k+1}\setminus\{j_0,\dots,j_0-me\}\right)
\sqcup\{q'_{i_1},\dots,q'_{i_1+k-1},q'_{i_1+k},x_{k+1},\dots, x_{m-1}\}.
$$
We have proved ($\dot b$) for $k+1$.

Finally, to prove ($\dot h$) for $k+1$, we have to show that
$x\not\in J_{i_1+k+1}$ and $x\not\in K_{i_1+k+1}$ when
$\ell(x)=k$ and $r(x)<r(q'_{i_1+k})$.
If $x\in J_{i_1+k+1}$ then $x\ne p'_{i_1+k}-e, q'_{i_1+k}$ implies
$x\in J_{i_1+k}$ and $x+e\in J_{i_1+k}$ by $x\not\in W(J_{i_1+k})$.
However, $x+e\not\in J_{i_1+k}$ if $k=0$, and
if $k\ge1$ then $\ell(x+e)=k-1$ and $r(x+e)<r(q'_{i_1+k})\le r(q'_{i_1+k-1})$
imply $x+e\not\in J_{i_1+k}$ by ($\dot h$).
We have proved $x\not\in J_{i_1+k+1}$. If $x\in K_{i_1+k+1}$ then
$x\in J_{i_1+k+1}$ or $x=q'_{i_1+k}$ by ($\dot b$) for $k+1$.
As both do not occur, $x\not\in K_{i_1+k+1}$.

We have proved the desired claims for $i_1\le i\le i_1+m$.
Note that we have also proved that $i_1,\dots,i_1+m\in U$.

Define $m'\ge m$ by $m'=i_2-i_1$ if $i_1+m,\dots,i_2\in U$, and by
$$
i_1+m,i_1+m+1,\dots,i_1+m'\in U\;\;\text{and}\;\;
i_1+m'+1\not\in U,
$$
otherwise. We want to show that the claims hold for
$i_1+m\le i\le i_1+m'+1$. To do this,
we show, for $m\le k\le m'+1$, that $p''_{i_1+j}=p'_{i_1+j}$,
for $0\le j\le k-1$, and
\begin{itemize}
\item[(\"a)]
$j_0-e\mathbb Z_{\ge0}\subset J_{i_1+k}$ and $\max J_{i_1+k}=j_0$.
\item[(\"b)]
$K_{i_1+k}=\left(J_{i_1+k}\setminus\{j_0,j_0-e,\dots,j_0-ke\}\right)
\sqcup\{q'_{i_1},\dots,q'_{i_1+k-1}\}$.
\item[(\"c)]
If $x\in J_{i_1+k}$ is such that $r(x)\le r(p'_{i_1+k})$ then
$x\not\in U(J_{i_1+k})$ unless $x=p'_{i_1+k}$.
\item[(\"d)]
If $x\in K_{i_1+k}$ is such that $r(x)\le r(p'_{i_1+k})$ then
$x\not\in U(K_{i_1+k})$ unless $x=p'_{i_1+k}$.
\item[(\"e)]
$\ell(q'_{i_1+j})=j$, for $0\le j\le k-1$.
\item[(\"f)]
$r(p'_{i_1})\ge r(q'_{i_1})\ge\cdots\ge r(q'_{i_1+k-1})$.
\item[(\"g)]
If there exists $x\in J_{i_1+k}$ such that $1\le\ell(x)\le k$ and
$r(p'_{i_1})<r(x)<e$ then
$(x+e\mathbb Z)\cap J_{\le j_0}\subset J_{i_1+k}$.
\item[(\"h)]
If $j_0-(j+1)e<x<q'_{i_1+j}$, for some $0\le j\le k-1$,
then $x\not\in J_{i_1+k}$ and $x\not\in K_{i_1+k}$.
\end{itemize}
By the same argument as before, these claims imply the desired
claims for $i_1+m\le i\le i_1+m'+1$.
Suppose that the claims hold for $k$ such that $m\le k\le m'$.
Thus $p'_{i_1+k}=p_{i_1}-ke$ and, by definition, $i_1+k<b$.
$i_1+k\in U$ implies $\ell(q'_{i_1+k})=\ell(p'_{i_1+k})=k$ and
$r(q'_{i_1+k})\le r(p'_{i_1})<e$. Thus (\"a) and (\"e) for $k+1$ follow.
If $m=0$ and $k=m$ then (\"f) for $k+1$ is clear.
Otherwise, $k\ge1$ and we have either $r(q'_{i_1+k-1})=r(p'_{i_1+k-1})$
or $r(q'_{i_1+k-1})<r(p'_{i_1+k-1})$ by $i_1+k-1\in U$.
Now the rest of the proof is entirely similar to the previous one.
The only difference is that we prove
$q''_{i_1+k}=j_0-(k+1)e$. To prove this, suppose that
$p'_{i_1+k}-e<x\le p'_{i_1+k}$. As $J_{i_1+k}=K_{i_1+k}$ on
runners $r(p'_{i_1})+1,\dots,e-1$, $x\in W(K_{i_1+k})$ implies
$x\ge j_0-(k+1)e$. As $j_0-(k+1)e\in W(K_{i_1+k})$ by (\"b),
we have $q''_{i_1+k}=j_0-(k+1)e$.

Note that we have also proved $m_{i_1+k}=k$, for $m\le k\le m'+1$.

If $i_1+m'+1=b$ then we have finished the proof. Suppose $i_1+m'+1<b$ and
$i_1+m'=i_2$. As $p'_{i_2+1}-e\not\in J_{i_2+1}$ implies
$p'_{i_2+1}-e\not\in J_{i_1+m'}$,
(\"g) for $k=m'$ implies $p'_{i_2+1}-je\not\in J_{i_1+m'}$, and thus
$p'_{i_2+1}-je\not\in J_{i_2+1}$, for
$1\le j\le m'$. Let $x=p'_{i_2+1}-(m'+1)e$. Then $x+e\not\in J_{i_1+m'}$ and
$p'_{i_1+m'}-e<x<p'_{i_1+m'}$. Thus, if $x\in J_{i_1+m'}$ then
$x\in W(J_{i_1+m'})$ and the minimality of $q'_{i_1+m'}$ is contradicted.
Therefore, $p'_{i_2+1}-je\not\in J_{i_2+1}$, for
$1\le j\le m'+1=m_{i_2+1}$.

To complete the proof of Lemma \ref{comparison}, we consider the
case $i_1+m'<i_2$.
Write $x'_k=q'_{i_1+k}$, for $0\le k\le m'$. We have
$$
r(p'_{i_1})\ge r(x'_0)\ge\cdots\ge r(x'_{m'})
$$
by (\"f) for $k=m'+1$. We show $i\not\in U$ and the claims (A) to (C) below,
for $i_1+m'+1\le i\le i_2+1$. They hold when $i=i_1+m'+1$. Suppose that
the claims hold for $i$ such that $i_1+m'+1\le i\le i_2$.
Thus $i\le i_2\le b-1$, $\ell(p'_i)\ge m'+1$, $r(p'_i)=r(p'_{i_1})$,
$i\not\in U$ and
\begin{itemize}
\item[(A)]
$j_0-e\mathbb Z_{\ge0}\subset J_i$ and $\max J_i=j_0$.
\item[(B)]
$K_i=\left(J_i\setminus\{j_0,j_0-e,\dots,j_0-(m'+1)e\}\right)
\sqcup\{x'_0,\dots,x'_{m'}\}$.
\item[(C)]
If $x\in\mathbb Z$ is such that $0\le r(x)\le r(p'_{i_1})$ then
$x\not\in U(J_i)$ and $x\not\in U(K_i)$ unless $x=p'_i$.
\end{itemize}
Note that (B) implies $p'_i\in U(K_i)$, and
(A), (B), (C) imply $U(J_i)=U(K_i)$ and $p''_i=p'_i$.

As $i\not\in U$, we have $p'_{i_1}-(i+1)e<q'_i<j_0-(i+1)e$ and
$q'_i<p'_i$ implies $q_i'-e\in J_i$ and $q'_i-e\in J_{i+1}$. Thus,
if $i+1\le i_2$ then
$q'_i-e\in W(J_{i+1})$ and it follows that
$$
p'_{i+1}-e<q'_{i+1}\le j_0-(i+2)e.
$$
Hence, $i+1\not\in U$.

(B) implies $q'_i\in W(K_i)$. Thus $q''_i\le q'_i$. Then,
$$
p'_i-e=p''_i-e<q''_i<j_0-(i+1)e
$$
and (B) implies $q''_i\in W(J_i)$, which proves $q''_i= q'_i$.
Therefore, (A), (B), (C) for $i+1$ follow.

We have proved $U(J_i)=U(K_i)$, $p'_i\ne x'_k$, for $0\le k\le m'$,
and (a), (b), (c), (d), (e), (f), for $i_1+m'+1\le i\le i_2+1$.
Now, $\ell(q'_i)=\ell(p'_i)+1\ge m'+2$ implies that
we do not touch the layer levels smaller than or equal to
$m'+1$ on runners $r(p'_{i_1})+1,\dots,e-1$.
Thus (g), for $i_1+m'+1\le i\le i_2+1$, follow.
Similarly, we do not touch the layer levels smaller than or equal to
$m'$ on runners $1,\dots, r(p'_{i_1})$.
Thus (h), for $i_1+m'+1\le i\le i_2+1$, follows.

If $i_2+1=b$ then we have finished the proof. Suppose $i_2+1<b$.
Then $i_2+1$ is a reset point and $r(p'_{i_1})<r(p'_{i_2+1})<e$. Thus,
(g) for $i=i_2$ implies
$p'_{i_2+1}-je\not\in J_{i_2+1}$, for $1\le j\le m_{i_2+1}$, since
$\ell(q'_{i_2})\ge m'+2$,
$m_{i_2+1}=m_{i_2}=m'+1$ and $p'_{i_2+1}-e\not\in J_{i_2+1}$.

Now, the induction on $i$ works and we have proved the claims for $0\le i\le b$.
\end{proof}

By Lemma \ref{comparison}, there exists $m\ge 1$ such that we may write
\begin{equation*}
K_b=\left(J_b\setminus\{j_0-ke \mid 0\le k\le m\}\right)
\cup\{x_k \mid 0\le k\le m-1\},
\end{equation*}
where $r(x_k)\le r(p'_b)$, for $0\le k\le m-1$. Consider
$K_b\cup\{j_0\}$. As $j_1\in K_b$ and $j_0-j_1\le e$, the
partition associated with $K_b\cup\{j_0\}$ is $e$-restricted. Note
that $U(J_b)=\emptyset$ and $U(K_b)=\emptyset$. Hence,
$U(K_b\cup\{j_0\})=\{j_0\}$ and explicit computation of
$\down^k(K_b\cup\{j_0\})$, for $k\ge0$, by using (a) to (h), shows
that we obtain $\down^{k+1}(K_b\cup\{j_0\})$ from
$\down^k(K_b\cup\{j_0\})$ by moving $x_k$ to $j_0-(k+1)e$, for
$0\le k\le m-1$. Thus we end up with $\base(K_b\cup\{j_0\})=J_b$.
Therefore,
$$
\base(\base(K)\cup\{j_0\})=\base(K_b\cup\{j_0\})=J_b=\base(J).
$$
We have now proved the following proposition.

\begin{prop}
\label{extension lemma}
Let $\lambda\in B(\Lambda_m)$ and $J$ the corresponding set of beta numbers
of charge $m$. Set $K=J_{\le j_0-1}$, where $j_0=\max J$.
Then, the partition associated with $\base(K)\cup\{j_0\}$ is
$e$-restricted and we have
$$
\base(J)=\base(\base(K)\cup\{j_0\}).
$$
\end{prop}

Let $\lambda\in B(\Lambda_m)$ and $J$ the corresponding set of
beta numbers of charge $m$. We delete the first row from $\lambda$
and we denote the resulting partition by $\mu$. Assume that
$\base(\mu)$ is already computed.
Then it is easy to compute $\base(\lambda)$ by using the above
proposition. It gives us an efficient inductive definition of base and
it is possible to generalize main results in section 8
to other types $A_{2n}^{(2)}$ and $D_{n+1}^{(2)}$.

\begin{cor}
\label{cor to ext lemma}
Let $\lambda\in B(\Lambda_m)$ and $J$ the corresponding set of
beta numbers of charge $m$. Let $j_0>\cdots>j_r$ be the largest
$r+1$ members of $J$. Define $J_{r+1}=J_{\le j_r-1}$ and
$J_k=\base(J_{k+1})\cup\{j_k\}$, for $k=r,\dots,0$. Then
the partition associated with $J_k$ is $e$-restricted and
$\base(J)=\base(J_0)$.
\end{cor}
\begin{proof}
We show by downward induction on $k$ that
$\max J_k=j_k$ and $\base(J_{\le j_k})=\base(J_k)$.
When $k=r+1$ there is nothing to prove. Suppose that the equations hold for
$k+1$. Then, $j_{k+1}\in\base(J_{k+1})$ and $j_k-j_{k+1}\le e$ imply
that the partition associated with $J_k$ is $e$-restricted.
Now, by Proposition \ref{extension lemma} and
the induction hypothesis,
\begin{equation*}
\begin{split}
\base(J_{\le j_k})&=\base(J_{\le j_{k+1}}\cup\{j_k\})
=\base(\base(J_{\le j_{k+1}})\cup\{j_k\})\\
&=\base(\base(J_{k+1})\cup\{j_k\})=\base(J_k).
\end{split}
\end{equation*}
Thus, $\base(J)=\base(J_0)$ follows.
\end{proof}

\section{Base Theorem}

Let $\lambda\in B(\Lambda_m)$ and $J$ the set of beta numbers of
charge $m$. Define
$$
M_i(\lambda)=M_i(J)=\max\{x\in J \mid x+e\mathbb Z=i\}.
$$

\begin{lemma}
\label{basic lemma}
Let $\lambda\in B(\Lambda_m)$.
\item[(1)]
If $M_i(\lambda)\le M_{i+1}(\lambda)$ then
$M_i(\down(\lambda))\le M_{i+1}(\down(\lambda))$.
In particular, if $M_i(\lambda)\le M_{i+1}(\lambda)$
then $s_i\base(\lambda)\le\base(\lambda)$.
\item[(2)]
If $\lambda$ is an $s_i$-core and $s_i\lambda\ge\lambda$ then
\begin{itemize}
\item[(i)]
$\down(\lambda)$ and $\down(s_i\lambda)$ are $s_i$-cores,
\item[(ii)]
$\down(s_i\lambda)=s_i\down(\lambda)$.
\end{itemize}
\item[(3)]
Suppose that $\lambda$ is an $s_i$-core and $s_i\lambda\ge\lambda$.
\begin{itemize}
\item[(a)]
If $s_i\base(\lambda)>\base(\lambda)$ then
$\base(s_i\lambda)=s_i\base(\lambda)>\base(\lambda)$.
\item[(b)]
If $s_i\base(\lambda)\le\base(\lambda)$ then
$\base(s_i\lambda)=\base(\lambda)$.
\end{itemize}
\item[(4)]
Suppose that $\lambda$ has an addable $i$-node
on the first row, and that if we delete the first row then
the resulting partition, which we denote by $\mu$, is an $e$-core.
\begin{itemize}
\item[(a)]
Suppose that $s_i\mu\ge\mu$. Then
$\base(\f_i^{\varphi_i(\lambda)-1}\lambda)=\base(\lambda)<s_i\base(\lambda)$.
\item[(b)]
Suppose that $s_i\mu\le\mu$. Then $\varphi_i(\lambda)=1$ and
$$
\base(\f_i^{\varphi_i(\lambda)}\lambda)=
s_i\base(\lambda)>\base(\lambda).
$$
\end{itemize}
\end{lemma}
\begin{proof}
(1) Let $J$ be the corresponding set of beta numbers of
charge $m$, and define $p'$ and $q'$ as in the definition
of $\down(J)$. Note that adding the bead $p'-e$ does not affect
$M_i(\lambda)$ or $M_{i+1}(\lambda)$ because if $q'<p'$ then
there exists a larger element $p'$ in $J$.
Thus it suffices to study the effect of moving $q'$.

First suppose that $q'\ne M_{i+1}(\lambda)$. Then
$$
M_{i+1}(\down(\lambda))=M_{i+1}(\lambda)\ge M_i(\lambda)
\ge M_i(\down(\lambda)).
$$
The last inequality is an equality when $q'\ne M_i(\lambda)$.
$M_i(\down(\lambda))\le M_{i+1}(\down(\lambda))$ holds.

Second suppose that $q'=M_{i+1}(\lambda)$.
In particular, $q'$ is on the $(i+1)^{th}$ runner.
Note that $p'$ cannot be on the $i^{th}$ runner:
if so then $p'\ge q'$ would imply $p'\ge q'+e-1$ and
$$
M_i(\lambda)\ge p'\ge q'+e-1>q'=M_{i+1}(\lambda),
$$
which contradicts our assumption.

We shall show $q'-1\not\in J$. Suppose on the contrary that
$q'-1\in J$. If $q'=p'$ then $q'-1>p'-e$ and
$q'-1+e\not\in J$ by $M_i(\lambda)\le M_{i+1}(\lambda)=q'$.
This implies that $q'-1\in W(J)$, which contradicts
$q'=\min W(J)$. If $q'<p'$ then we also have
$q'-1>p'-e$ and $q'-1+e\not\in J$, since
$p'-e=q'-1$ would imply that $p'$ is on the $i^{th}$ runner.
Hence we reach the contradiction $q'-1\in W(J)$ again.
We have proved that $q'-1\not\in J$.

Now we are ready to prove that
$M_i(\down(\lambda))\le M_{i+1}(\down(\lambda))$.
Since $q'-1\not\in J$ and $M_i(\lambda)\le q'$, we have
$M_i(\down(\lambda))\le q'-1-e$.

Suppose that $q'<p'$. Since $p'=\min U(J)$, we have
$q'-e\in J$ and
$$
M_{i+1}(\down(\lambda))=q'-e>M_i(\down(\lambda))
$$
follows. If $q'=p'$ then we have
$M_{i+1}(\down(\lambda))=q'-e$ by definition, and the result again follows.
We have proved the first half of the claim.

Now, define a decreasing sequence of partitions
$$
\lambda=\lambda^{(0)}>\cdots>\lambda^{(k)}>
\cdots>\lambda^{(s)}=\base(\lambda)
$$
by $\down(\lambda^{(k)})=\lambda^{(k+1)}$, for $0\le k<s$.
Then, by repeated use of the first half of the claim, we have
$M_i(\base(\lambda))\le M_{i+1}(\base(\lambda))$.
This implies that the $e$-core $\base(\lambda)$ does
not have an addable $i$-node. Thus
$s_i\base(\lambda)\le\base(\lambda)$.

(2) Note that $s_i\lambda$ is an $s_i$-core
by Lemma \ref{lemma for cores}(3).
Let $J$ be the set of beta numbers of charge $m$ associated
with $\lambda$.
As $\lambda$ is an $s_i$-core, $p'=\min U(J)$ cannot be
on the $i^{th}$ or the $(i+1)^{th}$ runners. Since
$s_i\lambda$ is obtained from $\lambda$ by the rule
given in Lemma \ref{lemma for cores}, both contain
$p'$, that is, $p'=\min U(J)=\min U(s_iJ)$.
Let $q'=\min W(J)$. Then $q'$ for $s_i\lambda$ is given by
$$
\min W(s_iJ)=
\begin{cases}
q'+1=M_i(\lambda)+1=M_{i+1}(s_i\lambda)
\quad&\text{if $q'+e\mathbb Z=i$.}\\
q'-1=M_{i+1}(\lambda)-1=M_i(s_i\lambda)
\quad&\text{if $q'+e\mathbb Z=i+1$.}\\
q'\quad&\text{otherwise.}
\end{cases}
$$
To see this, note that if $x<q'$ is located on a runner different
from the $i^{th}$ and the
$(i+1)^{th}$ runners and if $x$
satisfies $x\in J$, $p'-e<x$ and $x+e\not\in J$,
then $x\not\in W(J)$, which implies $x\not\in W(s_iJ)$.

Suppose that $q'+e\mathbb Z=i$. Then $q'<p'$ and
$q'=M_i(\lambda)\ge M_{i+1}(\lambda)$ implies that
$$
M_{i+1}(\lambda)-1\le q'-e\le p'-e.
$$
Thus $M_{i+1}(\lambda)-1\not\in W(s_iJ)$ and there is no element
of $W(s_iJ)$ on the $i^{th}$ runner. On the other hand, we have
$q'+1\in W(s_iJ)$ and $\min W(s_iJ)=q'+1$ follows.

If $q'+e\mathbb Z=i+1$ then $q'<p'$, $q'=M_{i+1}(\lambda)$ and
$\min W(s_iJ)=q'-1$ is easy to see.
Similarly, we have $\min W(s_iJ)=q'$ otherwise. Now it is clear that
$$
(i)\;\;\text{$\down(\lambda)$ and $\down(s_i\lambda)$ are $s_i$-cores},\quad
(ii)\;\;\down(s_i\lambda)=s_i\down(\lambda).
$$

(3) To prove (a) and (b), we consider two decreasing sequences
\begin{equation*}
\begin{split}
\lambda&=\lambda^{(0)}>\cdots>\lambda^{(k)}>
\cdots>\lambda^{(s)}=\base(\lambda)\\
s_i\lambda&=\mu^{(0)}>\cdots>\mu^{(k)}>
\cdots>\mu^{(t)}=\base(s_i\lambda)
\end{split}
\end{equation*}
where $\down(\lambda^{(k)})=\lambda^{(k+1)}$, for
$0\le k<s$, and
$\down(\mu^{(k)})=\mu^{(k+1)}$, for $0\le k<t$.

(a) We prove by induction on $k$ that
$$
(i)\;\;\text{$\lambda^{(k)}$ and $\mu^{(k)}$ are $s_i$-cores},\quad
(ii)\;\;\mu^{(k)}=s_i\lambda^{(k)},\quad
(iii)\;\;s_i\lambda^{(k)}\ge\lambda^{(k)},
$$
for $0\le k\le\min(s,t)$. This implies the desired result.
In fact, as $\lambda^{(k)}$ is an $e$-core if and only if
$\mu^{(k)}=s_i\lambda^{(k)}$ is an $e$-core by Lemma \ref{lemma for cores}(3),
we must have $s=t$.
Thus $\base(s_i\lambda)=s_i\base(\lambda)$ follows.

If $k=0$ then the claim holds by the hypothesis. Suppose that
the claim holds for $k$. Then (2) implies
\begin{itemize}
\item[(i)]
$\lambda^{(k+1)}=\down(\lambda^{(k)})$ and
$\mu^{(k+1)}=\down(\mu^{(k)})=\down(s_i\lambda^{(k)})$ are $s_i$-cores,
\item[(ii)]
$\mu^{(k+1)}=\down(\mu^{(k)})=\down(s_i\lambda^{(k)})
=s_i\down(\lambda^{(k)})=s_i\lambda^{(k+1)}$.
\end{itemize}

If $M_i(\down(\lambda^{(k)}))<M_{i+1}(\down(\lambda^{(k)}))$ then
(1) implies that $s_i\base(\lambda)\le\base(\lambda)$, contradicting
the hypothesis. Thus,
$M_i(\down(\lambda^{(k)}))\ge M_{i+1}(\down(\lambda^{(k)}))$ and
this and (i) imply
\begin{itemize}
\item[(iii)]
$s_i\lambda^{(k+1)}\ge\lambda^{(k+1)}$.
\end{itemize}

(b) If $s_i\lambda^{(0)}=\lambda^{(0)}$
then the result is obvious.
Suppose that $s_i\lambda^{(0)}>\lambda^{(0)}$. As
$s_i\lambda^{(t)}\le\lambda^{(t)}$,
the same induction argument as in (a) proves that
there exists the maximal $1\le k_0\le t$ such that
$$
(i)\;\;\text{$\lambda^{(k)}$ and $\mu^{(k)}$ are $s_i$-cores},\quad
(ii)\;\;\mu^{(k)}=s_i\lambda^{(k)},\quad
(iii)\;\;s_i\lambda^{(k)}>\lambda^{(k)},
$$
for $0\le k\le k_0-1$. Then (i) for $k=k_0-1$ and
$s_i\lambda^{(k_0-1)}>\lambda^{(k_0-1)}$ imply
$$
M_i(\lambda^{(k_0-1)})>M_{i+1}(\lambda^{(k_0-1)}).
$$
Applying (2) once more, we also have
$$
(i)\;\;\text{$\lambda^{(k_0)}$ is an $s_i$-core},\quad
(ii)\;\;\mu^{(k_0)}=s_i\lambda^{(k_0)}.
$$

Let $J$ be the set of beta numbers of charge $m$ associated
with $\lambda^{(k_0-1)}$. Then,
$\lambda^{(k_0-1)}$ and $\mu^{(k_0-1)}$ both have
$p'=\min U(J)=\min U(s_iJ)$. Consider $q'=\min W(J)$.
Assume that $q'+e\mathbb Z\ne i$. Then
$$
M_{i+1}(\lambda^{(k_0-1)})<M_i(\lambda^{(k_0-1)})=M_i(\lambda^{(k_0)})
\le M_{i+1}(\lambda^{(k_0)})
$$
and $M_{i+1}(\lambda^{(k_0)})$ is either $M_{i+1}(\lambda^{(k_0-1)})-e$
or $M_{i+1}(\lambda^{(k_0-1)})$. In either case, we have a contradiction, and
we conclude that $q'+e\mathbb Z=i$. Then
$$
M_{i+1}(\lambda^{(k_0-1)})<M_i(\lambda^{(k_0-1)})
=M_i(\lambda^{(k_0)})+e\le M_{i+1}(\lambda^{(k_0)})+e
$$
and $M_{i+1}(\lambda^{(k_0)})+e=M_{i+1}(\lambda^{(k_0-1)})+e$.

As $M_{i+1}(\lambda^{(k_0-1)})<M_i(\lambda^{(k_0-1)})$
implies $M_{i+1}(\lambda^{(k_0-1)})+e-1\le M_i(\lambda^{(k_0-1)})$
and $M_i(\lambda^{(k_0)})\ne M_{i+1}(\lambda^{(k_0)})$,
we have
$M_i(\lambda^{(k_0)})+1= M_{i+1}(\lambda^{(k_0)})$.
Since $\lambda^{(k_0)}$ is also an $s_i$-core, this implies
$s_i\lambda^{(k_0)}=\lambda^{(k_0)}$.
Hence, we have $\mu^{(k_0)}=\lambda^{(k_0)}$, which implies
$\base(s_i\lambda)=\base(\lambda)$.

(4) (a) Since $s_i\mu\ge\mu$, $\lambda$ does not
have a removable $i$-node.
Let $J$ be the set of beta numbers associated with
$\f_i^{\varphi_i(\lambda)-1}\lambda$, and let $K$
be the set of beta numbers associated with $\lambda$.
We have $\max J=j_0=\max K$. Then
\begin{itemize}
\item[(i)]
By deleting the first row from $\f_i^{\varphi_i(\lambda)-1}\lambda$,
we obtain $\f_i^{max}\mu=s_i\mu$.
\item[(ii)]
The set of beta numbers associated with
$\f_i^{\varphi_i(\lambda)-1}\lambda$
is $s_i(J\setminus\{j_0\})\cup\{j_0\}$.
\end{itemize}
If $\varphi_i(\lambda)=1$ then the claim
$\base(\f_i^{\varphi_i(\lambda)-1}\lambda)=\base(\lambda)$
is obvious.
Assume that $\varphi_i(\lambda)>1$. Then
$\f_i^{\varphi_i(\lambda)-1}\lambda$ has
both an addable $i$-node and a removable $i$-node, thus
it cannot be an $e$-core. This implies $U(J)\ne\emptyset$ and
we have $U(J)=\{j_0\}$, $p'=\min U(J)=j_0$.

Note that the abacus displays of $J$ and $K$ have the following form
by (i) and (ii) above.

\begin{equation*}
J:\;\;
\begin{aligned}
\cdots& \hphantom{=}& \times \;     & \times      & \hphantom{=}&\cdots\\
\cdots& \hphantom{=}& \times \;     & \times      & \hphantom{=}&\cdots\\
\cdots& \hphantom{=}& \hphantom{=}\;& \times      & \hphantom{=}&\cdots\\
\cdots& \hphantom{=}& \hphantom{=}\;& \times      & \hphantom{=}&\cdots\\
\cdots& \hphantom{=}& \hphantom{=}\;& \times      & \hphantom{=}&\cdots\\
\cdots& \hphantom{=}& \hphantom{=}\;& \hphantom{=}& \hphantom{=}&\cdots\\
\cdots& \hphantom{=}& j_0         \;& \hphantom{=}& \hphantom{=}&\cdots
\end{aligned}
\qquad K:\;\;
\begin{aligned}
\cdots& \hphantom{=}& \times \;     & \times      & \hphantom{=}&\cdots\\
\cdots& \hphantom{=}& \times \;     & \times      & \hphantom{=}&\cdots\\
\cdots& \hphantom{=}& \times \;     & \hphantom{=}& \hphantom{=}&\cdots\\
\cdots& \hphantom{=}& \times \;     & \hphantom{=}& \hphantom{=}&\cdots\\
\cdots& \hphantom{=}& \times \;     & \hphantom{=}& \hphantom{=}&\cdots\\
\cdots& \hphantom{=}& \hphantom{=}\;& \hphantom{=}& \hphantom{=}&\cdots\\
\cdots& \hphantom{=}& j_0         \;& \hphantom{=}& \hphantom{=}&\cdots
\end{aligned}
\end{equation*}

Thus, there exists $k_0$ such that, for
$0\le k\le k_0$, $\down^{k+1}(\f_i^{\varphi_i(\lambda)-1}\lambda)$
and $\down^{k+1}(\lambda)$ are obtained from
$\down^k(\f_i^{\varphi_i(\lambda)-1}\lambda)$
and $\down^k(\lambda)$ by moving the maximal element of
the $j^{th}$ runner, for some $j\ne i,i+1$, to the $i^{th}$ runner,
respectively. Note that $j$ is the same for
$\f_i^{\varphi_i(\lambda)-1}\lambda$ and $\lambda$ in each step
$k$. At $k=k_0$, we reach the following form.

\begin{equation*}
\down^{k_0}(J):\;\;
\begin{aligned}
\cdots& \hphantom{=}& \times \;     & \times      & \hphantom{=}&\cdots\\
\cdots& \hphantom{=}& \times \;     & \times      & \hphantom{=}&\cdots\\
\cdots& \hphantom{=}& \hphantom{=}\;& \times      & \hphantom{=}&\cdots\\
\cdots& \hphantom{=}& \hphantom{=}\;& \times      & \hphantom{=}&\cdots\\
\cdots& \hphantom{=}& \hphantom{=}\;& \times      & \hphantom{=}&\cdots\\
\cdots& \hphantom{=}& \times      \;& \hphantom{=}& \hphantom{=}&\cdots\\
\cdots& \hphantom{=}& j_0         \;& \hphantom{=}& \hphantom{=}&\cdots
\end{aligned}
\qquad \down^{k_0}(K):\;\;
\begin{aligned}
\cdots& \hphantom{=}& \times \;     & \times      & \hphantom{=}&\cdots\\
\cdots& \hphantom{=}& \times \;     & \times      & \hphantom{=}&\cdots\\
\cdots& \hphantom{=}& \times \;     & \hphantom{=}& \hphantom{=}&\cdots\\
\cdots& \hphantom{=}& \times \;     & \hphantom{=}& \hphantom{=}&\cdots\\
\cdots& \hphantom{=}& \times \;     & \hphantom{=}& \hphantom{=}&\cdots\\
\cdots& \hphantom{=}& \times \;     & \hphantom{=}& \hphantom{=}&\cdots\\
\cdots& \hphantom{=}& j_0         \;& \hphantom{=}& \hphantom{=}&\cdots
\end{aligned}
\end{equation*}

Note that $\down^{k_0}(K)=\base(\lambda)$. In particular, we have
$s_i\base(\lambda)>\base(\lambda)$. By computing
$\down^k(J)$, for $k>k_0$, we conclude that
$\base(\f_i^{\varphi_i(\lambda)-1}\lambda)=\base(\lambda)$.

(b) Since $s_i\mu\le\mu$, $\lambda$ has the unique addable $i$-node,
which is the addable $i$-node on the first row.
Thus, $\varphi_i(\lambda)=1$ and we compare $\base(\f_i\lambda)$
and $\base(\lambda)$.
Let $J$ and $K$ be the corresponding sets of beta numbers, respectively.
Then the abacus displays of $J$ and $K$ have the following form, where
$j_0'=j_0-1$.

\begin{equation*}
J:\;\;
\begin{aligned}
\cdots& \hphantom{=}& \times \;     & \times      & \hphantom{=}&\cdots\\
\cdots& \hphantom{=}& \times \;     & \times      & \hphantom{=}&\cdots\\
\cdots& \hphantom{=}& \hphantom{=}\;& \times      & \hphantom{=}&\cdots\\
\cdots& \hphantom{=}& \hphantom{=}\;& \times      & \hphantom{=}&\cdots\\
\cdots& \hphantom{=}& \hphantom{=}\;& \times      & \hphantom{=}&\cdots\\
\cdots& \hphantom{=}& \hphantom{=}\;& \hphantom{=}& \hphantom{=}&\cdots\\
\cdots& \hphantom{=}& \hphantom{=}\;& j_0         & \hphantom{=}&\cdots
\end{aligned}
\qquad K:\;\;
\begin{aligned}
\cdots& \hphantom{=}& \times \;     & \times      & \hphantom{=}&\cdots\\
\cdots& \hphantom{=}& \times \;     & \times      & \hphantom{=}&\cdots\\
\cdots& \hphantom{=}& \hphantom{=}\;& \times      & \hphantom{=}&\cdots\\
\cdots& \hphantom{=}& \hphantom{=}\;& \times      & \hphantom{=}&\cdots\\
\cdots& \hphantom{=}& \hphantom{=}\;& \times      & \hphantom{=}&\cdots\\
\cdots& \hphantom{=}& \hphantom{=}\;& \hphantom{=}& \hphantom{=}&\cdots\\
\cdots& \hphantom{=}& j_0'        \;& \hphantom{=}& \hphantom{=}&\cdots
\end{aligned}
\end{equation*}

By a similar argument as above, there exists $k_0$ such that
$\down^{k_0}(J)$ and $\down^{k_0}(K)$ have the following form.

\begin{equation*}
\down^{k_0}(J):\;\;
\begin{aligned}
\cdots& \hphantom{=}& \times \;     & \times      & \hphantom{=}&\cdots\\
\cdots& \hphantom{=}& \times \;     & \times      & \hphantom{=}&\cdots\\
\cdots& \hphantom{=}& \hphantom{=}\;& \times      & \hphantom{=}&\cdots\\
\cdots& \hphantom{=}& \hphantom{=}\;& \times      & \hphantom{=}&\cdots\\
\cdots& \hphantom{=}& \hphantom{=}\;& \times      & \hphantom{=}&\cdots\\
\cdots& \hphantom{=}& \hphantom{=}\;& \times      & \hphantom{=}&\cdots\\
\cdots& \hphantom{=}& \hphantom{=}\;& j_0         & \hphantom{=}&\cdots
\end{aligned}
\qquad \down^{k_0}(K):\;\;
\begin{aligned}
\cdots& \hphantom{=}& \times \;     & \times      & \hphantom{=}&\cdots\\
\cdots& \hphantom{=}& \times \;     & \times      & \hphantom{=}&\cdots\\
\cdots& \hphantom{=}& \hphantom{=}\;& \times      & \hphantom{=}&\cdots\\
\cdots& \hphantom{=}& \hphantom{=}\;& \times      & \hphantom{=}&\cdots\\
\cdots& \hphantom{=}& \hphantom{=}\;& \times      & \hphantom{=}&\cdots\\
\cdots& \hphantom{=}& \times      \;& \hphantom{=}& \hphantom{=}&\cdots\\
\cdots& \hphantom{=}& j_0'        \;& \hphantom{=}& \hphantom{=}&\cdots
\end{aligned}
\end{equation*}

Thus, $\base(J)=\down^{k_0}(J)$, and by computing $\down^k(K)$,
for $k>k_0$, we have $\base(J)=s_i\base(K)>\base(K)$.
\end{proof}

\begin{lemma}
\label{second basic lemma}
Let $\lambda\in B(\Lambda_m)$ and $J$ the corresponding set of
beta numbers of charge $m$. Suppose that $\f_i\lambda\ne0$
and that $\f_iJ$ is obtained from $J$ by moving $x$ to $x+1$.
\begin{itemize}
\item[(1)]
$s_i\base(J_{\le x-1})\le\base(J_{\le x-1})$.
\item[(2)]
$\base((\f_iJ)_{\le x+1})=s_i\base(J_{\le x+1})>\base(J_{\le x+1})$.
\item[(3)]
Suppose that $\{z\in J_{\ge x+1} \mid z+e\mathbb Z=i\}\ne\emptyset$.
We denote
$$
y=\min\{z\in J_{\ge x+1} \mid z+e\mathbb Z=i\}.
$$
Then we have either
\begin{itemize}
\item[(i)]
$\base((\f_iJ)_{\le y-1})=s_i\base(J_{\le y-1})>\base(J_{\le y-1})$, or
\item[(ii)]
$\base((\f_iJ)_{\le y-1})=\base(J_{\le y-1})$.
\end{itemize}
\end{itemize}
\end{lemma}
\begin{proof}
(1) Since $\f_iJ$ is obtained from $J$ by moving $x$ to $x+1$,
$x$ is the smallest addable $i$-integer which corresponds to
a normal $i$-node. Note that all the elements in
$$
\{x-ke\in J \mid k\in\mathbb Z_{\ge1},\;x-ke+1\not\in J,\;
x-ke>M_{i+1}(J_{\le x-1})\}
$$
correspond to addable normal $i$-nodes. Thus, it
must be empty and we have
$$
M_i(J_{\le x-1})\le M_{i+1}(J_{\le x-1}).
$$
Now Lemma \ref{basic lemma}(1) implies the result.

(2) Note that $J_{\le x-1}=(\f_iJ)_{\le x-1}$ and
$$
J_{\le x+1}=J_{\le x-1}\cup\{x\},\;\;
(\f_iJ)_{\le x+1}=(\f_iJ)_{\le x-1}\cup\{x+1\}.
$$
Thus Proposition \ref{extension lemma} implies
$$
\begin{cases}
\base(J_{\le x+1})&=\base(\base(J_{\le x-1})\cup\{x\}),\\[4pt]
\base((\f_iJ)_{\le x+1})&=\base(\base(J_{\le x-1})\cup\{x+1\}).
\end{cases}
$$
As $\base(J_{\le x-1})$ is the set of beta numbers of an $e$-core,
say $\mu$, and $s_i\mu\le\mu$ by (1), and the partition associated
with $J_{\le x+1}$ has an addable $i$-node on the first row,
we are in the situation
of Lemma \ref{basic lemma}(4)(b). Note that the addable $i$-node
on the first row is the lowest addable normal $i$-node. Thus,
$(\f_iJ)_{\le x+1}=\f_i\,J_{\le x+1}$ and
$$
\base((\f_iJ)_{\le x+1})=s_i\base(J_{\le x+1})>\base(J_{\le x+1}).
$$

(3) Denote $(\f_iJ)_{\le y-1}\cap\mathbb Z_{\ge x+2}=
J_{\le y-1}\cap\mathbb Z_{\ge x+2}$ by $L$.
$L$ does not contain beads on the $i^{th}$ and the $(i+1)^{th}$
runners. The former follows from the definition of $y$.
To see the latter, observe that
there is no bead between $x$ and $y$
on the $i^{th}$ runner. Thus, if there was a bead between $x+1$ and $y-e+1$ on
the $(i+1)^{th}$ runner,
then $RA$-deletion would occur between $x$ and the bead, contradicting
the fact that $x$ corresponds to a normal $i$-node. Hence the claim follows.

Write $L=\{j_s,\dots,j_{s+r}\}$ and set
$J'_{s+r+1}=J_{\le x+1}$, $J''_{s+r+1}=(\f_i J)_{\le x+1}$.
Define $J'_k$ and $J''_k$, for $k=s+r,\dots,s$, by
$$
J'_k=\base(J'_{k+1})\cup\{j_k\}\;\;\text{and}\;\;
J''_k=\base(J''_{k+1})\cup\{j_k\}.
$$
We have
\begin{itemize}
\item[($\sharp$)]
$\base(J''_{s+r+1})=s_i\base(J'_{s+r+1})>\base(J'_{s+r+1})$ by (2).
\item[($\sharp$)]
$\base(J_{\le y-1})=\base(J'_s)$ and
$\base((\f_iJ)_{\le y-1})=\base(J''_s)$ by
Corollary \ref{cor to ext lemma}.
\end{itemize}

Suppose that $\base(J'_k)=\base(J''_k)$, for some $k$. Then we have
$$
\base((\f_iJ)_{\le y-1})=\base(J''_s)=\base(J'_s)=\base(J_{\le y-1}).
$$
Next suppose that $\base(J'_k)\ne\base(J''_k)$, for all $k$. We prove
by downward induction on $k$ that $\base(J''_k)=s_i\base(J'_k)>\base(J'_k)$.
If $k=s+r+1$ then there is nothing to prove. Suppose that the assertion holds
for $k+1$. Let
$$
J'_{k+1,t}=\down^t(\base(J'_{k+1})\cup\{j_k\})\;\;\text{and}\;\;
J''_{k+1,t}=\down^t(\base(J''_{k+1})\cup\{j_k\}),
$$
for $t\ge0$. We show that
$$
(i)\;\; M_i(J'_{k+1,t})>M_{i+1}(J'_{k+1,t}).\quad
(ii)\;\;s_iJ'_{k+1,t}=J''_{k+1,t}.
$$
When $t=0$ (i) and (ii) follow from
$\base(J''_{k+1})=s_i\base(J'_{k+1})>\base(J'_{k+1})$.

Suppose (i) and (ii) for $t$ and apply the down operation to
$J'_{k+1,t}$ and $J''_{k+1,t}$. Then,
$p'$ is the same for both and
it lies on the same runner as $j_k$. Consider $q'$ for
$J'_{k+1,t}$. Then we have one of the following.
\begin{itemize}
\item[(a)]
If $q'$ is not on the $i^{th}$ or the $(i+1)^{th}$ runners,
then
$J'_{k+1,t+1}$ and $J''_{k+1,t+1}$ are obtained by
moving $q'$ to $p'-e$ respectively.
\item[(b)]
If $q'$ is on the $i^{th}$ runner, then $J'_{k+1,t+1}$
is obtained by moving $q'$ to $p'-e$ and
$J''_{k+1,t+1}$ is obtained by moving $q'+1$ to $p'-e$.
\item[(c)]
If $q'$ is on the $(i+1)^{th}$ runner, then $J'_{k+1,t+1}$
is obtained by moving $q'$ to $p'-e$ and
$J''_{k+1,t+1}$ is obtained by moving $q'-1$ to $p'-e$.
\end{itemize}
In all the cases, we have (ii) for $t+1$. Now suppose that
(i) breaks down at $t+1$. Then
we have
$$
M_i(J'_{k+1,t})>M_{i+1}(J'_{k+1,t})\;\;\text{and}\;\;
M_i(J'_{k+1,t+1})\le M_{i+1}(J'_{k+1,t+1}).
$$
The equality does not hold in the latter, since
they are on different runners. Thus, we have
$M_i(J'_{k+1,t+1})=M_i(J'_{k+1,t})-e$ and
$M_{i+1}(J'_{k+1,t+1})=M_{i+1}(J'_{k+1,t})$, and
$$
M_i(J'_{k+1,t})-e< M_{i+1}(J'_{k+1,t+1})\le M_i(J'_{k+1,t})-e+1
$$
implies that $M_{i+1}(J'_{k+1,t+1})=M_i(J'_{k+1,t+1})+1$. Hence we
conclude that $J''_{k+1,t+1}=s_iJ'_{k+1,t+1}=J'_{k+1,t+1}$.
However, this implies $\base(J'_k)=\base(J''_k)$, contradicting our
assumption. Hence, (i) holds for $t+1$.

Therefore, $\base(J''_k)=s_i\base(J'_k)>\base(J'_k)$ holds.
By setting $k=s$ and using $\base(J_{\le y-1})=\base(J'_s)$ and
$\base((\f_iJ)_{\le y-1})=\base(J''_s)$, we have proved
$$
\base((\f_iJ)_{\le y-1})=s_i\base(J_{\le y-1})>\base(J_{\le y-1})
$$
in this case.
\end{proof}

\begin{lemma}
\label{almost base lemma}
Let $\lambda\in B(\Lambda_m)$ and $J$ the corresponding set of beta
numbers of charge $m$. Suppose that $\f_i\lambda\ne0$ and
$\f_iJ$ is obtained from $J$ by moving $x\in J$ to
$x+1\in \f_iJ$.
\begin{itemize}
\item[(1)]
Suppose that $\{z\in J_{\ge x+1} \mid z+e\mathbb Z=i\}=\emptyset$.
\begin{itemize}
\item[(a)]
If $s_i\base(\lambda)>\base(\lambda)$ then
$\base(\f_i\lambda)=s_i\base(\lambda)>\base(\lambda)$.
\item[(b)]
If $s_i\base(\lambda)\le\base(\lambda)$ then
$\base(\f_i\lambda)=\base(\lambda)$.
\end{itemize}
\item[(2)]
If $\{z\in J_{\ge x+1} \mid z+e\mathbb Z=i\}\ne\emptyset$
then $\base(\f_i\lambda)=\base(\lambda)$.
\end{itemize}
\end{lemma}
\begin{proof}
(1) Write $J_{\ge x+2}=\{j_0,\dots,j_r\}$. Set
$J'_{r+1}=J_{\le x+1}$ and $J''_{r+1}=(\f_iJ)_{\le x+1}$.
Then define $J'_k$ and $J''_k$, for $k=r,\dots,0$, by
$$
J'_k=\base(J'_{k+1})\cup\{j_k\}\;\;\text{and}\;\;
J''_k=\base(J''_{k+1})\cup\{j_k\}.
$$
Then Corollary \ref{cor to ext lemma} implies
$$
\base(J)=\base(J'_0)\;\;\text{and}\;\;
\base(\f_iJ)=\base(J''_0).
$$
There is no element of $J_{\ge x+2}$ on the $i^{th}$ runner
because $M_i(J)=x$. Suppose that there is an element of
$J_{\ge x+2}$ on the $(i+1)^{th}$ runner. We denote by
$y$ the minimal such. Then $J$ has the following layers.
\begin{equation*}
\begin{aligned}
\cdots& \hphantom{=}& x           \;& \hphantom{=}& \hphantom{=}&\cdots\\
\cdots& \hphantom{=}& \hphantom{=}\;& \hphantom{=}& \hphantom{=}&\cdots\\
\cdots& \hphantom{=}& \hphantom{=}\;& y           & \hphantom{=}&\cdots
\end{aligned}
\end{equation*}
This implies that $RA$-deletion occurs between $x$ and $y$, which
is a contradiction. Thus,
there is  also no element of $J_{\ge x+2}$ on the $(i+1)^{th}$ runner.

By Lemma \ref{second basic lemma}(2), we have
$$
\base(J''_{r+1})=s_i\base(J'_{r+1})>\base(J'_{r+1}).
$$
Hence, $J''_r=s_iJ'_r$ are $s_i$-cores and
$M_i(J'_r)>M_{i+1}(J'_r)$.

We prove by downward induction on $k$ that
\begin{itemize}
\item[(a)]
If $s_i\base(J'_k)>\base(J'_k)$ then
$\base(J''_k)=s_i\base(J'_k)>\base(J'_k)$.
\item[(b)]
If $s_i\base(J'_k)\le\base(J'_k)$ then
$\base(J''_k)=\base(J'_k)$.
\end{itemize}
When $k=r$, (a) and (b) follow from Lemma \ref{basic lemma}(3).
Suppose that (a) and (b) hold for $k+1$.
Then we have either
\begin{itemize}
\item[(a')]
$J''_k=s_iJ'_k$ are $s_i$-cores and
$M_i(J'_k)>M_{i+1}(J'_k)$, or
\item[(b')]
$J''_k=J'_k$ is an $s_i$-core and $M_i(J'_k)\le M_{i+1}(J'_k)$.
\end{itemize}

Suppose that $s_i\base(J'_k)>\base(J'_k)$. Then (b') does
not occur by Lemma \ref{basic lemma}(1).
Thus (a') must occur and Lemma \ref{basic lemma}(3) implies
$$
\base(J''_k)=s_i\base(J'_k)>\base(J'_k).
$$
Suppose that $s_i\base(J'_k)\le\base(J'_k)$. If (b') occurs
then $\base(J''_k)=\base(J'_k)$ obviously holds, so we may
assume that (a') occurs. Then, Lemma \ref{basic lemma}(3)
implies $\base(J''_k)=\base(J'_k)$ also.
We have proved that (a) and (b) hold for $k$.

Setting $k=0$ and using $\base(J)=\base(J'_0)$ and
$\base(\f_iJ)=\base(J''_0)$, we have the desired result.

(2) Define $y=\min\{z\in J_{\ge x+1} \mid z+e\mathbb Z=i\}$ as before.
Then, by Proposition \ref{extension lemma},
$$
\begin{cases}
\base(J_{\le y})&=\base(\base(J_{\le y-1})\cup\{y\}),\\
\base((\f_iJ)_{\le y})&=
\base(\base((\f_iJ)_{\le y-1})\cup\{y\}).
\end{cases}
$$
Let $J'=\base(J_{\le y-1})\cup\{y\}$ and
$J''=\base((\f_iJ)_{\le y-1})\cup\{y\}$.
By Lemma \ref{second basic lemma}(3) we have either
\begin{itemize}
\item[(i)]
$\base((\f_iJ)_{\le y-1})=s_i\base(J_{\le y-1})>\base(J_{\le y-1})$, or
\item[(ii)]
$\base((\f_iJ)_{\le y-1})=\base(J_{\le y-1})$.
\end{itemize}
Let $\lambda'$ be the partition whose set of beta numbers of charge
$m$ is $J'$.
If (i) occurs then $J''=\f_i^{\varphi_i(\lambda')-1}J'$ and
we are in the situation
of Lemma \ref{basic lemma}(4)(a). Thus we have
$\base(J')=\base(J'')$. If (ii) occurs then $J'=J''$ and
we have $\base(J')=\base(J'')$ again. Thus,
$\base((\f_i J)_{\le y})=\base(J_{\le y})$ in both cases. Now
Corollary \ref{cor to ext lemma} implies
$\base(\f_iJ)=\base(J)$.
\end{proof}

The next theorem is the counterpart to
Theorem \ref{roof lemma}, the \lq\lq roof lemma\rq\rq
of \cite{KLMW}.

\begin{thm}
Let $\lambda\in B(\Lambda_m)$. Then
$$
\base(\f_i^{max}\lambda)=\begin{cases}
s_i\base(\lambda)\quad
(\text{if $\base(\lambda)$ has an addable $i$-node})\\
\base(\lambda) \quad(\text{otherwise})\end{cases}
$$
and $\base(\f_i^t\lambda)=\base(\lambda)$,
for $0\le t<\varphi_i(\lambda)$.
\end{thm}
\begin{proof}
The theorem is equivalent to the following two statements.
\begin{itemize}
\item[(a)]
If $\varphi_i(\lambda)=1$ and $s_i\base(\lambda)>\base(\lambda)$
then
$$
\base(\f_i\lambda)=s_i\base(\lambda)>\base(\lambda).
$$
\item[(b)]
Otherwise $\base(\f_i\lambda)=\base(\lambda)$.
\end{itemize}

Suppose that the assumption in (a) holds.
Then $s_i\base(\lambda)>\base(\lambda)$
implies $M_i(\lambda)>M_{i+1}(\lambda)$ by Lemma \ref{basic lemma}.
Thus $M_i(\lambda)$ corresponds to an addable normal $i$-node.
As $\varphi_i(\lambda)=1$, $\f_i\lambda$ is obtained from
$\lambda$ by adding this node. We apply Lemma \ref{almost base lemma}.
Then $x=M_i(\lambda)$ and (1)(a) applies. Hence the result follows.

Suppose that the assumption in (b) holds. Then we have
$s_i\base(\lambda)\le\base(\lambda)$ or $\varphi_i(\lambda)\ge2$.
In the former case, either (1)(b) or (2) of Lemma \ref{almost base lemma}
applies. In the latter case, Lemma \ref{almost base lemma}(2) applies.
Hence $\base(\f_i\lambda)=\base(\lambda)$ follows in both cases.
\end{proof}

\begin{cor}
\label{base is floor}
$\base(\lambda)=\floor(\lambda)$.
\end{cor}
\begin{proof}
Note that Lemma \ref{rule to compute paths} implies that
$$
\floor(\f_i^{max}\lambda)=\begin{cases}
s_i\floor(\lambda)\quad
(\text{if $\floor(\lambda)$ has an addable $i$-node})\\
\floor(\lambda) \quad(\text{otherwise})\end{cases}
$$
and $\floor(\f_i^t\lambda)=\floor(\lambda)$,
for $0\le t<\varphi_i(\lambda)$.
Thus induction on the size of
$\lambda$ proves the result.
\end{proof}

The next theorem follows from Theorem \ref{floor theorem} and
Corollary \ref{base is floor}.

\begin{thm}
\label{AKT}
In the partition realization of $B(\Lambda_m)$, we have
$$
B^w(\Lambda_m)=\{\lambda\in B(\Lambda_m) \mid
\base(\lambda)\supset w\emptyset_m\}.
$$
\end{thm}

Recall that $w_m$ is the longest element of $W_m$.

\begin{cor}
\label{base lemma}
Write $\base(\lambda)=w_\lambda\emptyset_m$, for a unique
$w_\lambda\in W/W_m$. Then
$$
w_\lambda w_m=\max\,\{w\in W \mid \lambda\in B^w(\Lambda_m)\}
$$
with respect to the Bruhat-Chevalley order.
\end{cor}

\section{Kleshchev multipartitions}

Recall that $M_i(\lambda)=\max\{x\in J \mid x+e\mathbb Z=i\}$.

\begin{defn}
Let $\lambda\in B(\Lambda_0)$ be an $e$-core, $J$ the
corresponding set of beta numbers of charge $0$.
Write $\{M_i(\lambda)\}_{i\in\mathbb Z/e\mathbb Z}$
in descending order
$$
M_{i_1}(\lambda)>M_{i_2}(\lambda)>\cdots>M_{i_e}(\lambda).
$$
Then define
$\tau_m(J)=J\cup\{ M_{i_k}(\lambda)+e\}_{1\le k\le m}$,
and denote the corresponding $e$-restricted partition
by $\tau_m(\lambda)\in B(\Lambda_m)$. If $m=0$ then
$\tau_m(\lambda)=\lambda$.
\end{defn}

Recall from the definition
of \r{W} in \cite[p.74]{Kc} and \cite[Proposition 6.5]{Kc} that
$W$ is the semidirect product of $W_0$ and $T$, where
$T=\{t_\alpha \mid \alpha\in\oplus_{i=1}^{e-1}\mathbb Z\alpha_i\}$, and
$T$ acts on weights by
$$
t_\alpha\Lambda=\Lambda+\Lambda(c)\alpha-
((\Lambda,\alpha)+\frac{1}{2}|\alpha|^2\Lambda(c))\delta.
$$
See \cite[(6.5.2)]{Kc}. Thus, any weight in the
$W$-orbit $W\Lambda_0$ is of the form
$t_\alpha\Lambda_0$, for some $t_\alpha\in T$.
Note that $t_\alpha$ is not necessarily a distinguished
coset representative.

\begin{lemma}
\label{content sequence}
Suppose that $\lambda\in B(\Lambda_0)$ is an $e$-core, and
write $\lambda=t_\alpha\emptyset_0$, for
$\alpha=\sum_{i=1}^{e-1}m_i\alpha_i$. Then
$m_i=N_0(\lambda)-N_i(\lambda)$, for $1\le i\le e-1$.
\end{lemma}
\begin{proof}
As $wt(\lambda)=t_\alpha\Lambda_0=\Lambda_0+\alpha-\frac{1}{2}|\alpha|^2\delta$,$$
\sum_{i=0}^{e-1}N_i(\lambda)\alpha_i=
\Lambda_0-t_\alpha\Lambda_0=\frac{1}{2}|\alpha|^2\delta-\alpha.
$$
Thus $N_0(\lambda)=\frac{1}{2}|\alpha|^2$ and
$N_i(\lambda)=\frac{1}{2}|\alpha|^2-m_i$, for
$1\le i\le e-1$. The result follows.
\end{proof}

\begin{prop}
\label{translation}
Let $\lambda=w\emptyset_0\in B(\Lambda_0)$ and let
$\mu=w'\emptyset_m\in B(\Lambda_m)$, where
$w\in W/W_0$ and $w'\in W/W_m$. Then $ww_0\ge w'$ if
and only if $\tau_m(\lambda)\supset\mu$.
\end{prop}
\begin{proof}
We use Proposition \ref{two orders} throughout freely, without
comment.

We may write $\lambda=t_\alpha\emptyset_0$, for
$\alpha=\sum_{i=1}^{e-1}m_i\alpha_i$, and
$t_\alpha=wv$, for $v\in W_0$. Then $ww_0\emptyset_m=t_\alpha u\emptyset_m$ for
$u=v^{-1}w_0\in W_0$.
On the other hand, if $u\in W_0$ then $t_\alpha u\le ww_0$, which implies
$t_\alpha u\emptyset_m\subset ww_0\emptyset_m$. Thus
$$
ww_0\emptyset_m=\max\{t_\alpha u\emptyset_m \mid u\in W_0\}.
$$
If $ww_0\ge w'$ then $ww_0\emptyset_m\supset w'\emptyset_m=\mu$,
and conversely, if $ww_0\emptyset_m\supset\mu$ then
$ww_0\ge w'$. Thus we want to show $ww_0\emptyset_m=\tau_m(\lambda)$.

Suppose that $m=0$. Then $ww_0\emptyset_m=w\emptyset_m=\lambda$ and
$ww_0\emptyset_m=\tau_m(\lambda)$ is trivial.

Suppose that $m\ne0$. Fix $u\in W_0$ and write
$u\Lambda_m=\Lambda_m-\beta$, for some
$\beta\in\sum_{i=1}^{e-1}\mathbb Z_{\ge0}\alpha_i$.
Then
$$
\begin{cases}
t_\alpha\Lambda_m=\Lambda_m+\alpha-
((\Lambda_m,\alpha)+\frac{1}{2}|\alpha|^2)\delta,\\
t_\alpha\beta=\beta-(\beta,\alpha)\delta.
\end{cases}
$$
We also have
$t_\alpha\Lambda_0=\Lambda_0+\alpha-\frac{1}{2}|\alpha|^2\delta$,
which implies
$\sum_{i=0}^{e-1}N_i(\lambda)\alpha_i
=\frac{1}{2}|\alpha|^2\delta-\alpha$ as before.
Therefore,
\begin{equation*}
\begin{split}
\sum_{i=0}^{e-1}N_i(t_\alpha u\emptyset_m)\alpha_i&=
\Lambda_m-t_\alpha u\Lambda_m=
\Lambda_m-t_\alpha(\Lambda_m-\beta)\\
&=\bigl((\Lambda_m,\alpha)+\frac{1}{2}|\alpha|^2-(\beta,\alpha)\bigr)\delta
-\alpha+\beta\\
&=(\Lambda_m-\beta,\alpha)\delta+\beta+\sum_{i=0}^{e-1}N_i(\lambda)\alpha_i.
\end{split}
\end{equation*}

As $t_\alpha u\emptyset_m\subset ww_0\emptyset_m$, for all $u$,
the height of $(\Lambda_m-\beta,\alpha)\delta+\beta$ must
attain a maximum value at $ww_0\emptyset_m$.

As $u\in W_0$, we may compute $u\Lambda_m$ by restricting the weights
to $\mathfrak g(A_{e-1})$. Hence we consider the restricted weights
for the moment, and, by abuse of notation, we use the same $u\Lambda_m$.
Then, $\Lambda_m$ may be considered as
the weight $\epsilon_1+\cdots+\epsilon_m$ of
$\mathfrak g(A_{e-1})=sl(e,\mathbb C)$, where the weight lattice of
$sl(e,\mathbb C)$ is realized as $\oplus_{i=1}^{e-1}
\mathbb Z\epsilon_i$ with $\sum_{i=1}^{e-1}\epsilon_i=0$ as usual,
and the simple roots are $\{\alpha_i=\epsilon_i-\epsilon_{i+1}\}_{1\le i<e}$.
Thus,
$$
u\Lambda_m=\Lambda_m-\beta\in\{\epsilon_{i_1}+\cdots+\epsilon_{i_m} \mid
1\le i_1<\cdots<i_m\le e\}.
$$
Write $u\Lambda_m=\sum_{k=1}^m\epsilon_{i_k}$. Note that
$(\epsilon_i,\epsilon_j)=\delta_{ij}$ and
we may compute $(\Lambda_m-\beta,\alpha)$ by using the restricted weights.
Thus, by Lemma \ref{content sequence},
$$
(\Lambda_m-\beta,\alpha)=\sum_{k=1}^m (\epsilon_{i_k},\alpha)=
\sum_{k=1}^m (m_{i_k}-m_{i_k-1})=
\sum_{k=1}^m(N_{i_k-1}(\lambda)-N_{i_k}(\lambda)).
$$
As $\beta=\sum_{k=1}^m (\epsilon_k-\epsilon_{i_k})$, the height of
$\beta$ is $\sum_{k=1}^m (i_k-k)$. Therefore, the value to be
maximized is
$$
\sum_{k=1}^m (N_{i_k-1}(\lambda)-N_{i_k}(\lambda))e+(i_k-k).
$$
Define $L_i=(N_{i-1}(\lambda)-N_i(\lambda))e+i$, for
$1\le i\le e$. Here, we understand that $N_e(\lambda)=N_0(\lambda)$.
It is important that the range for $i$ is
not $0\le i\le e-1$ but $1\le i\le e$.
Let $J$ be the set of beta numbers of charge $0$
associated with $\lambda$ and
$M_i(\lambda)=\max\{x\in J \mid x+e\mathbb Z=i\}$ as before.
Then,
$$
\sum_{i=0}^{e-1}N_i(t_\alpha u\emptyset_m)\alpha_i
=\bigl(\sum_{k=1}^m\frac{L_{i_k}-i_k}{e}\bigr)\delta+
\sum_{k=1}^m(\alpha_k+\cdots+\alpha_{i_k-1})
+\sum_{i=0}^{e-1}N_i(\lambda)\alpha_i.
$$

We claim that $L_i=M_i(\lambda)+e$, for $1\le i\le e$.
Recall how to read
$N_i(\lambda)$ from the abacus. We explain this by an example.
Let $\lambda=(4,2)$ and $e=6$. Then the corresponding
$J$ is displayed as follows.

\medskip
\begin{center}
\begin{tabular}{rrrrrr}
 $-12$ & $-11$ & $-10$ & $-9$ & $-8$ & $-7$ \\
 $-6$  & $-5$  & $-4$  & $-3$ & $-2$ &      \\
       & $1$   &       &      & $4$  &      \\
       &       &       &      &      &
\end{tabular}
\end{center}

We read the numbers on the abacus from $-\infty$ and
with initial value $0$, and increment the value by $1$
at each number which does not belong to $J$. Equivalently,
the value at $x$ is $|\{y\le x \mid y\not\in J\}|$. We obtain

\medskip
\begin{center}
\begin{tabular}{cccccc}
 $0$ &  $0$ &  $0$ & $0$ & $0$ & $0$ \\
 $0$ &  $0$ &  $0$ & $0$ & $0$ & $1$ \\
 $2$ &  $2$ &  $3$ & $4$ & $4$ & $5$ \\
 $6$ &  $7$ &  $8$ & $9$ & $10$& $11$
\end{tabular}
\end{center}

\medskip
We consider the same for the empty partition. Then we have

\medskip
\begin{center}
\begin{tabular}{cccccc}
 $0$ & $0$ & $0$ & $0$ & $0$ & $0$ \\
 $0$ & $0$ & $0$ & $0$ & $0$ & $0$ \\
 $0$ & $1$ & $2$ & $3$ & $4$ & $5$ \\
 $6$ & $7$ & $8$ & $9$ & $10$& $11$
\end{tabular}
\end{center}

\medskip
We compute the difference and obtain:

\medskip
\begin{center}
\begin{tabular}{cccccc}
 $0$ & $0$ & $0$ & $0$ & $0$ & $0$ \\
 $0$ & $0$ & $0$ & $0$ & $0$ & $1$ \\
 $2$ & $1$ & $1$ & $1$ & $0$ & $0$ \\
 $0$ & $0$ & $0$ & $0$ & $0$ & $0$
\end{tabular}
\end{center}

\medskip
Then $N_i(\lambda)$ is the summation of the
entries on the $i^{th}$ runner.
$$
N_0(\lambda)=2,\;N_1(\lambda)=1,\;N_2(\lambda)=1,\;
N_3(\lambda)=1,\;N_4(\lambda)=0,\;N_5(\lambda)=1.
$$
In this example, we have
$$
L_1=7,\;L_2=2,\;L_3=3,\;L_4=10,\;L_5=-1,\;L_6=0.
$$
The proof of this rule is by induction on the size of $\lambda$.
If $x\in J$ moves to $x+1$ when adding a node, then, as is
explained in Example \ref{J records contents}, the box to be
added has the content $x$. Then observe that
$|\{y\le x \mid y\not\in J\}|$ increases by $1$ at $x$.

Let $\alpha$ and $\beta=\alpha+1$ be two consecutive numbers
such that $\alpha\in i-1+e\Z$ and $\beta\in i+e\Z$.

Suppose that
$\alpha\ge0$. Then, by the above rule for computing
$N_i(\lambda)$, we have
\begin{itemize}
\item[(a)]
If $\beta\in J$ then the values at $\alpha$ and $\beta$
are the same. Thus, they contribute $1$
to $N_{i-1}(\lambda)-N_i(\lambda)$.
\item[(b)]
If $\beta\not\in J$ then the value at $\beta$ is
greater than the value at $\alpha$ by $1$. Thus,
they do not contribute
to $N_{i-1}(\lambda)-N_i(\lambda)$.
\end{itemize}
Similarly, if $\alpha<0$, then we have.
\begin{itemize}
\item[(a)]
If $\beta\not\in J$ then they contribute $-1$
to $N_{i-1}(\lambda)-N_i(\lambda)$.
\item[(b)]
If $\beta\in J$ then they do not contribute
to $N_{i-1}(\lambda)-N_i(\lambda)$.
\end{itemize}

Suppose that $M_i(\lambda)\ge1$.
We have, for example,
\begin{equation*}
\begin{aligned}
\cdot &\hphantom{=}&\hphantom{=}& \times      \;& \times  \\
\cdot &\hphantom{=}&\hphantom{=}& \hphantom{=}\;& \times  \\
0     &\hphantom{=}&\hphantom{=}& \hphantom{=}\;& \times  \\
\cdot &\hphantom{=}&\hphantom{=}& \hphantom{=}\;& \times  \\
\cdot &\hphantom{=}&\hphantom{=}& \hphantom{=}\;& \times  \\
\end{aligned}
\end{equation*}
Then only those $\beta\in J$ with $\alpha\ge0$ contribute and
the number of such is $\frac{M_i(\lambda)+e-i}{e}$. Hence
$L_i=\frac{M_i(\lambda)+e-i}{e}e+i=M_i(\lambda)+e$. Next
suppose that $M_i(\lambda)\le0$.

\begin{equation*}
\begin{aligned}
\cdot &\hphantom{=}&\hphantom{=}& \times  \;& \times  \\
\cdot &\hphantom{=}&\hphantom{=}& \times  \;&\hphantom{=}\\
0     &\hphantom{=}&\hphantom{=}& \times  \;&\hphantom{=}\\
\cdot &\hphantom{=}&\hphantom{=}& \times  \;&\hphantom{=}\\
\cdot &\hphantom{=}&\hphantom{=}& \times  \;&\hphantom{=}\\
\end{aligned}
\end{equation*}
Then only those $\beta\not\in J$ with $\alpha<0$ contribute and
the number of such is $\frac{i-e-M_i(\lambda)}{e}$. Hence
$L_i=-\frac{i-e-M_i(\lambda)}{e}e+i=M_i(\lambda)+e$.
We have proved $L_i=M_i(\lambda)+e$.

Recall that we want to maximize $\sum_{k=1}^m L_{i_k}$. This is
achieved precisely when $\{L_{i_k}-e \mid 1\le k\le m\}$ consists
of the largest $m$ numbers of $\{M_i(\lambda) \mid 1\le i\le e\}$.
From now on, we suppose that
$$
\{ M_{i_1}(\lambda), M_{i_2}(\lambda),\dots, M_{i_m}(\lambda)
\mid 1\le i_1<\cdots<i_m\le e \}
$$
are the largest $m$ numbers of $\{M_i(\lambda) \mid 1\le i\le e\}$.
We write $M_{i_k}$ for $M_{i_k}(\lambda)$. Then
$$
\sum_{i=0}^{e-1}N_i(ww_0\emptyset_m)\alpha_i=
\bigl(\sum_{k=1}^m\frac{M_{i_k}+e-i_k}{e}\bigr)\delta+
\sum_{k=1}^m(\alpha_k+\cdots+\alpha_{i_k-1})
+\sum_{i=0}^{e-1}N_i(\lambda)\alpha_i.
$$
We compute $\Lambda_0-wt(\lambda)$ and
$\Lambda_m-wt(\tau_m(\lambda))$. For the computation,
it is helpful to
view a partition as a difference of two diagrams both of which
extend infinitely to the left.
Let $\mu\in B(\Lambda_m)$ and define two subsets of $\Z^2$ by
$$
A=\{(i,j) \mid i\ge -m,\;j< \mu_{i+m}\}\;\;\text{and}\;\;
B=\{(i,j) \mid i\ge -m,\;j< 0\},
$$
where the $i$-coordiate increases downward as in English convention.
We also define the residue of $x=(i,j)\in\Z^2$ by
$\res(x)=-i+j+e\Z\in\Z/e\Z$. Then
$$
\Lambda_m-wt(\mu)=\sum_{x\in A\setminus B}\alpha_{\res(x)}
=\sum_{x\in A}\alpha_{\res(x)}-\sum_{x\in B}\alpha_{\res(x)}.
$$
We can justify the rightmost by considering the region
$D=\{(i,j) \mid i\le N, j\ge N'\}$, for sufficiently large
$N$ and $-N'$, and understand it as
$$
\sum_{x\in A\cap D}\alpha_{\res(x)}-\sum_{x\in B\cap D}\alpha_{\res(x)}.
$$
Let $k_0>k_1>\cdots$ be the beta numbers of $\mu$. Thus,
$k_j=\mu_j+m-j$. We may read them from $\mu$ as Example \ref{J records contents}.
Then we may write
$$
\sum_{x\in A}\alpha_{\res(x)}=\sum_{j\ge0}\sum_{s<k_j}\alpha_s,\;\;\text{and}\;\;
\sum_{x\in B}\alpha_{\res(x)}=\sum_{j\ge0}\sum_{s<m-j}\alpha_s.
$$
They do not make sense, but their difference does.
Note that we can rearrange the order of a finite number of rows of $A$ or $B$
to compute $\Lambda_m-wt(\mu)$.

Now we compare $\Lambda_0-wt(\lambda)$ and $\Lambda_m-wt(\tau_m(\lambda))$.
Let $A=\{(i,j) \mid i\ge 0,\;j<\lambda_i\}$ and
$B=\{(i,j) \mid i\ge 0,\;j<0\}$. Then
$$
\Lambda_0-wt(\lambda)=\sum_{x\in A}\alpha_{\res(x)}-
\sum_{x\in B}\alpha_{\res(x)}.
$$
Define $A', B'\subset\Z^2$ by
$$
A'=\{(-k,j) \mid 1\le k\le m,\; j<M_{i_k}+e-k\},\;\;
B'=\{(-k,j) \mid 1\le k\le m,\; j<0\}.
$$
Then
$$
\Lambda_m-wt(\tau_m(\lambda))=
\sum_{x\in A\cup A'}\alpha_{\res(x)}-\sum_{x\in B\cup B'}\alpha_{\res(x)}.
$$
Thus $(\Lambda_m-wt(\tau_m(\lambda)))-(\Lambda_0-wt(\lambda))$
is given by
$$
\sum_{x\in A'}\alpha_{\res(x)}-\sum_{x\in B'}\alpha_{\res(x)}.
$$
Observe that the first term is given by
$\sum_{k=1}^m\sum_{j<M_{i_k}+e}\alpha_j$ and the second term is given by
$\sum_{k=1}^m\sum_{j<k}\alpha_j$.
Thus, for a sufficiently large $N$, we have
$$
\Lambda_m-wt(\tau_m(\lambda))=\sum_{k=1}^m
\bigl(\sum_{j=-N}^{M_{i_k}+e-1}\alpha_j
-\sum_{j=-N}^{k-1}\alpha_j\bigr)+\Lambda_0-wt(\lambda),
$$
and each term in the sum is equal to
$$
\frac{M_{i_k}+e-i_k}{e}\delta+(\alpha_k+\cdots+\alpha_{i_k-1}).
$$
Hence $\Lambda_m-wt(\tau_m(\lambda))=
\Lambda_m-wt(ww_0\emptyset_m)$, which implies
$ww_0\emptyset_m=\tau_m(\lambda)$.
\end{proof}

We may describe $\tau_m(\lambda)$, for $1\le m<e$, by
Young diagrammatic terms. To see this, let
$\ell=\ell(\lambda)$ be the length of
$\lambda=(\lambda_0,\lambda_1,\dots)$ and define
$$
\nu_i=\begin{cases} \lambda_i+e-m\quad&(0\le i<m)\\
\min\{\lambda_i+e-m,\lambda_{i-m}\}&(m\le i)\end{cases}
$$
We have $\nu_i=0$ if and only if $i\ge\ell+m$.
It is clear that $\nu_0\ge\cdots\ge\nu_{m-1}$ and
$\nu_m\ge\nu_{m+1}\ge\cdots$. As $\nu_{m-1}<\nu_m$ would imply
$\lambda_{m-1}+e-m<\lambda_m+e-m$, we have
$\nu_{m-1}\ge\nu_m$. Hence, $\nu$ is a partition.

Let $\operatorname{shift}^m(\lambda)=(0^m,\lambda_0,\dots,\lambda_{\ell-1},0,\dots)$.
We denote by $a^b$ the partition $(a^b,0,\dots)$. The sum of
partitions is defined by $\lambda+\mu=(\lambda_0+\mu_0,\lambda_1+\mu_1,\dots)$.
The following proposition shows that
$$
\tau_m(\lambda)=
(\lambda+(e-m)^{\ell+m})\cap((\lambda_1+e-m)^m+\operatorname{shift}^m(\lambda)).
$$
In particular, we have
$$
a(\tau_m(\lambda))=a(\lambda)+e-m\;\;\text{and}\;\;
\ell(\tau_m(\lambda))=\ell(\lambda)+m.
$$

\begin{prop}
Let $\lambda$ be an $e$-core, and define $\nu$ as above. Then
$$
\nu=(\nu_0,\dots,\nu_{\ell+m-1},0,\dots)=\tau_m(\lambda).
$$
\end{prop}
\begin{proof}
Let $J$ be the set of beta numbers of charge $0$ associated
with $\lambda$, and let $K$ be the set of beta numbers of charge
$m$ associated with $\nu$. Then
$$
k_{i+m}=\min\{\lambda_{i+m}+e-m-i,\lambda_i-i\}=\min\{j_{i+m}+e, j_i\},
$$
for $i\ge0$. We also have $k_i=j_i+e$, for $0\le i<m$.
Hence, to obtain $K$ from $J$, we start with $J+e$, namely
we slide down all the beads by one on the abacus, and move
$j_{i+m}+e$ to $j_i$ when $j_{i+m}+e>j_i$, for $i\ge0$.
Since $\nu$ is a partition, $j_i=j_{i'+m}+e>j_{i'}$, for some
$i'$, when it occurs.

Our aim is to prove that $K=J\cup\{M_{i_k}(\lambda)+e\}_{1\le k\le m}$.
First we show that $x\in J$ implies $x\in K$.
Suppose that $x=j_i$ and $x\not\in K$. Since $x$ must move,
$j_i=j_{i'+m}+e>j_{i'}$, for some $i'\ge0$. Thus $i<i'$ and
$j_{i+m}+e>j_{i'+m}+e=j_i$. Hence $j_{i+m}+e$ moves to $x$,
which contradicts the assumption $x\not\in K$.

Next consider $x\in\{M_i(\lambda)+e\}_{i\in\Z/e\Z}$.
As $x\not\in J$, no $j_{i'+m}+e\in J+e$ moves to $x$.
Hence $x\not\in K$ if and only if $x=j_{i+m}+e>j_i$, for
some $i$. Let $x=j_{i+m}+e$. We have to show that
$j_{i+m}+e>j_i$ if and only if
$x\not\in\{M_{i_k}(\lambda)+e\}_{1\le k\le m}$.
If $j_{i+m}+e>j_i$ then $j_i>j_{i+1}>\cdots>j_{i+m-1}\ge j_{i+m}+1$ implies
$$
\{j_{i+m-1}, j_{i+m-2},\dots, j_i\}\subset
\{j_{i+m}+1,j_{i+m}+2,\dots,j_{i+m}+e-1\}\cap J.
$$
Hence $j_{i+m-1}+e,\dots,j_i+e$ are in pairwise distinct runners
and all of them are greater than $x$. We have proved
$x\not\in\{M_{i_k}(\lambda)+e\}_{1\le k\le m}$.
If $j_{i+m}+e\le j_i$ then there exists $i+m-1\ge i'>i$ such that
$$
\{j_{i+m-1}, j_{i+m-2},\dots, j_{i'}\}=
\{j_{i+m}+1,j_{i+m}+2,\dots,j_{i+m}+e-1\}\cap J.
$$
In fact, it is clear that $j_{i+m-1}$ is the minimal element of the
right hand side. Denote the maximal element by $j_{i'}$. Then
$j_{i'}<j_{i+m}+e\le j_i$ implies $i'>i$.

These beads are in pairwise distinct runners.
Each of the $i+m-i'(<m)$ runners has a bead which is greater than $x$,
but the remaining runners do not have such a bead. Hence
$x\in\{M_{i_k}(\lambda)+e\}_{1\le k\le m}$.
\end{proof}

We are now prepared to prove the following.

\begin{thm}
\label{main result}
Let $\lambda\otimes\mu\in B(\Lambda_0)\otimes B(\Lambda_m)$.
Then $\lambda\otimes \mu\in B(\Lambda_0+\Lambda_m)$ if and
only if
$$
\tau_m(\base(\lambda))\supset\roof(\mu).
$$
\end{thm}
\begin{proof}
Suppose that $m=0$. By Corollary \ref{roof is ceil} and
Corollary \ref{base is floor}, $\base(\lambda)\supset\roof(\mu)$
is equivalent to $\floor(\lambda)\supset\ceil(\mu)$.
Write $\floor(\lambda)=w\emptyset_0$ and
$\base(\mu)=w'\emptyset_0$, for $w, w'\in W/W_0$.
Then $\floor(\lambda)\supset\ceil(\mu)$ is equivalent to $w\ge w'$,
which is further equivalent to
$$
f(\lambda)=w\Lambda_0\ge w'\Lambda_0=i(\mu).
$$
Hence Corollary \ref{the result to use} for $r=d=2$ implies
the result.

Suppose that $m\ne0$. Write $\base(\lambda)=w\emptyset_0$ and
$\roof(\mu)=w'\emptyset_m$, for $w\in W/W_0$ and
$w'\in W/W_m$ respectively. Then Corollary \ref{the result to use}
for $r=2, d=1$ implies that
$\lambda\otimes\mu\in B(\Lambda_0+\Lambda_m)$ if and only if
$ww_0\ge w'$. This is equivalent to
$\tau_m(\base(\lambda))\supset\roof(\mu)$ by
Proposition \ref{translation}.
\end{proof}

Let $\H_n$ be the cyclotomic Hecke algebra defined by
$(T_0+1)^d(T_0+q^m)^{r-d}=0$, $(T_i-q)(T_i+1)=0$, for
$1\le i<n$, and the type $B$ braid relations.
As was mentioned in the introduction, a complete set
of simple $\H_n$-modules is given by
the set of nonzero $D^{(\lambda^{(r)},\dots,\lambda^{(1)})}$'s,
where $D^{(\lambda^{(r)},\dots,\lambda^{(1)})}$ is obtained
from the Specht module $S^{(\lambda^{(r)},\dots,\lambda^{(1)})}$
by factoring out the radical of the invariant symmetric
bilinear form defined on it. The complete set is naturally
a $\mathfrak g(A^{(1)}_{e-1})$-crystal $B(\Lambda)$, where
$\Lambda=d\Lambda_0+(r-d)\Lambda_m$. See \cite{AM} and \cite{A2},
or \cite{A1}. Note that when $r=2$ and
$Q=-q^m$, we obtain the Hecke algebra $\H_n(Q,q)$ of type
$B$ as special cases. Theorem \ref{main result} combined
with the results explained in the introduction
gives the following.

\begin{cor}
\label{main result2}
Let $\underline\lambda=\lambda^{(1)}\otimes\cdots\otimes
\lambda^{(r)}\in
B(\Lambda_0)^{\otimes d}\otimes B(\Lambda_m)^{\otimes r-d}$.
Then the following are equivalent.
\begin{itemize}
\item[(i)]
$D^{(\lambda^{(r)},\dots,\lambda^{(1)})}\ne0$.
\item[(ii)]
$\underline\lambda\in B(d\Lambda_0+(r-d)\Lambda_m)$.
\item[(iii)]
The following three conditions hold.
\begin{itemize}
\item[(a)]
$\base(\lambda^{(k)})\supset\roof(\lambda^{(k+1)})$, for $1\le k<d$,
\item[(b)]
$\tau_m(\base(\lambda^{(d)}))\supset\roof(\lambda^{(d+1)})$,
\item[(c)]
$\base(\lambda^{(k)})\supset\roof(\lambda^{(k+1)})$, for $d<k<r$.
\end{itemize}
\end{itemize}
\end{cor}

Recall that $a(\lambda)$ is the length of the first row, and
$\ell(\lambda)$ is the length of the first column.
For $\underline\lambda=\lambda^{(1)}\otimes\cdots\otimes
\lambda^{(r)}$, define
$a_i(\underline\lambda)=a(\lambda^{(i)})-\ell(\lambda^{(i+1)})$.
Mathas proved the following result.

\begin{prop}
Suppose that $e=2$ and let
$$
\underline\lambda=\lambda^{(1)}\otimes\cdots\otimes
\lambda^{(r)}\in
B(\Lambda_{m_1})\otimes\cdots\otimes B(\Lambda_{m_r})=
B(\Lambda_0)^{\otimes d}\otimes B(\Lambda_m)^{\otimes r-d}.
$$
Then $\underline\lambda\in B(d\Lambda_0+(r-d)\Lambda_m)$ if and only if
$a_i(\underline\lambda)\ge \delta_{m_im_{i+1}}-1$, for $1\le i<r$.
\end{prop}

Observe that any $2$-core $\lambda$ is of the form
$(c,c-1,\dots,1)$ and $a(\lambda)=\ell(\lambda)=c$.
Using the closed formulas for $\ceil(\lambda)$ and $\floor(\lambda)$ for
a partition $\lambda$ which is given in
Proposition \ref{floor and ceil for e=2},
we have
\begin{itemize}
\item[(i)]
If $m_i=m_{i+1}$ then
$\floor(\lambda^{(i)})\supset\ceil(\lambda^{(i+1)})$
is equivalent to
$$
a(\lambda^{(i)})\ge\ell(\lambda^{(i+1)}).
$$
\item[(ii)]
If $m_i\ne m_{i+1}$ then
$\tau_1(\floor(\lambda^{(i)}))\supset\ceil(\lambda^{(i+1)})$
is equivalent to
$$
a(\lambda^{(i)})+1\ge\ell(\lambda^{(i+1)}).
$$
\end{itemize}
Thus, Mathas' result follows from our results.

Now consider $e=3$. Recently, in the spirit similar to Mathas' result in
$e=2$, Fayers has obtained a necessary and sufficient condition
for $(\lambda,\mu)$ to be a Kleshchev bipartition \cite{F}.
According to him, the condition may be restated as follows.

\begin{prop}
\label{Fayers}
Suppose that $e=3$ and let
$\lambda\otimes\mu\in B(\Lambda_0)\otimes B(\Lambda_m)$.
\begin{itemize}
\item[(i)]
If $m=0$ then $\lambda\otimes\mu\in B(\Lambda_0+\Lambda_m)$ if and
only if
$$
a(\lambda)\ge \ell(m(\mu))\;\;\text{and}\;\;a(m(\lambda))\ge \ell(\mu).
$$
\item[(ii)]
If $m=1$ then $\lambda\otimes\mu\in B(\Lambda_0+\Lambda_m)$ if and
only if
$$
a(\lambda)\ge \ell(m(\mu))-2\;\;\text{and}\;\;a(m(\lambda))\ge \ell(\mu)-1.
$$
\item[(iii)]
If $m=2$ then $\lambda\otimes\mu\in B(\Lambda_0+\Lambda_m)$ if and
only if
$$
a(\lambda)\ge \ell(m(\mu))-1\;\;\text{and}\;\;a(m(\lambda))\ge \ell(\mu)-2.
$$
\end{itemize}
\end{prop}

Recall that $\ell(\roof(\mu))=\ell(\mu)$ by Lemma \ref{up operation}(3),
and $a(\base(\lambda))=a(\lambda)$ by Lemma \ref{down operation}(3).
By Proposition \ref{ceil and floor under Mullineux},
we have the following equalities.
\begin{itemize}
\item[(i)]
$a(\lambda)=a(\base(\lambda))$ and
$\ell(m(\mu))=\ell(\roof(m(\mu)))=a(\roof(\mu))$.
\item[(ii)]
$a(m(\lambda))=a(\base(m(\lambda)))=\ell(\base(\lambda))$ and
$\ell(\mu)=\ell(\roof(\mu))$.
\end{itemize}

Thus, his condition is precisely
$$
a(\tau_m(\base(\lambda)))\ge a(\roof(\mu))\;\;\text{and}\;\;
\ell(\tau_m(\base(\lambda)))\ge \ell(\roof(\mu)).
$$

Note that any $3$-core $\lambda$ is of the form
$(c,c-2,\dots,c-2r+2,d^2,(d-1)^2,\dots,1^2)$, where
$d=c-2r$ or $d=c-2r+1$.\footnote{The number of $i$ such that
$\lambda_i=\lambda_{i+1}+2$ is $r$ in the former case, and
$r-1$ in the latter case.}
In particular, $\lambda$ is
determined by $a(\lambda)$ and $\ell(\lambda)$, because
$a(\lambda)=c$ and $\ell(\lambda)=r+2d=2c-3r$ or $2c-3r+2$
imply
$$
r=-\bigl[\frac{\ell(\lambda)-2a(\lambda)}{3}\bigr],\;\;
d=\bigl[\frac{2\ell(\lambda)-a(\lambda)}{3}\bigr].
$$
Hence, the above condition is equivalent to
$\tau_m(\base(\lambda))\supset\roof(\mu)$.

As a conclusion, we may deduce Proposition \ref{Fayers} from our results,
and conversely, we may restate our results Theorem \ref{main result}
and Corollary \ref{main result2} in $e=3$ by using his more explicit
numerical conditions, which we do not mention here.

\end{document}